\documentclass[10pt,reqno]{amsart}
\usepackage{amssymb,version,graphicx,subfigure,fancybox,mathrsfs,pifont,booktabs,wrapfig,color}
\usepackage{url,hyperref}
\usepackage{subfigure}
\usepackage{subeqnarray}
\usepackage[ruled,vlined]{algorithm2e}
\usepackage{leftidx}
\usepackage{bm,amsfonts}
\usepackage[version=4]{mhchem}
\usepackage{mathrsfs}
\usepackage{amsmath}
\usepackage{amsfonts}
\usepackage{multirow}
\usepackage{arydshln} 
\usepackage{booktabs}
\usepackage{graphicx}
\usepackage{flafter}
\usepackage[section]{placeins}

\usepackage{booktabs}
\usepackage{multirow}
\usepackage{colortbl}
\definecolor{tabcolor}{rgb}{.105,.410,.113}

\allowdisplaybreaks
\catcode`\@=11
\@addtoreset{equation}{section}   
\renewcommand{\theequation}{\arabic{section}.\arabic{equation}}
\@addtoreset{table}{section}   
\renewcommand\thefigure{\thesection.\@arabic\c@figure}
\@addtoreset{figure}{section}   
\renewcommand\thetable{\thesection.\@arabic\c@table}

\newcommand{\new}{\newcommand*}
\new{\rnew}{\renewcommand*}
\new{\newe}{\newenvironment*}
\new{\stl}{\setlength}
\stl{\arraycolsep}{0.5mm}

\allowdisplaybreaks[4]

\newtheorem{thm}{\bf Theorem}

\newenvironment{theorem}{\begin{thm}} {\end{thm}}
\newtheorem{lmm}{\bf Lemma}

\newenvironment{lemma}{\begin{lmm}}{\end{lmm}}
\newtheorem{rem}{\bf Remark}

\newenvironment{remark}{\begin{rem}}{\end{rem}}
\newtheorem{defff}{\bf Definition}

\theoremstyle{definition}
\newtheorem{exm}{\bf Example} 
\newtheorem{pro}{\bf Proposition} 
\newenvironment{proposition}{\begin{pro}} {\end{pro}}
\newtheorem{cor}{\bf Corollary} 

\renewcommand{\theequation}{\arabic{section}.\arabic{equation}}
\newcommand{\refe}[1]{{\rm (\ref{#1})}}
\newcommand \dint {\displaystyle\int}

\renewcommand \dfrac {\displaystyle\frac}

\def\ball{\mathbb{B}}
\def\d{\mathrm{d}}
\def\BB^d{\mathbb{B}^d}

\title[M\"{u}ntz ball polynomials]{M\"{u}ntz ball polynomials and M\"{u}ntz spectral-Galerkin methods for singular eigenvalue problems}
\author[X. Yang, L. Wang, H. Li \& C. Sheng]{
		\;\; Xiu Yang${}^{\dag,\ddag}$, \;\;  Li-Lian Wang${}^{*}$, \;\; Huiyuan Li${}^{\S}$ \; and\; Changtao Sheng${}^{||}$
		}
	\thanks{${}^{\dag}$School of Mathematics and Statistics, Shandong University, Weihai 264209, Shandong,  China. The work of this author is partially supported by the Natural Science Foundation of Shandong Province  (No. ZR2021QA023), the National Natural Science Foundation of China (No. 12171284). Email: yangxiu0204@sdu.edu.cn(X. Yang). \\
	\indent${}^{\ddag}$School of Mathematics, Shandong University, Jinan 250100, Shandong, China. \\
	\indent ${}^{*}$Corresponding author. Division of Mathematical Sciences, School of Physical and Mathematical Sciences, Nanyang Technological University, 637371, Singapore. The research of this author is partially supported by Singapore MOE AcRF Tier 1 Grant: MOE-Tier1-RG15/21.  Email: lilian@ntu.edu.sg (L. Wang).\\
	\indent ${}^{\S}$ State Key Laboratory of Computer Science/Laboratory of Parallel Computing, Institute of Software, Chinese Academy of Sciences, 100190, Beijing, China. The research of this authors is is partially supported  by the National Natural Science Foundation of China (Nos. 12131005, 11871145 and 11971016). Email: huiyuan@iscas.ac.cn (H. Li). \\
	\indent ${}^{||}$School of Mathematics, Shanghai University of Finance and Economics, 200433, Shanghai, China. The work of this author is partially supported by the National Natural Science Foundation of China (Nos. 12201385 and 12271365), Shanghai Pujiang Program 21PJ1403500, the Fundamental Research Funds for the Central Universities 2021110474 and Shanghai Post-doctoral Excellence Program  2021154. Email: ctsheng@sufe.edu.cn (C. Sheng). \\
	\indent The first author would like to acknowledge the support of China Scholarship Council (CSC, No.  201906220068) for the visit of NTU to work on this topic, and thank Beijing Computational Science Research Center for hosting the visit where this
	research topic was  further explored.}

\subjclass{ 33C45, 33C55,  65N25,  47A75, 35Q40}
\keywords{M\"{u}ntz ball polynomials, M\"{u}ntz spectral-Galerkin methods, singular eigenvalue problems}

\begin{document}
	
\begin{abstract}
In this paper, we introduce a new family of orthogonal systems, termed as the M\"{u}ntz ball polynomials (MBPs), which are orthogonal with respect to the weight function: $\|\bm{x}\|^{2\theta+2\mu-2} (1-\|\bm{x}\|^{2\theta})^{\alpha}$ with the parameters $\alpha>-1, \mu>- 1/2$ and $\theta>0$ in 
 the $d$-dimensional unit ball $\bm x\in {\mathbb B}^d=\big\{\bm{x}\in\mathbb{R}^d: r=\|\bm{x}\|\leq1\big\}.$  We then develop efficient and spectrally accurate MBP spectral-Galerkin methods for singular eigenvalue problems including degenerating elliptic problems with perturbed ellipticity and Schr\"odinger's operators with fractional potentials. 
We demonstrate that the use of such non-standard basis functions can not only tailor to the singularity of the solutions but also  lead to sparse linear systems which can  be solved efficiently. 
\end{abstract}	
	
\maketitle


\vspace*{-10pt}
\section{Introduction}\label{section1}
Over  the past  decade, there has been a growing research interest in the construction of orthogonal polynomials and nonstandard basis functions  in spherical or related geometries with many interesting applications  (see, e.g., \cite{Dai2013approximation,Dunkl2014Orthogonal,Li2014spectral,Dyda2017Eigenvalues,Dyda2017Fractional,Atkinson2019Book,Olver2020Orthogonal,Zhang2020ACHA} and the references therein).  In particular, the multivariate polynomials/functions built upon Jacobi polynomials and spherical harmonics, are appealing basis functions for developing efficient spectral methods.
Sheng et al.\!    \cite{Sheng2021Nontensorial}   introduced nontensorial generalised Hermite polynomials/functions in arbitrary dimensions and developed efficient and accurate spectral method for solving PDEs in $\mathbb{R}^d$. Dai and Xu  provided in the book \cite{Dai2013approximation}  a cohesive account of the approximation theory and harmonic analysis on spheres and balls.  In particular,  the ball polynomials discussed therein
\begin{equation}\label{ballpoly20}
	P_{k, \ell}^{\alpha, n}(\boldsymbol{x})=P_{k}^{(\alpha, n+\frac{d}{2}-1)}(2\|\boldsymbol{x}\|^{2}-1)\, Y_{\ell}^{n}(\bm{\hat{x}}),\;\;\;\; \bm x \in \mathbb B^d,
\end{equation}
for $\alpha>-1$, $0\le \ell\le n $ and $k\ge 0,$  defined  in the $d$-dimensional unit ball  $\mathbb B^d, $ are orthogonal with respect to the weight function $\omega_{\alpha}(\boldsymbol{x}):=(1-\|\boldsymbol{x}\|^{2})^{\alpha}.$ Here,   $P_{k}^{(\alpha, \beta)}$ is the classical Jacobi polynomial and $Y_\ell^n$ is the spherical harmonic function.  The spectral algorithms using ball polynomials for  PDEs in balls have advantages over the usual spectral method based on mixed Legendre and spherical harmonic approximation (cf.\!  Li and Xu \cite{Li2014spectral}).  Dyda et al. \cite{Dyda2017Fractional}  proved constructively that the ball polynomials are eigenfunctions of the ``weighted''  fractional Laplacian operator
\begin{equation}\label{ballpoly2}
	(-\Delta)^{\alpha}\big((1-\|\boldsymbol{x}\|^{2})^{\alpha}_+\, P_{k, \ell}^{\alpha, n}(\boldsymbol{x})\big)= C_{d,k}^{\alpha, n}\, P_{k, \ell}^{\alpha, n}(\boldsymbol{x}),\quad \bm x\in {\mathbb B}^d,
\end{equation}
for $\alpha\in (0,1),$ where $a_{+}=\max\{a,0\}$ and  $C_{d,k}^{\alpha, n}$ is an explicit constant  (cf.\! \cite[(6)]{Dyda2017Fractional}). Dyda et al. \cite{Dyda2017Eigenvalues} further proposed  an efficient numerical scheme  based on this construction to study the fractional Laplacian eigenvalue problems in the unit ball with global homogeneous Dirichlet boundary condition.  Recently, Olver and Xu \cite{Olver2020Orthogonal} constructed orthogonal polynomials  on some quadratic surfaces of revolution such as cone, paraboloid and hyperboloid of revolution.  Consider for example the double cone  in $\mathbb R^{d+1}$:
\[\partial {\mathcal C}:=\big\{(\bm{x}, t)\,:\,  \|\bm x\|=|t|, \;\; \bm{x} \in \mathbb{R}^{d}, \,\, t\in [-1,1]\big\},\]
which has the apex at the origin and $t$-axis as its revolution  axis.  According to  \cite{Olver2020Orthogonal}, the orthogonal polynomials on the surface $\partial {\mathcal C}$ are given by
\begin{equation}\label{coneorth}
	S_{k, \ell}^{\mu, n}(\boldsymbol{x}, t)=C_{k-n}^{(\mu, n+\frac{d-1}{2})}(t)\, Y_{\ell}^{n}(\bm{\hat{x}}),
\end{equation}
where $C_{n}^{(\mu, \nu)}(t)$ is the generalized Gegenbauer polynomial of degree $n$.

In this paper, we introduce a new family of nonstandard basis functions, termed as M\"{u}ntz ball polynomials (MBPs), and also develop efficient and spectrally accurate M\"{u}ntz spectral-Galerkin methods for a class of singular eigenvalue problems.  Different from \eqref{ballpoly20} and \eqref{coneorth}, the MBPs $\big\{\mathcal{S}_{k,\ell,n}^{\alpha,\mu,\theta}(\{\beta_n\};\bm{x})\big\}$ (see \eqref{MBPs}),
 are composed of   M\"{u}ntz polynomials in the radial direction and  
orthogonal with respect to the weight function: $\|\bm{x}\|^{2\theta+2\mu-2} (1-\|\bm{x}\|^{2\theta})^{\alpha},\; \bm x\in {\mathbb B}^d$ with $\alpha>-1, \mu>- 1/2$ and $\theta>0$. 
 In addition to the parameters $\mu,\theta,$  they are also equipped with a general M\"untz  sequence  $\{\beta_n\}_{n=0}^ \infty$ that are free to choose. These provide great flexibility that can tailor to singularities or other properties of the underlying solutions. It is important to remark that 
 the $d$-dimensional MBPs with $\mu=0$, $\theta=1$ and $\beta_n=n+d/2-1$ reduce  to the ball polynomials in \cite{Dai2013approximation}. Moreover,  when $d=1$, $\theta=1$ and $\beta_n=n+d/2-1+\mu$, the corresponding orthogonal polynomials known as the generalised ultraspherical polynomials  which were first investigated by Chihara \cite[p.\!\! 156]{Chihara1978an}.
 
In this work, we focus on  the MBPs with
$$
\beta_n:=\beta_{n,c}^{\mu,\theta}(d)=\sqrt{c+(n+d/2-1)^2+\mu(\mu+d-2)}/\theta,
$$
and show that they are  the basis of choice for a class of  degenerate eigenvalue problems with singular potentials and  Schr\"{o}dinger eigenvalue problems with fractional power potentials (see  Table \ref{severalbasis} below). For such singular problems, it is challenging to construct spectrally accurate methods, though it is desirable.   

\begin{table}[!htb]
		\renewcommand*{\arraystretch}{1.8}
		\caption{\small MBP basis functions for eigenvalue problems with different operators} \label{severalbasis}
		\vspace*{-6pt} {\small
	\begin{tabular}{|ccccc|}
		\hline	
		\multicolumn{1}{|c|}{Operators} & \multicolumn{1}{c|}{Basis functions} & \multicolumn{1}{c|}{$\theta$} & \multicolumn{1}{c|}{$\mu$} \\ \hline
		\multicolumn{4}{|c|}{$\beta_{n}={\theta}^{-1}{\sqrt{c+(n+\frac{d}{2}-1)^2+\mu(\mu+d-2)}}$}                                                                 \\ \hline
		\multicolumn{1}{|c|}{$-\Delta +\dfrac{c}{\|\bm{x}\|^2}$ 
}
& \multicolumn{1}{c|}{${P}_{k}^{(-1,\beta_{n})}(2\|\bm{x}\|^2-1)\|\bm{x}\|^{\theta\beta_{n}+1-\frac{d}{2}} Y_\ell^n(\bm{\hat{x}})$}      & \multicolumn{1}{c|}{1} &  \multicolumn{1}{c|}{0}   \\ \hline
		\multicolumn{1}{|c|}{$-\Delta +\dfrac{c}{\|\bm{x}\|^2}+\dfrac{z}{\|\bm{x}\|}$ 
}        & \multicolumn{1}{c|}{${P}_{k}^{(-1,\beta_{n})}(2\|\bm{x}\|-1)\|\bm{x}\|^{\theta\beta_{n}+1-\frac{d}{2}} Y_\ell^n(\bm{\hat{x}})$}      & \multicolumn{1}{c|}{$\dfrac{1}{2}$} &  \multicolumn{1}{c|}{0}     \\ \hline
		\multicolumn{1}{|c|}{$-\Delta +\dfrac{c}{\|\bm{x}\|^2}+z\|\bm{x}\|^{\frac{2\nu-2\eta}{\eta+1}}$       } & \multicolumn{1}{c|}{${P}_{k}^{(-1,\beta_{n})}(2\|\bm{x}\|^{\frac{2}{\eta+1}}-1)\|\bm{x}\|^{\theta\beta_{n}+1-\frac{d}{2}} Y_\ell^n(\bm{\hat{x}})$}      & \multicolumn{1}{c|}{$\dfrac{1}{\eta+1}$} &  \multicolumn{1}{c|}{0}  \\ \hline		
		\multicolumn{1}{|c|}{$-\nabla(\|\bm{x}\|^{2\mu}\nabla)  + c\|\bm{x}\|^{2\mu-2} $ }       & \multicolumn{1}{c|}{$P_{k}^{(-1,\beta_{n})} \left(2 \|\bm{x}\|^{2-2\mu}-1\right) \|\bm{x}\|^{\theta\beta_{n}+1-\frac{d}{2}-\mu} Y_\ell^n(\bm{\hat{x}})$}      & \multicolumn{1}{c|}{$1-\mu$} & \multicolumn{1}{c|}{$> -\frac{1}{2}$}	\\ \hline
	\end{tabular}}
\end{table}

The paper is organised as follows. In Section \ref{sect2:preli}, we make necessary preparations by reviewing some properties of Jacobi polynomials and spherical harmonics. 
In Section \ref{sect3:MBPs}, we define the MBPs and present their important properties that pave the way for developing the spectral algorithms.  In Section 4, we construct efficient and accurate MBP spectral-Galerkin methods for some singular eigenvalue problems.  The final section is some concluding remarks. 

\section{Preliminaries on Jacobi polynomials and spherical harmonics functions}\label{sect2:preli}

In this section, we introduce some notation and review some relevant properties of Jacobi polynomials and spherical harmonics that  are necessary for the definition of the new spectral basis functions.
\subsection{Jacobi and  generalised Jacobi polynomials}\label{subsection:pre}

 For $\alpha\in\mathbb{R}$, the rising factorial in the Pochhammer symbol, binomial factor and Gamma function are related by
\[(\alpha)_k=\alpha(\alpha+1)(\alpha+2)\cdots(\alpha+k-1)=\frac{\Gamma(\alpha+k)}{\Gamma(\alpha)}.   \]

We follow the definition and normalization in Szeg\"{o} \cite{Szego1975orthogonal}.  For $\alpha,\beta>-1,$  the classical Jacobi polynomials are mutually orthogonal with respect to the weight function $\omega^{(\alpha, \beta)}(x)=(1-x)^{\alpha}(1+x)^{\beta}$ on $I=(-1,1):$
\begin{equation}\label{GGPeq9}
	\int_{-1}^1P_n^{(\alpha,\beta)}(x)P_m^{(\alpha,\beta)}(x)\omega^{(\alpha, \beta)}(x)\,\mathrm{d}x=\frac{2^{\alpha+\beta+1}\Gamma(n+\alpha+1)\Gamma(n+\beta+1)}{(2n+\alpha+\beta+1)\Gamma(n+1)\Gamma(n+\alpha+\beta+1)}\delta_{nm},
\end{equation}
which can be defined in terms of the hypergeometric function
\begin{equation}\label{GGPeq41}
	\begin{aligned}
		P_{n}^{(\alpha, \beta)}(x)
		&={{n+\alpha} \choose n}\,_{2} F_{1}\Big(\!-n, n+\alpha+\beta+1 ; \alpha+1 ;  \frac{1-x}{2}\Big)\\
		&=\sum_{k=0}^{n} \frac{(n+\alpha+\beta+1)_{k}(\alpha+k+1)_{n-k}}{k !(n-k) !}\Big(\frac{x-1}{2}\Big)^{k}.
	\end{aligned}
\end{equation}
They are the eigenfunctions of the Sturm-Liouville problem
\begin{equation}\label{JacobiSL}
	\mathcal{L}_{x}^{(\alpha, \beta)} P_{n}^{(\alpha, \beta)}(x):=-\frac{1}{\omega^{(\alpha, \beta)}(x)} \partial_{x}\big(\omega^{(\alpha+1, \beta+1)}(x) \partial_{x} P_{n}^{(\alpha, \beta)}(x)\big)=\lambda_{n}^{(\alpha, \beta)} P_{n}^{(\alpha, \beta)}(x), \quad x \in I,
\end{equation}
with the corresponding eigenvalues $\lambda_{n}^{(\alpha, \beta)}=n(n+\alpha+\beta+1)$. Hereafter, we sometimes use $\partial_x$ to denote the ordinary derivative $\frac{\rm d}{\mathrm{d}x}$.

By \cite[p.\!\! 304]{Andrews1999special}, the Jacobi polynomials  satisfy  the following properties
\begin{align}
	&\partial_x P_n^{(\alpha,\beta)}(x)=\frac{n+\alpha+\beta+1}{2}P_{n-1}^{(\alpha+1,\beta+1)}(x),\label{LemEQ1}\\
	&(2n+\alpha+\beta+1)P_n^{(\alpha,\beta)}(x)=(n+\alpha+\beta+1)P_n^{(\alpha+1,\beta)}(x)-(n+\beta)P_{n-1}^{(\alpha+1,\beta)}(x),\label{LemEQ2}\\
	& \Big(n+\frac{\alpha+\beta}{2}+1\Big)(1-x)P_n^{(\alpha+1,\beta)}(x)= (n+\alpha+1)P_n^{(\alpha,\beta)}(x) - (n+1)P_{n+1}^{(\alpha,\beta)}(x)\label{LemEQ3}.
\end{align}

 The Jacobi polynomials   can be generalized to  cases with general $\alpha,\beta\in {\mathbb R}$ as in \cite{Szego1975orthogonal,Li2010Optimal,Cagliero2015Explicit}.
Assume that $P_{n}^{(\alpha, \beta)}(x)$ does not vanish identically in $x$.
\begin{itemize}
	\item[(i)] $P_{n}^{(\alpha, \beta)}(x)$ has degree $<n$ in $x$ if and only if (iff)
	$-n-\alpha-\beta\in\{1,2,\ldots,n\}$. Then the degree is $-n-\alpha-\beta-1$ and
	\begin{equation}\label{GJP1}
		P_{n}^{(\alpha, \beta)}(x)=\frac{\Gamma(n+\alpha+1)\Gamma(-n-\alpha-\beta)}{\Gamma(n+1)\Gamma(-n-\beta)} P_{-n-\alpha-\beta-1}^{(\alpha, \beta)}(x) .
	\end{equation}
	\item[(ii)]$P_{n}^{(\alpha, \beta)}(1)=0$ iff
	$-\alpha\in\{1,2,\ldots,n\}$. Then the zero at 1 has multiplicity $-\alpha$ and
	\begin{equation}\label{GJP2}
		P_{n}^{(\alpha, \beta)}(x)=\frac{\Gamma(n+\alpha+1)\Gamma(n+\beta+1)}{\Gamma(n+1)\Gamma(n+\alpha+\beta+1)}\Big(\frac{x-1}{2}\Big)^{-\alpha} P_{n+\alpha}^{(-\alpha, \beta)}(x).
	\end{equation}
   	\item[(iii)]$P_{n}^{(\alpha, \beta)}(-1)=0$ iff
	$-\beta\in\{1,2,\ldots,n\}$. Then the zero at $-1$ has multiplicity $-\beta$ and
	\begin{equation}\label{GJP3}
		P_{n}^{(\alpha, \beta)}(x)=\frac{\Gamma(n+\alpha+1)\Gamma(n+\beta+1)}{\Gamma(n+1)\Gamma(n+\alpha+\beta+1)}\Big(\frac{x+1}{2}\Big)^{-\beta} P_{n+\beta}^{(\alpha, -\beta)}(x).
	\end{equation}
	\item[(iv)]$P_{n}^{(\alpha, \beta)}(\pm 1)=0$  iff
$-\alpha, -\beta$ or $-\alpha-\beta \in\{1,2,\ldots,n\}$. Then the zero at 1 has multiplicity $-\alpha$ and at $-1$ has multiplicity $-\beta$
	\begin{equation}
		P_{n}^{(\alpha, \beta)}(x)= \Big(\frac{x-1}{2}\Big)^{-\alpha}\Big(\frac{x+1}{2}\Big)^{-\beta}P_{n+\alpha+\beta}^{(-\alpha, -\beta)}(x).
	\end{equation}
\end{itemize}
In particular, if $\alpha=-1$, we directly obtain from \eqref{GJP2} that
\begin{equation}\label{GGPeqJapobi}
	P_{0}^{(-1, \beta)}(x)=1, \quad P_{n}^{(-1, \beta)}(x)=\frac{n+\beta}{n} \frac{x-1}{2} P_{n-1}^{(1, \beta)}(x), \quad n \geq 1,\,\, \beta>-1 .
\end{equation}

\subsection{Spherical harmonics} Here, we adopt the same definition and notation as in  Dai and Xu \cite{Dai2013approximation}. 	Let $\mathbb N$  be the set of positive integers and ${\mathbb N}_0:=\{0\}\cup {\mathbb N}$. Let $\mathbb{R}^d$ be the $d$-dimensional Euclidean space. For any $\bm{x},\bm{y}\in\mathbb{R}^d$, we define the inner product and norm of $\mathbb{R}^d$ as $\langle\bm{x},\bm{y}\rangle:=\sum_{i=1}^dx_iy_i$, and  $\|\bm{x}\|:=\sqrt{\langle\bm{x},\bm{x}\rangle}$, respectively. Denote the unit vector along $\bm{x}\not=\bm 0$ by $\bm{\hat{x}}=\bm{x}/\|\bm{x}\|$. The unit sphere $\mathbb{S}^{d-1}$ and the unit ball $\mathbb{B}^d$ of $\mathbb{R}^d$ are respectively defined by \[\mathbb{S}^{d-1}:=\big\{\bm{\hat{x}}\in\mathbb{R}^d:\|\bm{\hat{x}}\|=1\big\},\quad \mathbb{B}^d:=\big\{\bm{x}\in\mathbb{R}^d:\|\bm{x}\|\leq1\big\}.\]
Define the inner product
\begin{equation}\label{sphhar}
	\langle f,g\rangle_{\mathbb{S}^{d-1}}:=\int_{\mathbb{S}^{d-1}}f(\bm{\hat{x}})g(\bm{\hat{x}})\, {\rm d}\sigma(\bm{\hat{x}}),\end{equation}
where ${\rm d}\sigma$ is the surface measure, under which the spherical harmonics of different degree are orthogonal to each other (cf. \cite[Thm. 1.1.2]{Dai2013approximation}).

Let $\mathcal{P}_n^d$ be the space of multivariate homogeneous polynomials of degree $n$:
\[\mathcal{P}_n^d=\mathrm{span}\big\{\bm{x^k}=x_1^{k_1}x_2^{k_2}\cdots x_d^{k_d}:k_1+k_2+\cdots+k_d=n,\, k_i\in{\mathbb N}_0, i=1, \ldots, d\big\}.\]
As a subspace of $\mathcal{P}_{n}^{d},$  the  space of all real harmonic polynomials of degree $n$ is defined as
\[\mathcal{H}_n^d:=\big\{P\in\mathcal{P}_n^d:\Delta P(\bm{x})=0\big\}.\]
It is known that the dimensionality are
\begin{equation}\label{dim}
	\mathrm{dim}\big(\mathcal{P}_n^d\big)={{n+d-1} \choose n}; \quad\;\; a_n^d:=\mathrm{dim} \big(\mathcal{H}_n^d\big)={{n+d-1} \choose n}-{{n+d-3}\choose {n-2}},
\end{equation}
where for $n=0,1$, the second binomial coefficient in $ a_n^d$ vanishes (cf. \cite[(1.1.5)]{Dai2013approximation}).

Note that (cf. \cite[Ch. 1]{Dai2013approximation}): for any $Y(\bm{x})\in \mathcal{H}_n^d$,
\begin{equation}\label{Ybm}
	Y(\bm{x})=\|\bm{x}\|^nY(\bm{x}/\|\bm{x}\|)=\|\bm x\|^nY(\bm{\hat{x}}),\quad \bm{\hat{x}}\in \mathbb{S}^{d-1}.
\end{equation}
Define the spherical component of the gradient operator in spherical-polar coordinates $(r,\bm{\hat{x}})$,
\begin{equation}\label{nabla0}
	\nabla_{\!0}=r(\nabla-\bm{\hat{x}}\partial_r), \quad \bm{x}=r\bm{\hat{x}},\,\, r=\|\bm{x}\|, \,\,\bm{\hat{x}}\in\mathbb{S}^{d-1}.
\end{equation}
Then the Laplace operator takes the form
\begin{equation}\label{GGPeq37}
	\Delta = \partial^2_r+ \frac{d-1}{r} \partial_r + \frac{1}{r^2} \Delta_{0},
\end{equation}
where $\Delta_{0}:=\nabla_{\!0}\cdot\nabla_{\!0}$ is the Laplace-Beltrami operator on $\mathbb S^{d-1}$.
Consequently,
$$
\Delta Y(\boldsymbol{x})=\Delta\big[r^n Y(\hat{\boldsymbol{x}})\big]=n(n+d-2) r^{n-2} Y(\hat{\boldsymbol{x}})+r^{n-2} \Delta_0 Y(\hat{\boldsymbol{x}}).
$$
The spherical harmonics are eigenfunctions of the Laplace-Beltrami operator,
\begin{equation}\label{GGPeq38}
\Delta_0 Y(\hat{\boldsymbol{x}})=-n(n+d-2) Y(\hat{\boldsymbol{x}}), \quad Y \in \mathcal{H}_n^d, \quad \hat{\boldsymbol{x}} \in \mathbb{S}^{d-1} .
\end{equation}
For $n\in \mathbb{N}_0$, let $\{Y_\ell^n(\bm {\hat x}): 1 \le \ell \le a_n^d\}$ be the (real) normalised spherical harmonic basis  of  $\mathcal{H}_n^d$ on $\mathbb{S}^{d-1}$,   so that
\begin{equation}\label{GGPeq8}
	\langle Y_\ell^n,Y_\iota^m\rangle_{\mathbb{S}^{d-1}}=\delta_{nm}\delta_{\ell\iota}, \quad (\ell,n), (\iota,m)\in \Upsilon_{\infty}^d.
\end{equation}
Here, we introduce the index sets:  
\[\begin{array}{l}
	\Upsilon_{\infty}^d=\{(\ell,n):1\leq\ell\leq a_n^d, 0\leq n < \infty, \ell,n\in\mathbb{N}_0\},\vspace{1ex}\\
	\Upsilon_{N}^d=\{(\ell,n):1\leq\ell\leq a_n^d, 0\leq n\leq N, \ell,n\in\mathbb{N}_0\}.
\end{array}\]
Finally, let
$L_{\omega}^2(\mathbb{B}^d)$ be  the space of square integrable  functions weighted with a generic non-negative weight function $\omega(\bm{x})$ in the ball  $\mathbb{B}^d$ furnished with the inner product and norm
\begin{equation}\label{innerprod}
(u, v)_{\omega}:=\int_{\mathbb{B}^d} u(\boldsymbol{x}) v(\boldsymbol{x}) \omega(\boldsymbol{x}) \mathrm{d} \boldsymbol{x} ,\quad \|u\|_{\omega}=\sqrt{(u,u)_{\omega}}\,.
\end{equation}
When $\omega(\bm{x})=1$, we drop the weight in the above notation.

\section{M\"untz ball polynomials}\label{sect3:MBPs}
\setcounter{equation}{0}
\setcounter{lmm}{0}
\setcounter{thm}{0}

In dynamical systems, there are several system with singular potential of the type $1/\|\bm x\|^{2\mu}$, e.g., the Coulomb's potential in the Schr\"{o}dinger equation. There are many difficulties in numerical simulations for those singular terms. In order to capture the singularity and construct efficient numerical schemes,  we introduce the M\"{u}ntz ball polynomials and present various appealing properties of this new family of basis functions.

\subsection{Definition of M\"{u}ntz ball polynomials}
We first review the definition of the M\"untz polynomials in one dimension.
Given an increasing sequence of distinct non-negative real numbers
$$\Lambda:=\big\{ 0 \leq \lambda_0<\lambda_1<\cdots <\lambda_k<\cdots \to \infty \big\},$$
we define the  space of M\"untz polynomials
$$
M(\Lambda)=\operatorname{span}\left\{x^{\lambda_0}, x^{\lambda_1}, x^{\lambda_2}, \ldots\right\},\quad x\ge 0.
$$
As a generalisation of the  Weierstrass approximation theorem,  the  M\"untz Theorem \cite{Cheney1998Introduction} states
that any continuous function on a closed and bounded interval can be uniformly approximated by  M\"untz polynomials,  if and only if
$\lambda_0=0$ and $\sum_{k=1}^{\infty}\lambda_k^{-1}=\infty.$

We next introduce the new basis function that are  warped products of the  M\"{u}ntz polynomials in radial direction and spherical harmonics in the angular directions.
This notation also generalizes the  ball polynomials (cf.\! \cite{Dai2013approximation}). For this reason, we term them as M\"{u}ntz ball polynomials.
\begin{defff}[{\bf M\"{u}ntz ball polynomials}] \label{MBPs-defn} Given  real $\alpha>-1 $,  $\mu>-1/2$, $\theta>0$, and given an increasing sequence $\{\beta_n\}_{n=0}^ \infty$ with  $\beta_0>-1$, the M\"{u}ntz ball polynomials in a $d$-dimensional unit ball $\ball^d$ are defined as
\begin{equation}\label{MBPs}	\mathcal{S}_{k,\ell,n}^{\alpha,\mu,\theta}(\bm{x};\{\beta_n\})={P}_{k}^{(\alpha,\beta_n)}(2\|\bm{x}\|^{2\theta}-1)\|\bm{x}\|^{\theta\beta_n+1-\frac{d}{2}-\mu} Y_\ell^n(\bm{\hat{x}}),\;\;\; \bm{x}\in \ball^d,
\end{equation}
for $k\in\mathbb{N}_0 $ and $(\ell,n)\in \Upsilon_{\infty}^d$.
Moreover, the MBPs are also defined for $\alpha\le -1$ with the understanding of the Jacobi polynomials in radial direction as the generalised ones in Subsection \ref{subsection:pre}.
\end{defff}

When $\mu=0$, $\theta=1$ and $\beta_n=n+d/2-1$, the MBPs reduce to the ball polynomials (cf.\! \cite{Dai2013approximation}):
for $\alpha>-1$,
\begin{equation}\label{ballpoly}
P_{k, \ell}^{\alpha, n}(\boldsymbol{x})=P_k^{(\alpha, n+\frac{d}{2}-1)}(2\|\boldsymbol{x}\|^2-1) Y_{\ell}^n(\boldsymbol{x}), \quad \boldsymbol{x} \in \mathbb{B}^d, \;\; k\in\mathbb{N}_0,\;\; (\ell,n)\in \Upsilon_{\infty}^d.
\end{equation}

With the above selection of the parameters, the so-defined MBPs are orthogonal in the following sense.
\begin{theorem}The MBPs defined in \eqref{MBPs}  are mutually orthogonal with respect to the weight function $\|\bm{x}\|^{2\theta+2\mu-2} (1-\|\bm{x}\|^{2\theta})^{\alpha}$, namely,
\begin{align}
	&\int_{\ball^d}\mathcal{S}_{k,\ell,n}^{\alpha,\mu,\theta}(\bm{x};\{\beta_n\})\mathcal{S}_{j,\iota,m}^{\alpha,\mu,\theta}(\bm{x};\{\beta_n\}) \|\bm{x}\|^{2\theta+2\mu-2} (1-\|\bm{x}\|^{2\theta})^{\alpha}\d \bm{x} \nonumber\\
		& \qquad =\frac{\Gamma(k+\alpha+1)\Gamma(k+\beta_n+1)}{2\theta(2k+\alpha+\beta_n+1)\Gamma(k+1)\Gamma(k+\alpha+\beta_n+1) } \delta_{kj}\delta_{nm}  \delta_{\ell \iota}.\label{orthp2}
\end{align}
\end{theorem}
\begin{proof} Noting that
\begin{equation}\label{ballint}
\int_{\mathbb{B}^d} f(\bm{x}) \mathrm{d} \bm{x}=\int_0^1 r^{d-1}  \int_{\mathbb{S}^{d-1}} f(r \bm{\hat{x}}) \mathrm{d} \sigma(\bm{\hat{x}})\, \mathrm{d} r,
\end{equation}
we obtain  from \eqref{sphhar} that
\begin{align}
	&\int_{\ball^d}\mathcal{S}_{k,\ell,n}^{\alpha,\mu,\theta}(\bm{x};\{\beta_n\})\mathcal{S}_{j,\iota,m}^{\alpha,\mu,\theta}(\bm{x};\{\beta_n\})\|\bm{x}\|^{2\theta+2\mu-2} (1-\|\bm{x}\|^{2\theta})^{\alpha}\d \bm{x} \nonumber\\
	=&\int_{\mathbb{S}^{d-1}} Y_{\ell}^{n}(\bm{\hat{x}}) Y_{\iota}^{m}(\bm{\hat{x}}) \d \sigma(\bm{\hat{x}}) \int_{0}^1{P}_{k}^{(\alpha,\beta_n)}(2r^{2\theta}-1) {P}_{j}^{(\alpha,\beta_n)}(2r^{2\theta}-1) (1-r^{2\theta})^{\alpha} r^{2\theta\beta_n+2\theta-1}\d r\nonumber\\
	=\,&\delta_{nm} \delta_{\ell\iota} \int_{0}^1{P}_{k}^{(\alpha,\beta_n)}(2r^{2\theta}-1) {P}_{j}^{(\alpha,\beta_n)}(2r^{2\theta}-1) (1-r^{2\theta})^{\alpha} r^{2\theta\beta_n+2\theta-1}\d r \label{muntz}\\
	=\,&\delta_{nm} \delta_{\ell \iota} \times \frac{1}{4\theta}\int_{-1}^1{P}_{k}^{(\alpha,\beta_n)}(\rho) {P}_{j}^{(\alpha,\beta_n)}(\rho) \Big(\frac{1-\rho}{2}\Big)^{\alpha}\Big(\frac{1+\rho}{2}\Big)^{\beta_n}\d \rho \qquad (\text{note: } \rho=2r^{2\theta}-1) \nonumber\\
	=\,&\frac{\Gamma(k+\alpha+1)\Gamma(k+\beta_n+1)}{2\theta(2k+\alpha+\beta_n+1)\Gamma(k+1)\Gamma(k+\alpha+\beta_n+1) } \delta_{kj}\delta_{nm}  \delta_{\ell \iota}, \nonumber
\end{align}
where in the last step, we used the orthogonality \eqref{GGPeq9} of Jacobi polynomials.
\end{proof}

	The one-dimensional MBPs are of independent interest. For $d=1$, we derive from \refe{dim} that $a_0^1=a_1^1=1$ and $a_n^1=0$ for $n\geq 2$, so there exist only two orthonormal harmonic polynomials: $Y_1^0(x)={1}/{\sqrt{2}}$ and $Y_1^1(x)={x}/{\sqrt{2}}$. As a result, the one-dimensional MBPs with  $\beta_n={(n-1/2+\mu)}/{\theta}$ reduce to
		\begin{equation}\label{1dGGPs}
			\begin{split}
				&\mathcal{S}_{k,1,0}^{\alpha,\mu,\theta}(x;\beta_0)=\dfrac{1}{\sqrt{2}}P_k^{(\alpha,\frac{\mu-1/2}{\theta})}(2x^{2\theta}-1), \\  
				& \mathcal{S}_{k,1,1}^{\alpha,\mu,\theta}(x;\beta_1)=\dfrac{1}{\sqrt{2}}\,xP_k^{(\alpha,\frac{\mu+1/2}{\theta})}(2x^{2\theta}-1),
			\end{split}
		\end{equation}
 which are mutually orthogonal with respect to the weight function $|x|^{2\theta+2\mu-2}(1-x^{2\theta})^\alpha$. It is important to point out that these special MBPs with $\theta=1$ are closely related to the  generalised ultraspherical polynomials (GUPs) introduced in \cite[p. 156]{Chihara1978an}:
\begin{equation}\label{relationG2J}
	S_{2k+l}^{(\alpha,\mu)}(x)=x^lP_{k}^{(\alpha,\mu-1/2+l)}(2x^2-1),\quad  l=0,1.
\end{equation}
When $\theta=1$, one has
	\begin{equation}\label{1dGGPs1}
	\begin{split}
		&\mathcal{S}_{k,1,0}^{\alpha,\mu,1}(x)=\dfrac{1}{\sqrt{2}}P_k^{(\alpha,\mu-1/2)}(2x^{2}-1)=\dfrac{1}{\sqrt{2}}S_{2k}^{(\alpha,\mu)}(x), \\  
		& \mathcal{S}_{k,1,1}^{\alpha,\mu,1}(x)=\dfrac{1}{\sqrt{2}}\,xP_k^{(\alpha,\mu+1/2)}(2x^{2}-1)=\dfrac{1}{\sqrt{2}}S_{2k+1}^{(\alpha,\mu)}(x).
	\end{split}
\end{equation}
Thus these MBPs  are the GUPs with a different normalisation.

\subsection{The family of MBPs of interest} In what follows, we focus on the  MBPs \eqref{MBPs}  with the  specific parameters
\begin{equation}\label{betank}
\beta_n:=\beta_{n,c}^{\mu,\theta}(d)=\frac 1 {\theta}{\sqrt{c+\Big(n+\frac{d}{2}-1\Big)^2+\mu(\mu+d-2)}},\quad c>0.
\end{equation}
Accordingly,  we drop the dependence on $\beta_n$ and simply  denote
\begin{equation}\label{MBPsSk}
\mathcal{S}_{k,\ell,n,c}^{\alpha,\mu,\theta}(\bm{x})=\mathcal{S}_{k,\ell,n}^{\alpha,\mu,\theta}(\bm{x};\{\beta_n\}),\quad \bm x\in{\mathbb B}^d,
\end{equation}
where $\alpha\in {\mathbb R}$,  $\mu>-1/2$, $\theta>0,$   $k\in\mathbb{N}_0 $ and $(\ell,n)\in \Upsilon_{\infty}^d$ as before.

\begin{theorem}\label{theoSL} For $\alpha>-1$, $\mu>-1/2$ and $c,\theta>0$, we define the second-order differential operator:
\begin{align}
	\mathscr{D}_{c,\theta,\bm{x}}^{(\alpha,\mu)}&=-\|\bm{x}\|^{-2\mu}(1-\|\bm{x}\|^{2\theta})^{-\alpha} \nabla \cdot \big((\bm{I}_d-\|\bm{x}\|^{2\theta-2}\bm{x}\bm{x}^t)\|\bm{x}\|^{2\mu}(1-\|\bm{x}\|^{2\theta})^{\alpha}\big)\nabla+\frac{c}{\|\bm{x}\|^{2}}\\
	&=-\|\bm{x}\|^{-2\mu}(1-\|\bm{x}\|^{2\theta})^{-\alpha} \nabla \cdot (\|\bm{x}\|^{2\mu}(1-\|\bm{x}\|^{2\theta})^{\alpha+1})\nabla+\frac{c}{\|\bm{x}\|^{2}}-\|\bm{x}\|^{2\theta-2}\Delta_0.	
\end{align}
Then the MBPs $\{\mathcal{S}_{k,\ell,n,c}^{\alpha,\mu,\theta}(\bm{x})\}$ are the eigenfunctions of the operator
$\|\bm{x}\|^{2-2\theta}\mathscr{D}_{c,\theta,\bm{x}}^{(\alpha,\mu)},$ that is,
	\begin{equation}\label{eqSLMuntz} \|\bm{x}\|^{2-2\theta}\mathscr{D}_{c,\theta,\bm{x}}^{(\alpha,\mu)}\mathcal{S}_{k,\ell,n,c}^{\alpha,\mu,\theta}(\bm{x})=\chi_{\theta,k}^{\alpha,\mu}\,\mathcal{S}_{k,\ell,n,c}^{\alpha,\mu,\theta}(\bm{x}),
	\end{equation}
	where the eigenvalues are
\begin{equation}\label{eignvalueChi}
\chi_{\theta,k}^{\alpha,\mu}:=\chi_{\theta,k}^{\alpha,\mu}(\beta_n)=(2 \theta k+\theta \beta_n-\mu+1-d / 2)(2 \theta k+\theta \beta_n+2 \theta \alpha+2 \theta+\mu+d / 2-1).
\end{equation}
\end{theorem}
\begin{proof}
	Using the Leibniz rule for gradient and divergence, one derives the operator $\mathscr{D}_{c,\theta,\bm{x}}^{(\alpha,\mu)}$ has the  equivalent form:
	\[\begin{aligned}
		\mathscr{D}_{c,\theta,\bm{x}}^{(\alpha,\mu)}&=-\|\bm{x}\|^{-2\mu}(1-\|\bm{x}\|^{2\theta})^{-\alpha} \nabla \cdot (\|\bm{x}\|^{2\mu}(1-\|\bm{x}\|^{2\theta})^{\alpha+1})\nabla+\frac{c}{\|\bm{x}\|^{2}}-\|\bm{x}\|^{2\theta-2}\Delta_0 \\
		&=-\left(1-r^{2 \theta}\right) \partial_{r}^{2}-\frac{d-1+2 \mu}{r} \partial_{r}+[2 \theta(\alpha+1)+d-1+2 \mu] r^{2 \theta-1} \partial_{r}+\frac{c-\Delta_{0}}{r^{2}}\\
		&=-\frac{1}{r^{d+2\mu-1}(1-r^{2\theta})^\alpha}\partial_r\big[r^{d+2\mu-1}(1-r^{2\theta})^{\alpha+1}\partial_r\big]+\frac{c-\Delta_{0}}{r^{2}},
	\end{aligned}	\]
	where we used  \eqref{GGPeq37} and the identity: $\bm{x}\cdot\nabla=r\bm{\hat{x}}\cdot\nabla=r\partial_r$.
	
	In view of the definition of $\Delta_0$, one derives that in the $r$-direction,
	\begin{equation}\label{diffr}		\mathscr{D}_{c,\theta,r}^{(\alpha,\mu)}=-\frac{1}{r^{d+2\mu-1}(1-r^{2\theta})^\alpha}\partial_r\big[r^{d+2\mu-1}(1-r^{2\theta})^{\alpha+1}\partial_r\big]+\frac{c+n(n+d-2)}{r^{2}}.
	\end{equation}
	With a change of variable: $\rho=2r^{2\theta}-1$, we obtain that 
	\[\begin{aligned}
		&	\mathscr{D}_{c,\theta,r}^{(\alpha,\mu)}[r^{\theta\beta_n+1-\frac{d}{2}-\mu}P_k^{(\alpha,\beta_n)}(2r^{2\theta}-1)]\\		&=-\frac{r^{\theta\beta_n+2\theta-1-\frac{d}{2}-\mu}}{(1-r^{2\theta})^{\alpha}r^{2\theta\beta_n+2\theta-1}}\partial_r\big[(1-r^{2\theta})^{\alpha+1}r^{2\theta\beta_n+1}\partial_rP_k^{(\alpha,\beta_n)}(2r^2-1)\big]\\ &\quad+(\theta\beta_n+1-\frac{d}{2}-\mu)(2\theta\alpha+\theta\beta_n-1+\frac{d}{2}+\mu+2\theta)r^{\theta\beta_n+2\theta-1-\frac{d}{2}-\mu}P_k^{(\alpha,\beta_n)}(2r^2-1)\\
&=-4\theta^2r^{\theta\beta_n+2\theta-1-\frac{d}{2}-\mu}(1-\rho)^{-\alpha}(1+\rho)^{-\beta_n}\partial_\rho\big[(1-\rho)^{\alpha+1}(1+\rho)^{\beta_n+1}\partial_\rho P_k^{(\alpha,\beta_n)}(\rho)\big]\\		&\quad+(\theta\beta_n+1-\frac{d}{2}-\mu)(2\theta\alpha+\theta\beta_n-1+\frac{d}{2}+\mu+2\theta)r^{\theta\beta_n+2\theta-1-\frac{d}{2}-\mu}P_k^{(\alpha,\beta_n)}(2r^2-1).
	\end{aligned}\]
Thus using  \eqref{JacobiSL} leads to
	\[\begin{aligned}
		&\mathscr{D}_{c,\theta,r}^{(\alpha,\mu)}[r^{\theta\beta_n+1-\frac{d}{2}-\mu}P_k^{(\alpha,\beta_n)}(2r^2-1)]\\
		&=[4\theta^2k(k+\alpha+\beta_n+1)]r^{\theta\beta_n+2\theta-1-\frac{d}{2}-\mu}P_k^{(\alpha,\beta_n)}(2r^2-1)\\ &\quad+(\theta\beta_n+1-\frac{d}{2}-\mu)(2\theta\alpha+\theta\beta_n-1+\frac{d}{2}+\mu+2\theta)r^{\theta\beta_n+2\theta-1-\frac{d}{2}-\mu}P_k^{(\alpha,\beta_n)}(2r^2-1)\\
	    &=\chi_{\theta,k}^{\alpha,\mu}\, r^{\theta\beta_n+2\theta-1-\frac{d}{2}-\mu}P_k^{(\alpha,\beta_n)}(2r^2-1).
	\end{aligned}\]
 Noting the relation between $\mathscr{D}_{c,\theta,\bm{x}}^{(\alpha,\mu)}$ and $\mathscr{D}_{c,\theta,r}^{(\alpha,\mu)},$ we obtain
  \eqref{eqSLMuntz}  from \eqref{GGPeq38} and \eqref{MBPs} directly.
  \end{proof}

As a  consequence of Theorem \ref{theoSL}, we have the following property.
\begin{cor} \label{proorth}{\em Given the  parameters satisfying the same conditions as in Theorem \ref{theoSL},
	the MBPs $\{\mathcal{S}_{k,\ell,n,c}^{\alpha,\mu,\theta}(\bm{x})\}$ form a Sobolev orthogonal basis in the sense that
	\begin{align}
		&\big(\nabla \mathcal{S}_{k,\ell,n,c}^{\alpha,\mu,\theta} , \nabla \mathcal{S}_{j,\iota,m,c}^{\alpha,\mu,\theta} \big)_{r^{2\mu}(1-r^{2\theta})^{\alpha+1} }+c\big(\mathcal{S}_{k,\ell,n,c}^{\alpha,\mu,\theta} , \mathcal{S}_{j,\iota,m,c}^{\alpha,\mu,\theta}\big)_{r^{2\mu-2}(1-r^{2\theta})^{\alpha} } \nonumber\\
		&\quad -\big(\nabla_0 \mathcal{S}_{k,\ell,n,c}^{\alpha,\mu,\theta} , \nabla_0 \mathcal{S}_{j,\iota,m,c}^{\alpha,\mu,\theta} \big)_{r^{2\mu+2\theta-2}(1-r^{2\theta})^{\alpha} } \nonumber \\
		&=\frac{\chi_{\theta,k}^{\alpha,\mu} \ \Gamma(k+\alpha+1)\Gamma(k+\beta_n+1)}{2\theta(2k+\alpha+\beta_n+1)\Gamma(k+1)\Gamma(k+\alpha+\beta_n+1) } \delta_{kj}\delta_{nm}  \delta_{\ell \iota}.  \label{orthMBPs}
	\end{align}}
	\end{cor}
  \begin{proof}
 Using  \eqref{eqSLMuntz}, we derive from the orthogonality \eqref{orthp2} that
    	\begin{align}
    		&\big(\nabla \mathcal{S}_{k,\ell,n,c}^{\alpha,\mu,\theta} , \nabla \mathcal{S}_{j,\iota,m,c}^{\alpha,\mu,\theta} \big)_{r^{2\mu}(1-r^{2\theta})^{\alpha+1} }+c\big(\mathcal{S}_{k,\ell,n,c}^{\alpha,\mu,\theta} , \mathcal{S}_{j,\iota,m,c}^{\alpha,\mu,\theta}\big)_{r^{2\mu-2}(1-r^{2\theta})^{\alpha} } \nonumber\\
    		&\quad -\big(\nabla_0 \mathcal{S}_{k,\ell,n,c}^{\alpha,\mu,\theta} , \nabla_0 \mathcal{S}_{j,\iota,m,c}^{\alpha,\mu,\theta} \big)_{r^{2\mu+2\theta-2}(1-r^{2\theta})^{\alpha} } \nonumber\\
    		& =-(\|\bm{x}\|^{-2\mu}(1-\|\bm{x}\|^{2\theta})^{-\alpha} \nabla \cdot (\|\bm{x}\|^{2\mu}(1-\|\bm{x}\|^{2\theta})^{\alpha+1})\nabla\mathcal{S}_{k,\ell,n,c}^{\alpha,\mu,\theta}, \mathcal{S}_{j,\iota,m,c}^{\alpha,\mu,\theta})_{r^{2\mu}(1-r^{2\theta})^{\alpha} }\nonumber\\
    		&\quad +c(\|\bm{x}\|^{-2}\mathcal{S}_{k,\ell,n,c}^{\alpha,\mu,\theta},\mathcal{S}_{j,\iota,m,c}^{\alpha,\mu,\theta})_{r^{2\mu}(1-r^{2\theta})^{\alpha} }-(\|\bm{x}\|^{2\theta-2}\Delta_0\mathcal{S}_{k,\ell,n,c}^{\alpha,\mu,\theta},\mathcal{S}_{j,\iota,m,c}^{\alpha,\mu,\theta})_{r^{2\mu}(1-r^{2\theta})^{\alpha} }\nonumber\\
    		& =(\mathscr{D}_{c,\theta,\bm{x}}^{(\alpha,\mu)}\mathcal{S}_{k,\ell,n,c}^{\alpha,\mu,\theta},\mathcal{S}_{j,\iota,m,c}^{\alpha,\mu,\theta})_{r^{2\mu}(1-r^{2\theta})^{\alpha} } \nonumber\\
    		&=(\chi_{\theta,k}^{\alpha,\mu}\|\bm{x}\|^{2\theta-2}\mathcal{S}_{k,\ell,n,c}^{\alpha,\mu,\theta},\mathcal{S}_{j,\iota,m,c}^{\alpha,\mu,\theta})_{r^{2\mu}(1-r^{2\theta})^{\alpha} }\nonumber\\
    		&=\chi_{\theta,k}^{\alpha,\mu}(\mathcal{S}_{k,\ell,n,c}^{\alpha,\mu,\theta},\mathcal{S}_{j,\iota,m,c}^{\alpha,\mu,\theta})_{r^{2\mu+2\theta-2}(1-r^{2\theta})^{\alpha} }\nonumber,
    		\end{align}    	
    	   which completes the proof.
  	\end{proof}

 The following identity is essential for developing  the MBP spectral methods for the eigenvalue problems involving the degenerating elliptic  operator:  $-\nabla \cdot  \|\bm{x}\|^{2\mu} \nabla  + c\|\bm{x}\|^{2\mu-2}$. 
\begin{theorem}\label{theoMuntzII}
	Given the  parameters satisfying the same conditions as in Theorem \ref{theoSL},
	the MBPs $\{\mathcal{S}_{k,\ell,n,c}^{\alpha,\mu,\theta}(\bm{x})\}$ satisfy the   differential identity:
	\begin{align}
		(-\nabla\cdot \|\bm{x}\|^{2\mu} \nabla  & +c \|\bm{x}\|^{2\mu-2}) \mathcal{S}^{\alpha,\mu,\theta}_{k,\ell,n,c}(\bm{x})\nonumber\\
		&=-4\theta^2(k+\beta_n)(k+\alpha+\beta_n+1)\|\bm{x}\|^{2\theta-2+2\mu} \mathcal{S}^{\alpha+2,\mu,\theta}_{k-1,\ell,n,c}(\bm{x}).\label{eigMGGFI}
	\end{align}
\end{theorem}
\begin{proof}	
	Using the Leibniz rule for gradient and divergence, we derive
	\begin{equation}\label{rthetaA}
	\begin{aligned}
		-\nabla \cdot \|\bm{x}\|^{2\mu}\nabla +c \|\bm{x}\|^{2\mu-2}&= -\|\bm{x}\|^{2\mu}\nabla\cdot\nabla-2\mu\|\bm{x}\|^{2\mu-2}\cdot \bm{x}\nabla+c \|\bm{x}\|^{2\mu-2}\\	
		&=-r^{2\mu}\Big(\partial^2_r+ \frac{2\mu+d-1}{r} \partial_r + \frac{\Delta_0-c}{r^2}\Big).
	\end{aligned}
	\end{equation}
Thus applying  the operator $-\nabla \cdot \|\bm{x}\|^{2\mu}\nabla +c \|\bm{x}\|^{2\mu-2}$ to the MBPs, we  obtain from   \refe{GGPeq37} and \refe{GGPeq38} that in $r$-direction,
	\begin{equation*}
		\begin{aligned}
			&- \frac{r^{2\mu}}{r^{\theta\beta_n+1-\frac{d}{2}-\mu} } \Big[ \partial^2_r+ \frac{2\mu+d-1}{r} \partial_r - \frac{n(n+d-2)+c}{r^2} \Big] [ r^{\theta\beta_n+1-\frac{d}{2}-\mu}   P^{(\alpha,\beta_n)}_k(2r^{2\theta}-1) ]\\
			&\quad = \, -r^{2\mu}\bigg\{\Big[\partial_r^2   + \frac{2\theta\beta_n+1}{r} \partial_r\Big] P^{(\alpha,\beta_n)}_k(2r^{2\theta}-1)\\
 &			\qquad +\frac{(\theta\beta_n+1-d/2-\mu)(\theta\beta_n-1+d/2+\mu)-n(n+d-2)-c}{r^2} P^{(\alpha,\beta_n)}_k(2r^{2\theta}-1)\bigg\}\\
			&\quad =  -r^{2\mu} \Big[\partial_r^2   + \frac{2\theta\beta_n+1}{r} \partial_r  \Big]P^{(\alpha,\beta_n)}_k(2r^{2\theta}-1)			,
		\end{aligned}
	\end{equation*}
	where we chose the value of $\beta_n$ in \eqref{betank} so that the coefficient of $1/r^2$ vanishes.
	Making the change of  variable  $\rho =2r^{2\theta}-1$, we  derive from \eqref{LemEQ1} that
	\begin{align}
		&-r^{2\mu} \Big[\partial_r^2   + \frac{2\theta\beta_n+1}{r} \partial_r  \Big]P^{(\alpha,\beta_n)}_k(2r^{2\theta}-1)\nonumber\\
		&=-2^{2+\frac{1}{\theta}-\frac{\mu}{\theta}}\theta^2\Big[(1+\rho)^{2-\frac{1}{\theta}+\frac{\mu}{\theta}}\partial^2_\rho+(1+\beta_n)(1+\rho)^{1-\frac{1}{\theta}+\frac{\mu}{\theta}}\partial_\rho\Big] P^{(\alpha,\beta_n)}_k (\rho)\nonumber\\
		&=-2^{2+\frac{1}{\theta}-\frac{\mu}{\theta}}\theta^2(1+\rho)^{1-\frac{1}{\theta}+\frac{\mu}{\theta}}\Big[(1+\rho)\partial^2_\rho+(1+\beta_n)\partial_\rho\Big]P^{(\alpha,\beta_n)}_k (\rho)\nonumber\\
		&=-2^{\frac{1}{\theta}-\frac{\mu}{\theta}}\theta^2(1+\rho)^{1-\frac{1}{\theta}+\frac{\mu}{\theta}}(k+\alpha+\beta_n+1)\nonumber
		\\&\quad\times\Big[(k+\alpha+\beta_n+2)(1+\rho) P^{(\alpha+2,\beta_n+2)}_{k-2}(\rho)+2(1+\beta_n) P^{(\alpha+1,\beta_n+1)}_{k-1} (\rho)\Big]
		\label{GGPseq1}.
	\end{align}
A combination of \eqref{LemEQ2} and \eqref{LemEQ3}, together with the symmetry of Jacobi polynomials,  leads to
	\begin{align}
		\label{symEQ} (k+\frac{\alpha + \beta}{2} +1)(1+x)P_k^{(\beta, \alpha+1)}(x)=(k+\alpha+1)P_k^{(\beta, \alpha)}(x)+(k+1)P_{k+1}^{(\beta, \alpha)}(x), \\
		\label{symEQ1} (2k+ \alpha+\beta +1)P_k^{(\beta, \alpha)}(x)  = (k+\alpha+\beta+1)P_k^{(\beta, \alpha+1)}(x) + (k+\beta)P_{k-1}^{(\beta, \alpha+1)}(x).
	\end{align}
	Then using \eqref{LemEQ2} and \eqref{symEQ}-\eqref{symEQ1}, we  can further deduce from \eqref{GGPseq1} that
	\begin{align}
		&2^{1+\frac{1}{\theta}-\frac{\mu}{\theta}}\theta^2(1+\rho)^{1-\frac{1}{\theta}+\frac{\mu}{\theta}}\frac{(k+\alpha+\beta_n+1) (k+\beta_n) }{2k+\alpha+\beta_n+1} \\&\quad\times\Big[ (k+\alpha+\beta_n+2)P^{(\alpha+2,\beta_n+1)}_{k-1}(\rho)  +(k+\alpha+1)  P^{(\alpha+2,\beta_n+1)}_{k-2}(\rho) \Big]\nonumber\\
		&=2^{1+\frac{1}{\theta}-\frac{\mu}{\theta}}\theta^2(1+\rho)^{1-\frac{1}{\theta}+\frac{\mu}{\theta}}(k+\beta_n)(k+\alpha+\beta_n+1) P^{(\alpha+2,\beta_n)}_{k-1}(\rho).\nonumber
	\end{align}
In view of \eqref{rthetaA}, we obtain
\begin{align}
		&(-\nabla\cdot \|\bm{x}\|^{2\mu} \nabla +c\|\bm{x}\|^{2\mu-2}) \mathcal{S}^{\alpha,\mu,\theta}_{k,\ell,n,c}(\bm{x}) \nonumber\\
		& =-4\theta^2(k+\beta_n)(k+\alpha+\beta_n+1)\|\bm{x}\|^{2\theta-2+2\mu} \mathcal{S}^{\alpha+2,\mu,\theta}_{k-1,\ell,n,c}(\bm{x}). \nonumber
	\end{align}
	This ends the proof.
\end{proof}

\begin{rem}
It is noteworthy that the following two special cases  have been studied in literature.
\begin{itemize}
	\item[(i)] When $\mu=0$, $\theta=1$ and $c=0$, the MBPs  $\{\mathcal{S}_{k,\ell,n,0}^{\alpha,0,1}(\bm{x})\}$ degenerate to the ball polynomials  $\{P_{k,\ell}^{\alpha,n}(\bm{x})\}$ defined in \eqref{ballpoly},  and
	\begin{equation}\label{balleigen}
		\begin{aligned}
			\mathscr{D}_{0,1,\bm{x}}^{(\alpha,0)}\mathcal{S}_{k,\ell,n,0}^{\alpha,0,1}(\bm{x})&=\big[-(1-\|\bm{x}\|^2)^{-\alpha} \nabla \cdot ((1-\|\bm{x}\|^2)^{\alpha+1})\nabla-\Delta_0\big]\mathcal{S}_{k,\ell,n,0}^{\alpha,0,1}(\bm{x})\\
			&=(n+2k)(n+2k+2\alpha+d)\mathcal{S}_{k,\ell,n,0}^{\alpha,0,1}(\bm{x}).
		\end{aligned}
	\end{equation}
	We refer to   \cite[Proposition 11.1.15]{Dai2013approximation} for this type of results and more discussions.
\medskip
\item[(ii)] When $\mu=0$, $\theta=1/2$ and $c\neq0$, the MBPs  $\{\mathcal{S}_{k,\ell,n,c}^{\alpha,0,1/2}(\bm{x})\}$ reduce to the spectral basis functions $Q_{k,\ell}^{\alpha,n}$ studied  in \cite{Ma2018Efficient}.
\end{itemize}
\end{rem}
\medskip

 To conclude this section, we highlight the following property of the MBPs with 	 $\alpha=-1,$ which  will be useful for the algorithm developing in the applications to  eigenvalue problems. We provide the derivation in Appendix \ref{prolemmamuntzorth}.
\begin{lemma}\label{lemmuntzorth}
	For $\mu>-1/2$, $\theta>0$, and $k,j \in\mathbb{N}_0$, $(\ell,n), (\iota,m)\in \Upsilon_{\infty}^d$, the MBPs $\{\mathcal{S}_{k, \ell,n,c}^{-1,\mu,\theta}(\bm{x})\}$ form a Sobolev orthogonal basis in the sense that
	\begin{equation}\label{muntzstiffmatrix1}
		\begin{array}{lll}
			&	\big(\nabla \mathcal{S}_{k, \ell,n,c}^{-1,\mu,\theta}, \nabla \mathcal{S}_{j, \iota,m,c}^{-1,\mu,\theta}\big)_{r^{2\mu}}+c\big(\mathcal{S}_{k, \ell,n,c}^{-1,\mu,\theta},  \mathcal{S}_{j, \iota,m,c}^{-1,\mu,\theta}\big)_{r^{2\mu-2}}  \vspace{1ex}\\
			&\quad = \delta_{k j} \delta_{nm} \delta_{\ell \iota} \Big[\frac{2\theta(k+\beta_{n})^2}{2k+\beta_{n}}(1-\delta_{k0})+\big(\theta \beta_{n}+1-\frac{d}{2}-\mu\big) \delta_{k0}\Big],
		\end{array}
	\end{equation}
where $\{\beta_n\}$ are given in \eqref{betank}.
\end{lemma}

\section{MBP spectral-Galerkin methods for eigenvalue problems}\label{sect4:gal}
\setcounter{equation}{0}
\setcounter{lmm}{0}
\setcounter{thm}{0}

\subsection{Degenerate eigenvalue problems with singular potentials}
We start with  the eigenvalue problem: find $\{\lambda,u\}$ with $u\neq0$ such that
\begin{equation}\label{eigensingular}
\begin{cases}
		{\mathcal L}_\nu[u](\bm x):=-\nabla \cdot  \|\bm{x}\|^{2\mu} \nabla u(\bm{x}) + c\|\bm{x}\|^{2\mu-2}u (\bm{x})= \lambda u(\bm{x}),&\quad  \bm{x}  \in  {\mathbb B}^d,\\
		u(\bm{x})=0, &\quad  \bm{x} \in  {\mathbb S}^{d-1},
\end{cases}
\end{equation}
for given $\mu>-1/2$ and $c>0.$  The linear operator ${\mathcal L}_\nu$  is an elliptical operator with a perturbed ellipticity due to the involvement of
degenerate coefficients, as $\inf_{{\mathbb B}^d}\|\bm{x}\|^{2\mu}= 0$ or  $\sup_{{\mathbb B}^d}\|\bm{x}\|^{2\mu}= \infty$    (cf.\! \cite{Kufner1985Weighted}).
According to the standard theory, it admits a countable set of eigenvalues
$$
0<\lambda_1\le \lambda_2\le \cdots \to \infty,
$$
but the eigenfunctions are typically singular at the origin.
It is challenging but necessary  to develop spectrally accurate  method to compute as many as trustable eigenvalues as possible.
It is important to remark that the eigen-problem with $\mu=0,$ is known as the Schr\"odinger eigenvalue problem with inverse square potential.  The inverse square potential possesses the same homogeneity or ``differential order" as the Laplacian, while it usually invokes strong singularities of the Schrödinger eigenfunctions and thus cannot be treated as a lower-order perturbation term \cite{Cao2006Solutions,Felli2007Schorodinger,Felli2006Elliptic}. 
 On the other hand,  it is crucial to consider this problem in the special domain ${\mathbb B}^d,$
as it is an essential building block for the (non-standard) spectral-element method on a general domain (see Li and Zhang \cite{Li2017efficient}).

To further motivate why we propose the MBP spectral methods and how we properly select the parameters, we first conduct some analytical study.
We sketch the derivation in Appendix \ref{proanasolution}.
\begin{proposition}\label{anasolution} {\em The eigenvalues of the problem \eqref{eigensingular} are determined by the zeros of the Bessel function of order
$\beta_n/2$ with $\beta_n$ given in \eqref{betank}, i.e.,
\begin{equation}\label{anaeigvalue}
		J_{\frac{\beta_n}{2}}\Big(\frac{\sqrt{\lambda}}{1-\mu}\Big)=0,
	\end{equation}
and	the corresponding eigenfunction has the series expression
	\[
	u(\bm{x})=\sum_{n=0}^{\infty} \sum_{\ell=1}^{a_{n}^{d}} u_{\ell}^{n}(r) Y_{\ell}^{n}(\bm{\hat{x}}),
	\]
	where 
	\begin{equation}\label{anaeigfun}
		u_{\ell}^{n}(r) =r^{(1-\mu)(\beta_n/2+1)-d/2}\sum_{m=0}^{\infty}  c_m^n \lambda^{m+\beta_n/4} (r^{1-\mu})^{2m},
	\end{equation}
with the constant $c_m^n$ given in \eqref{analysissolution}.
}
\end{proposition}

It is evident from \eqref{anaeigfun} that  $u_{\ell}^{n}$  has  a singular  behaviour  of the form $r^{(1-\mu)(\beta_n/2+1)-d/2} p_\ell^n(r^{1-\mu})$ where $p_\ell^n(z)$ is a smooth function.
It is known that any approximation by polynomials  to such  singular functions as $u_{\ell}^{n}$ has a very limited convergence.
However, they can be best approximation by the MBPs by choosing the parameters so as to capture the singular factors.

With this in mind, we employ the MBP approximation with $\theta=1-\mu$ and $\alpha=-1$ (to meet the homogeneous boundary condition).
Define the approximation space
\begin{equation*}
	\mathcal{V}_{\!N,K}=\text{span}\big\{S_{k,\ell,n,c}^{-1,\mu,1-\mu}(\bm{x})\,: \, (\ell,n)\in \Upsilon_{N}^d,\, 1\leq k\leq K, \, k \in \mathbb{N}_{0}\big\}.
\end{equation*}
The MBP spectral-Galerkin scheme for \eqref{eigensingular} is to find $\lambda_{N,K} \in \mathbb{R}$ and $u_{N,K}\in \mathcal{V}_{\!N,K} \backslash \{0\}$, such that
\begin{equation}\label{eigprosingular1}
	(\|\bm x\|^{2\mu}\nabla u_{N,K},\nabla v_{N,K})+c(\|\bm x\|^{2\mu-2}u_{N,K},v_{N,K})=\lambda_{N,K}(u_{N,K}, v_{N,K}), \quad \forall\, v_{N,K} \in \mathcal{V}_{\!N,K}.
\end{equation}
In implementation, we write
\begin{equation}\label{uNKr}
	u_{N,K}(\bm{x})=\sum_{n=0}^N\sum_{\ell=1}^{a_n^d}\sum_{k=1}^{K}\hat{u}_{k,\ell}^nS_{k,\ell,n,c}^{-1,\mu,1-\mu}(\bm{x}),
\end{equation}
and denote
\begin{equation*}
	\begin{array}{l}
		\boldsymbol{u}=\big( \hat{\boldsymbol{u}}_{1}^{0}, \hat{\boldsymbol{u}}_{2}^{0}, \cdots, \hat{\boldsymbol{u}}_{a_{0}^{d}}^{0}, \cdots, \hat{\boldsymbol{u}}_{1}^{N}, \hat{\boldsymbol{u}}_{2}^{N}, \cdots, \hat{\boldsymbol{u}}_{a_{N}^{d}}^{N}\big)^{t}, \quad \hat{\boldsymbol{u}}_{\ell}^{n}=\big(\hat{u}_{1, \ell}^{n}, \hat{u}_{2, \ell}^{n}, \ldots, \hat{u}^n_{K, \ell}\big)^{t}.
	\end{array}
\end{equation*}
Accordingly, we denote the stiffness and the mass matrices by $\bm{S}$ and $\bm{M}$, respectively, with the corresponding  entries
\begin{align}
	&\bm{S} \to (\|\bm x\|^{2\mu}\nabla S_{k,\ell,n,c}^{-1,\mu,1-\mu},\nabla S_{k,\ell,n,c}^{-1,\mu,1-\mu})+c(\|\bm x\|^{2\mu-2}S_{k,\ell,n,c}^{-1,\mu,1-\mu},S_{k,\ell,n,c}^{-1,\mu,1-\mu}),\nonumber\\
	&\bm{M}\to ( S_{k,\ell,n,c}^{-1,\mu,1-\mu},S_{k,\ell,n,c}^{-1,\mu,1-\mu}),\nonumber
\end{align}
and the linear system of the problem \eqref{eigprosingular1} reads 
\[\bm{S}\boldsymbol{u}=\lambda_N\bm{M}\boldsymbol{u}.\]

Importantly, we can evaluate the stiffness matrix explicitly and provide some details  in Appendix \ref{prolemmamuntz1}.
\begin{lemma}\label{lemmamuntz1}
	For fixed $k,j \in\mathbb{N}_0$, $(\ell,n),(\iota,m)\in \Upsilon_{\infty}^d$, it holds that
	\begin{equation}\label{muntzmassmatrix1}
			\big(S_{k,\ell,n,c}^{-1,\mu,1-\mu},  S_{k,\ell,n,c}^{-1,\mu,1-\mu}\big)= \dfrac{1}{1-\mu}\delta_{nm} \delta_{\ell \iota}
			\times
			\begin{cases}
				\frac{1}{2(1+\beta_n)}, & k=j=0, \\
				\frac{(k+\beta_n)^2}{(2k+\beta_n-1)(2k+\beta_n)(2k+\beta_n+1)}, \quad & k=j\geq 1, \\[6pt]
				-\frac{(k+\beta_n)(k+\beta_n-1)}{2(2k+\beta_n-2)(2k+\beta_n-1)(2k+\beta_n)}, \quad & k=j+1, \\[6pt]
				-\frac{(k+\beta_n)(k+\beta_n+1)}{2(2k+\beta_n)(2k+\beta_n+1)(2k+\beta_n+2)}, \quad & k=j-1, \\ 
				0, & \text {otherwise}.
			\end{cases}
	\end{equation}
\end{lemma}

 From Lemma \ref{lemmuntzorth} with $\theta=1-\mu$, one have that the stiffness matrix $\bm{S}$ is  a diagonal matrix. While one can see from Lemma \ref{lemmamuntz1} that the mass matrix $\bm{M}$ is tridiagonal matrix.

\begin{figure*}[!htb]
	\centering
	{\includegraphics[width=7cm,height=5cm]{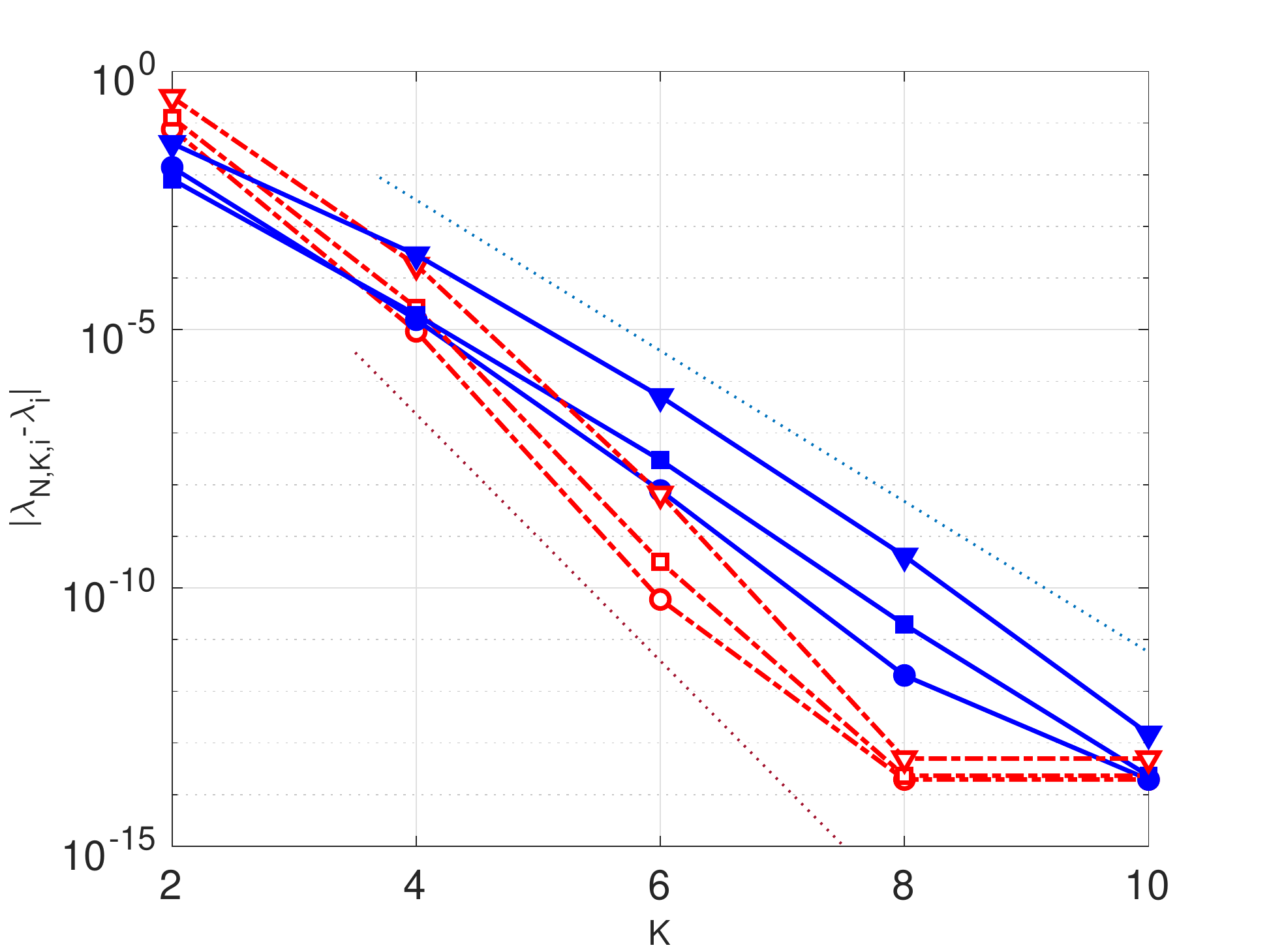}}
	{\includegraphics[width=7cm,height=5cm]{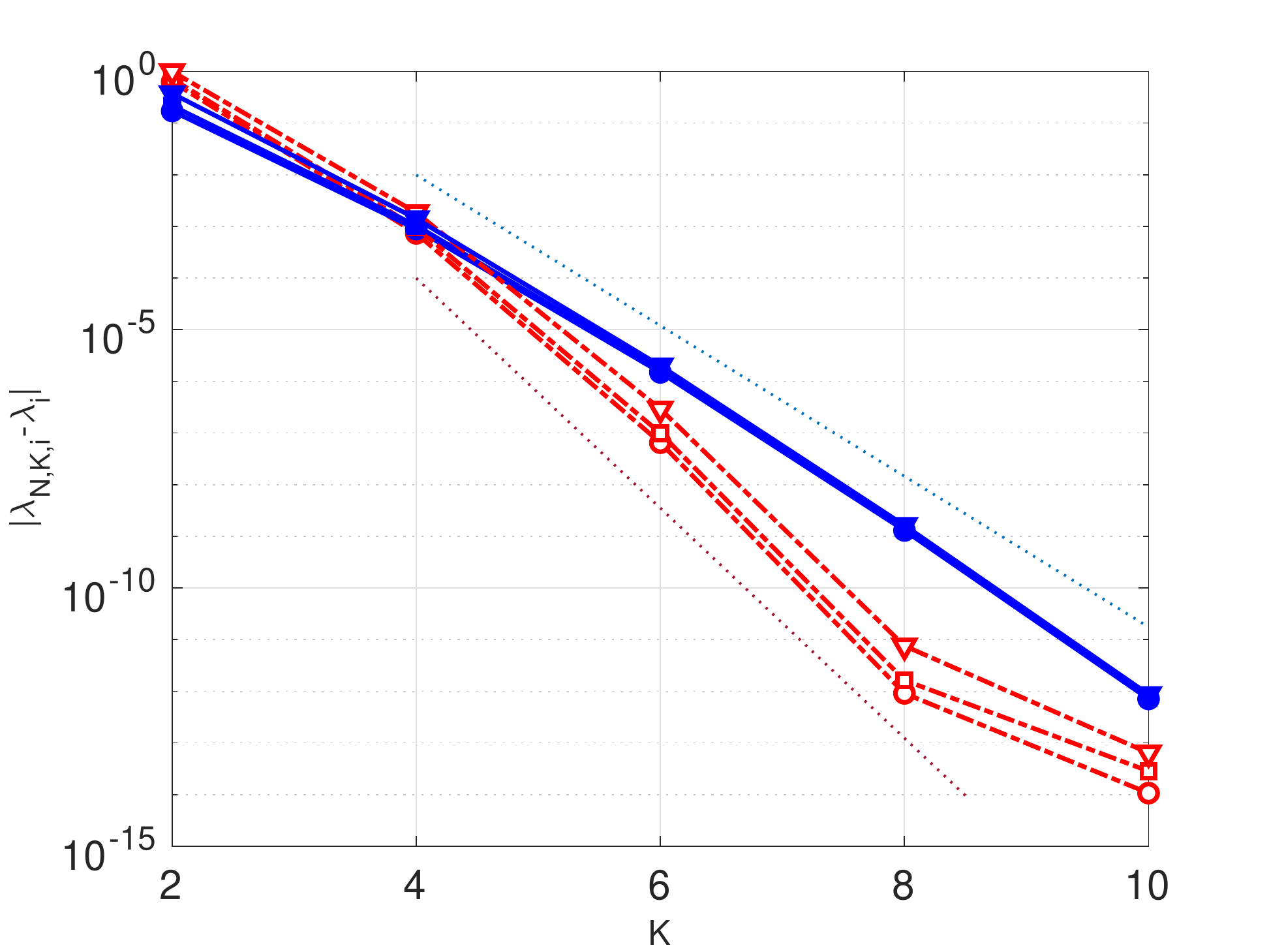}} 
	\caption{\small Approximation errors $|\lambda_{N,K,i}-\lambda_i|$ for  \eqref{eigensingular} with $\mu=1/2$ versus $K$ by the MBPs spectral method using basis $S_{k,\ell,n,c}^{-1,\mu,1-\mu}(\bm{x})$ (primitive markers with dotted line) and $S_{k,\ell,n,c}^{-1,\mu,(1-\mu)/2}(\bm{x})$ (filled markers with full line)   on the unit disk. $\circ:\lambda_1(n=0); \square:\lambda_2(n=1);\vartriangle:\lambda_3(n=2)$.( Left: $c=2$. Right: $c=10$).}
	\label{figerrdegenerate}
\end{figure*}

\begin{figure*}[!htb]
	\centering
	{\includegraphics[width=7cm,height=5cm]{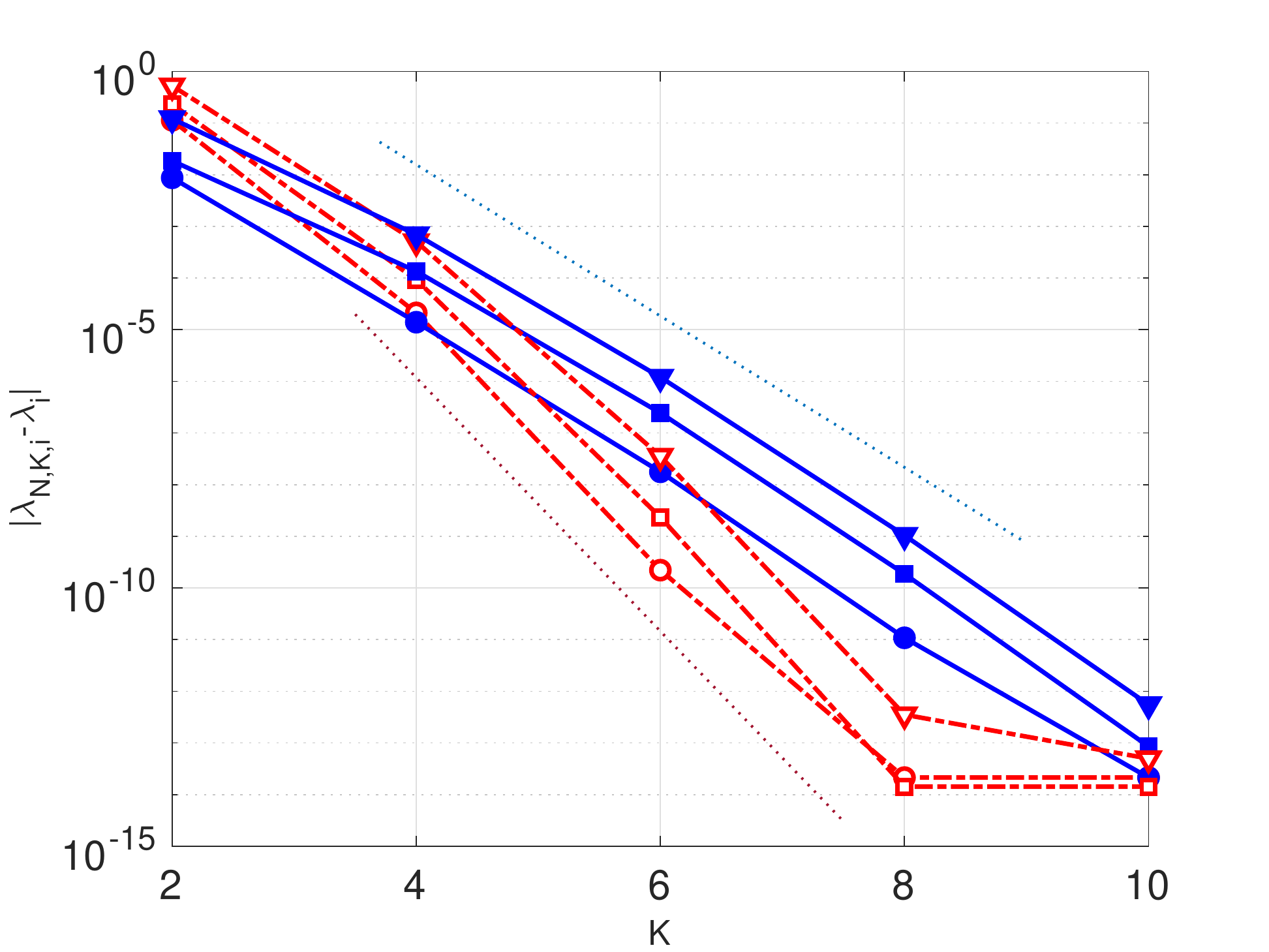}}
	{\includegraphics[width=7cm,height=5cm]{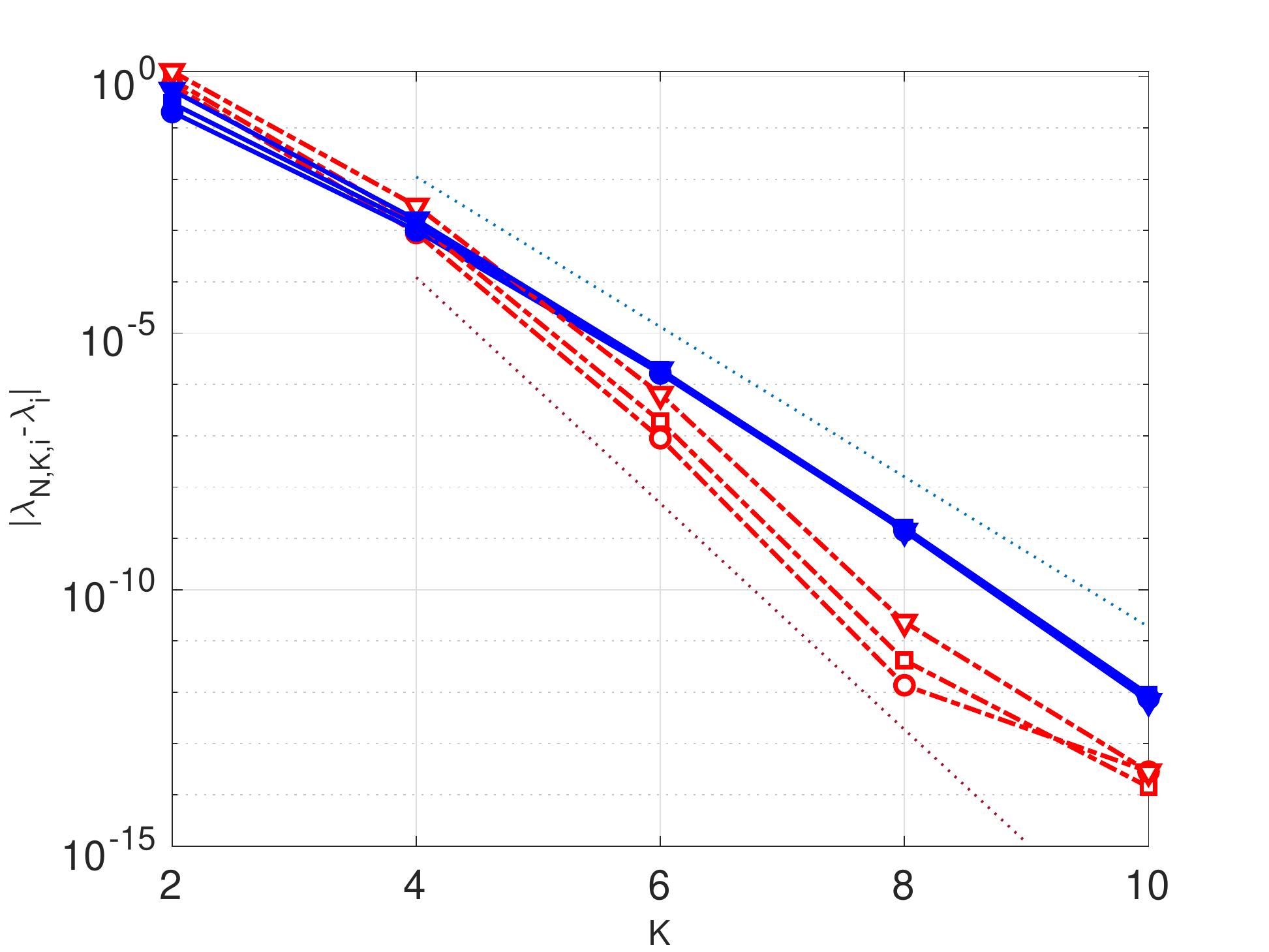}} 
	\caption{\small Approximation errors $|\lambda_{N,K,i}-\lambda_i|$ for  \eqref{eigensingular} with $\mu=1/2$ versus $K$ by the MBPs spectral method using basis $S_{k,\ell,n,c}^{-1,\mu,1-\mu}(\bm{x})$ (primitive markers with dotted line) and $S_{k,\ell,n,c}^{-1,\mu,(1-\mu)/2}(\bm{x})$ (filled markers with full line) on the unit ball. $\circ:\lambda_1(n=0); \square:\lambda_2(n=1);\vartriangle:\lambda_3(n=2)$.( Left: $c=2$. Right: $c=10$).}
	\label{figerrdegenerate3d}
\end{figure*}

\begin{remark} Observe  from \eqref{anaeigfun} that the non-singular series of $u_{\ell}^{n}(r)$
can  be accurately  approximated by the polynomials of $r^{1-\mu}$.  Thus we can also choose the MBPs $\{S_{k,\ell,n,c}^{-1,\mu,(1-\mu)/2}(\bm{x})\}$, i.e.,
with the parameter $\theta=(1-\mu)/2$ to capture the singularity in the eigenfunctions. Similar to Lemma \ref{lemmamuntz1}, we can show  that the mass matrix $\bm{M}$ with MBPs $\{S_{k,\ell,n,c}^{-1,\mu,(1-\mu)/2}(\bm{x})\}$ is pentadiagonal matrix, and from Lemma \ref{lemmuntzorth} with $\theta=(1-\mu)/2$, one can get that the stiffness matrix $\bm{S}$ is  a diagonal matrix.
 In Figures \ref{figerrdegenerate}-\ref{figerrdegenerate3d}, one will see that the convergence order of the method using basis function $S_{k,\ell,n,c}^{-1,\mu,1-\mu}(\bm{x})$  is slightly  higher than that of the method using basis function $S_{k,\ell,n,c}^{-1,\mu,(1-\mu)/2}(\bm{x})$.
	\end{remark}

\begin{table}[!htb]
	\caption{\small The eigenvalues $\lambda$ of \eqref{eigensingular}  with $\mu=1/2$ in $d=2$ }\small
	\begin{tabular}{|ccc|c|c|}
		\cline{1-5}
		$c$   & $n$ & $k$ & Exact &  Numerical   \\ \cline{1-5}
		2   & 1 & 0 &12.6566911210566      &  12.6566911210566    \\
		2   & 1 & 1 & 27.8493337022154       &   27.8493337022154   \\
		2   & 1 & 2 & 47.8938240898394      &   47.8938240898394 \\
		2   & 2 & 0 & 19.2347320834118       &  19.2347320834119   \\
		0.1 & 0 & 0 & 4.15524482735467      &  4.15524482735467    \\
		1   & 0 & 0 & 7.38102583563937      &7.38102583563937  \\
		4   & 0 & 0 & 14.9582601885591      & 14.9582601885591 \\
		10  & 0 & 0 & 27.0470413068364       &27.0470413068364  \\ \cline{1-5}
	\end{tabular}
	\label{tableeigenvalue2d}
\end{table}

\begin{table}[!htb]
	\caption{\small The eigenvalues $\lambda$ of \eqref{eigensingular} with  $\mu=1/2$ in $d=3$}\small
	\begin{tabular}{|ccc|c|c|}
		\cline{1-5}
		$c$   & $n$ & $k$ & Exact  & Numerical  \\ \cline{1-5}
		2   & 1 & 0 &16.6039682455504      &  16.6039682455504     \\
		2   & 1 & 1 & 34.0437745462078     &  34.0437745462078  \\
		2   & 1 & 2 & 56.2700086174140      &   56.2700086174140 \\
		2   & 2 & 0 & 24.6815681193123       &  24.6815681193123  \\
		0.1 & 0 & 0 & 6.91481384026372     &  6.91481384026372   \\
		1   & 0 & 0 & 9.51612890626288      &9.51612890626288   \\
		4   & 0 & 0 & 16.6039682455504      & 16.6039682455504\\
		10  & 0 & 0 & 28.4407806172599      &28.4407806172599  \\ \cline{1-5}
	\end{tabular}
	\label{tableeigenvalue3d}
\end{table}

\begin{table}[!htb]
	\caption{\small Radial normalized eigenfunctions of \eqref{eigensingular} with $n=0$, $\ell=1$, $\mu=1/2$  in $d=2$}\small
	\begin{tabular}{|cc|c|c|}
		\cline{1-4}
		$r$   & $c$ & Exact &  Numerical \\ \cline{1-4}
		0.1   & 1  &0.027609083125293     &  0.027609083125293   \\
		0.5   & 1  & 0.024094590264357      &  0.024094590264357  \\
		0.8   & 1  & 0.008294196243488      &   0.008294196243488 \\
		0.1   & 2 & 0.019361071632968      &  0.019361071632964  \\
		0.5   & 2  & 0.027859282291842       &  0.027859282291843  \\
		0.8   & 2  & 0.010413071837381       &  0.010413071837381  \\ \cline{1-4}
	\end{tabular}
	\label{tableeigen}
\end{table}

In Tables \ref{tableeigenvalue2d}- \ref{tableeigenvalue3d}, we tabulate the numerical eigenvalues for  $\mu=1/2$ obtained by the MBP spectral method and the analytical values of the eigenvalues obtained by the analytical expression \eqref{anaeigvalue} for various choices of $c$, $n$, $k$ in $d=2$, and $d=3$, respectively.  It is seen that our method is spectrally accurate. The approximations errors for the first third eigenvalues of the MBP spectral method with basis $S_{k,\ell,n,c}^{-1,\mu,1-\mu}(\bm{x})$ and $S_{k,\ell,n,c}^{-1,\mu,(1-\mu)/2}(\bm{x})$ are plotted in Figure \ref{figerrdegenerate} in simi-log scale for both $c=2$ and $c=10$ in $d=2$ dimensions for \eqref{eigensingular} with $\mu=1/2$, while in Figure \ref{figerrdegenerate3d}, we give the approximation errors for the first third eigenvalues in simi-log scale for both $c=2$ and $c=10$ in $d=3$ dimensions for \eqref{eigensingular} with $\mu=1/2$. We see that the two MBP spectral methods share the spectral accuracy.

\begin{figure}[!htb]
	\subfigure[$(\mu,n,k,\ell)=(0.5,0,0,1)$]{
		\begin{minipage}[t]{0.45\textwidth}
			\centering \rotatebox[origin=cc]{-0}{\includegraphics[width=2.8in,height=1.45in]{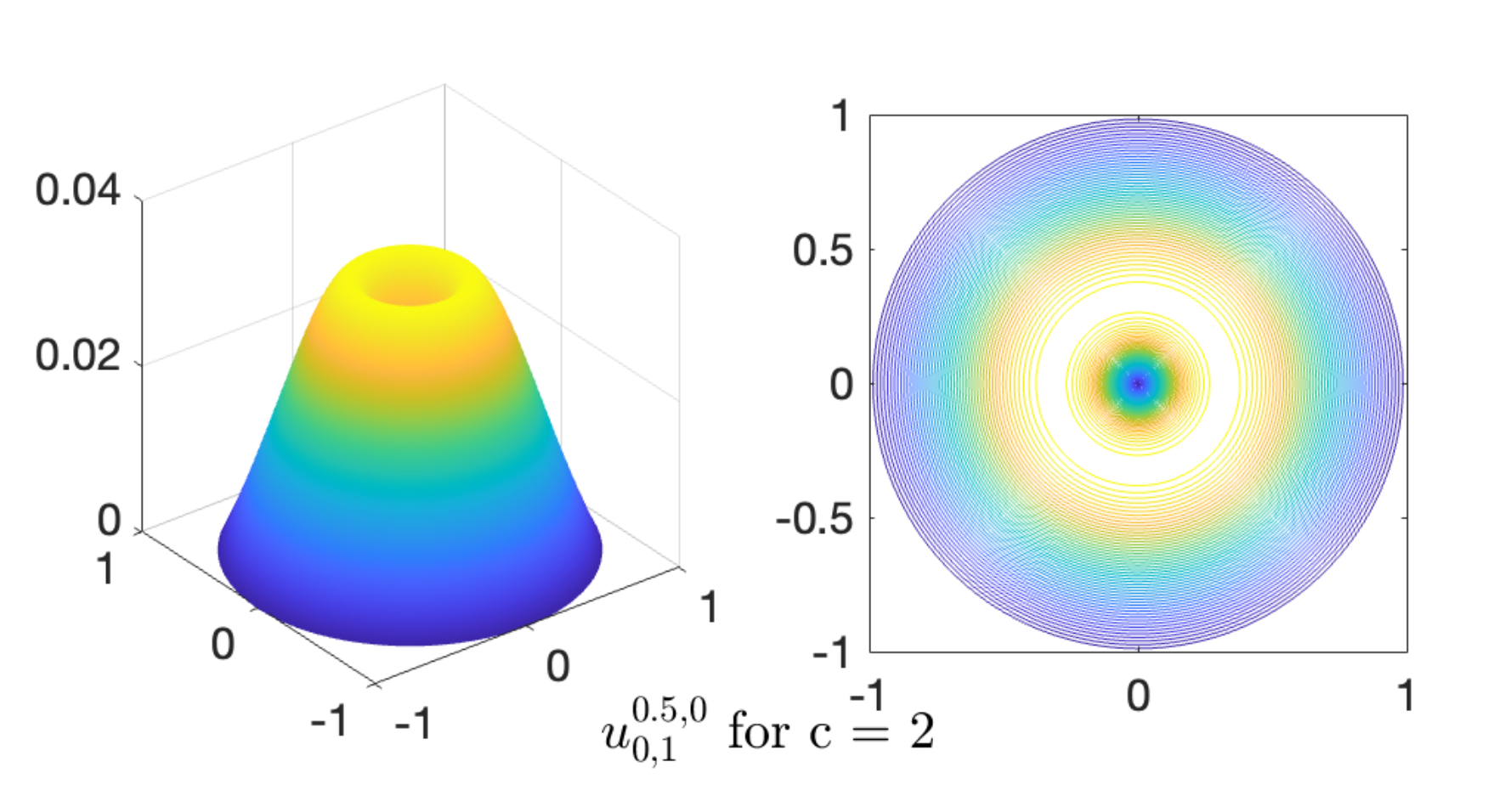}}
	\end{minipage}}%
	\subfigure[$(\mu,n,k,\ell)=(0.5,0,1,1)$]{
		\begin{minipage}[t]{0.45\textwidth}
			\centering \rotatebox[origin=cc]{-0}{\includegraphics[width=2.8in,height=1.45in]{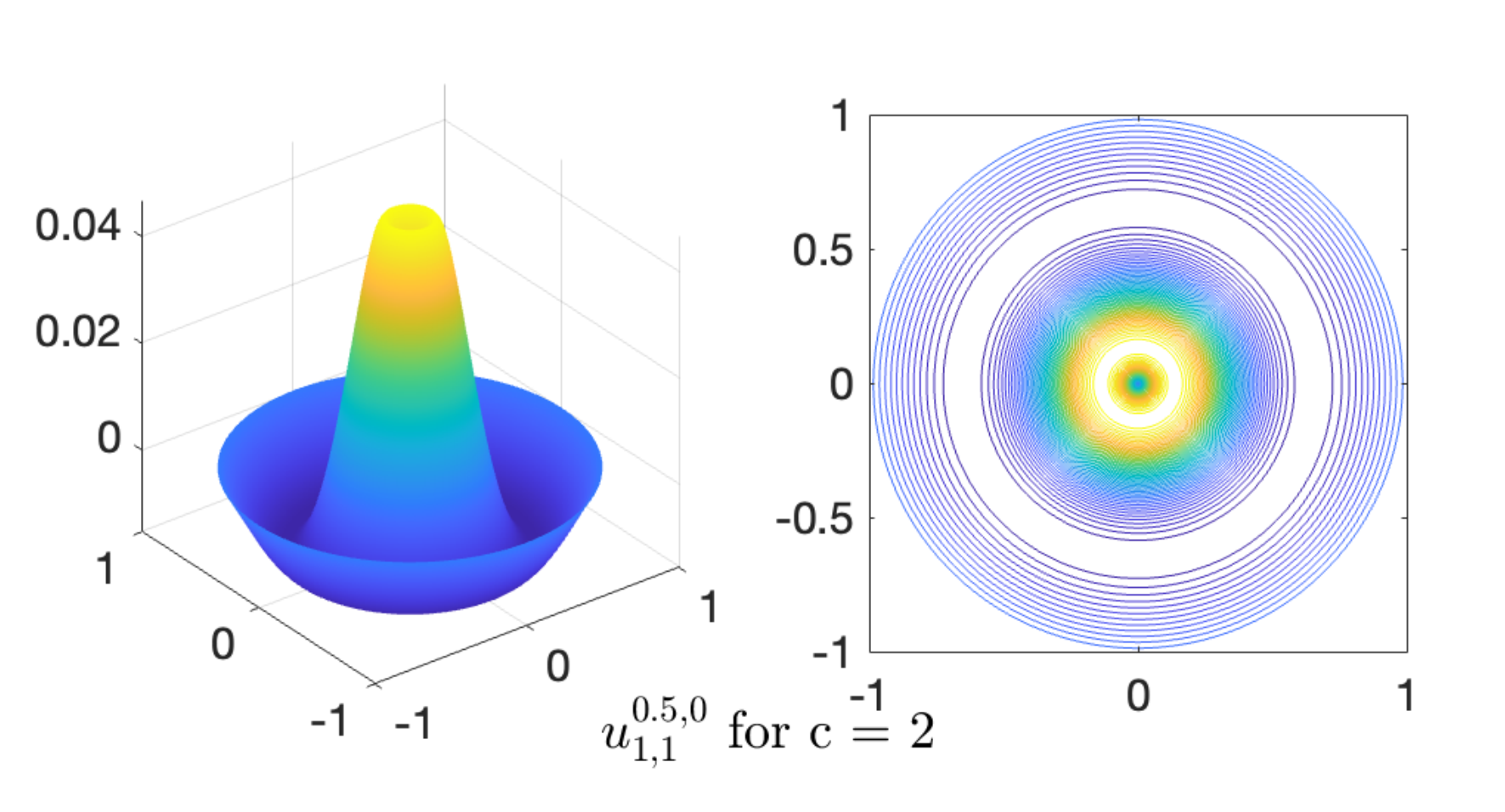}}
	\end{minipage}}
	\subfigure[$(\mu,n,k,\ell)=(0.5,0,2,1)$]{
		\begin{minipage}[t]{0.45\textwidth}
			\centering \rotatebox[origin=cc]{-0}{\includegraphics[width=2.8in,height=1.45in]{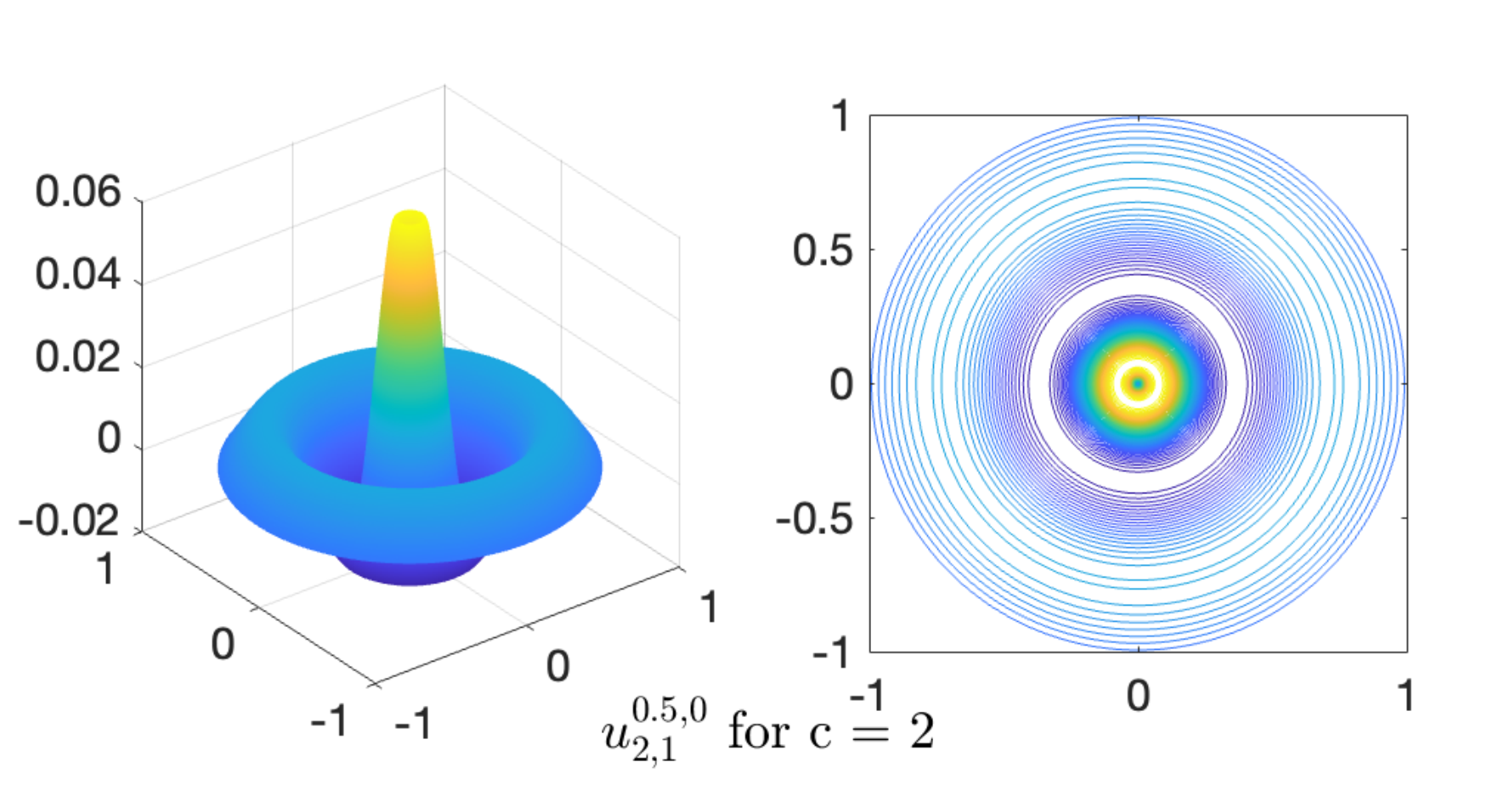}}
	\end{minipage}}%
	\subfigure[$(\mu,n,k,\ell)=(0.5,0,3,1)$]{
		\begin{minipage}[t]{0.45\textwidth}
			\centering \rotatebox[origin=cc]{-0}{\includegraphics[width=2.8in,height=1.45in]{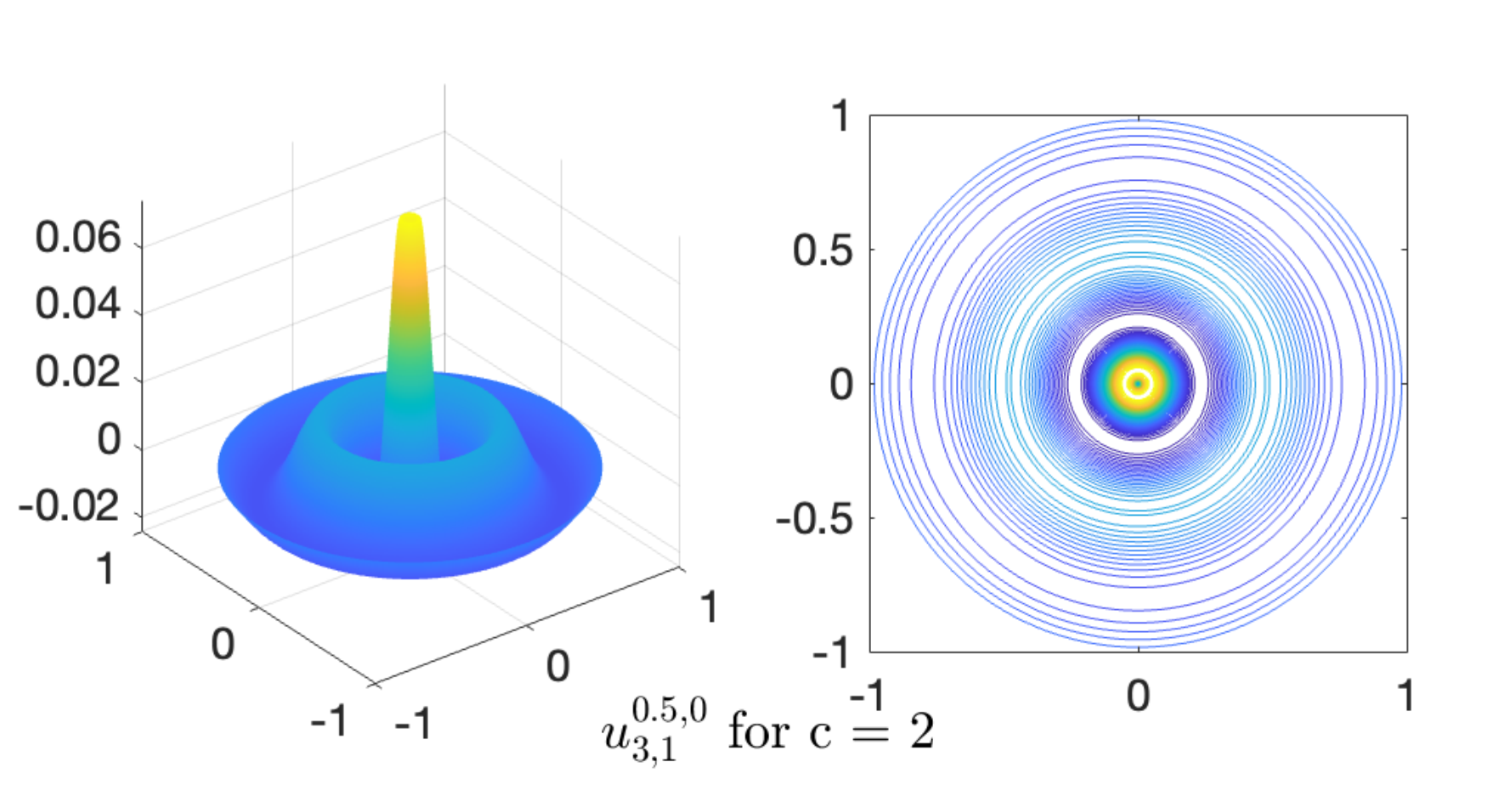}}
	\end{minipage}}%
	\caption {\small Normalized eigenfunctions $\{u_{k,\ell}^{\mu, n}(\bm{x},c)\}$ with  $c=2$, $\mu=1/2$ in $d=2$}\label{GGPs2d-1}
\end{figure}

\begin{figure}[!htb]
	\subfigure[$(\mu,n,k,\ell)=(0.5,1,0,1)$]{
		\begin{minipage}[t]{0.45\textwidth}
			\centering \rotatebox[origin=cc]{-0}{\includegraphics[width=2.8in,height=1.45in]{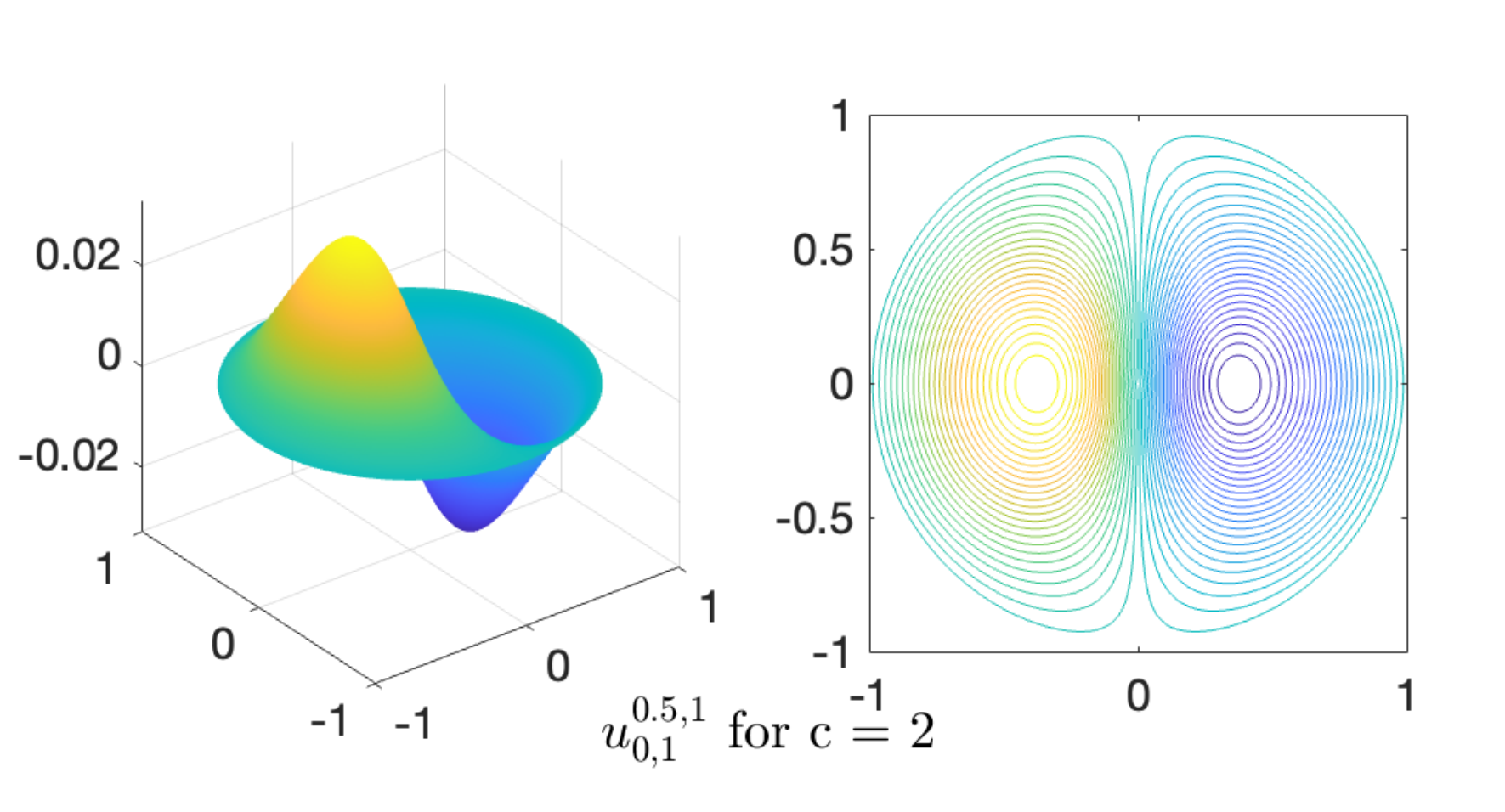}}
	\end{minipage}}%
	\subfigure[$(\mu,n,k,\ell)=(0.5,1,0,2)$]{
		\begin{minipage}[t]{0.45\textwidth}
			\centering \rotatebox[origin=cc]{-0}{\includegraphics[width=2.8in,height=1.45in]{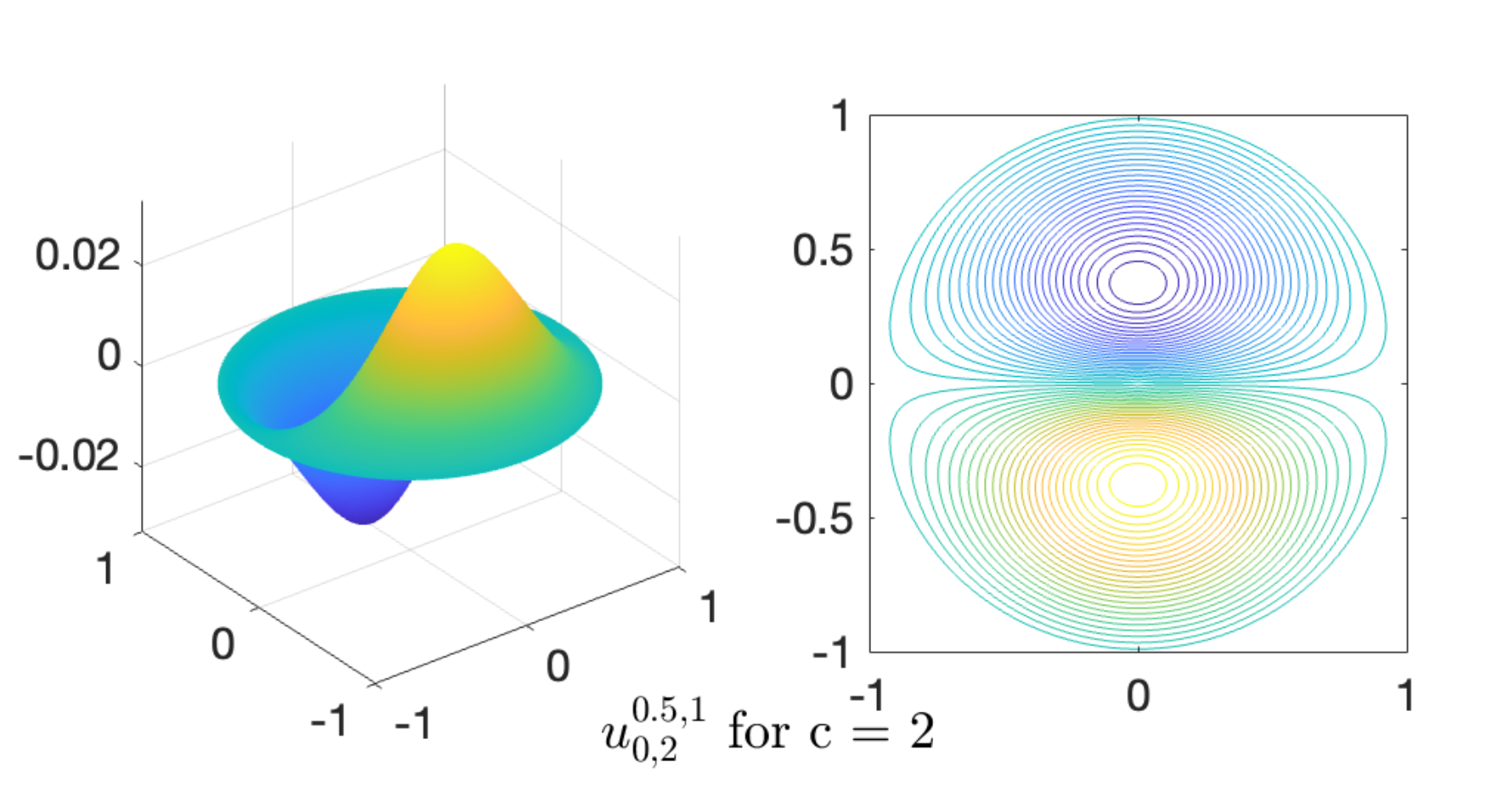}}
	\end{minipage}}
	\subfigure[$(\mu,n,k,\ell)=(0.5,2,0,1)$]{
		\begin{minipage}[t]{0.45\textwidth}
			\centering \rotatebox[origin=cc]{-0}{\includegraphics[width=2.8in,height=1.45in]{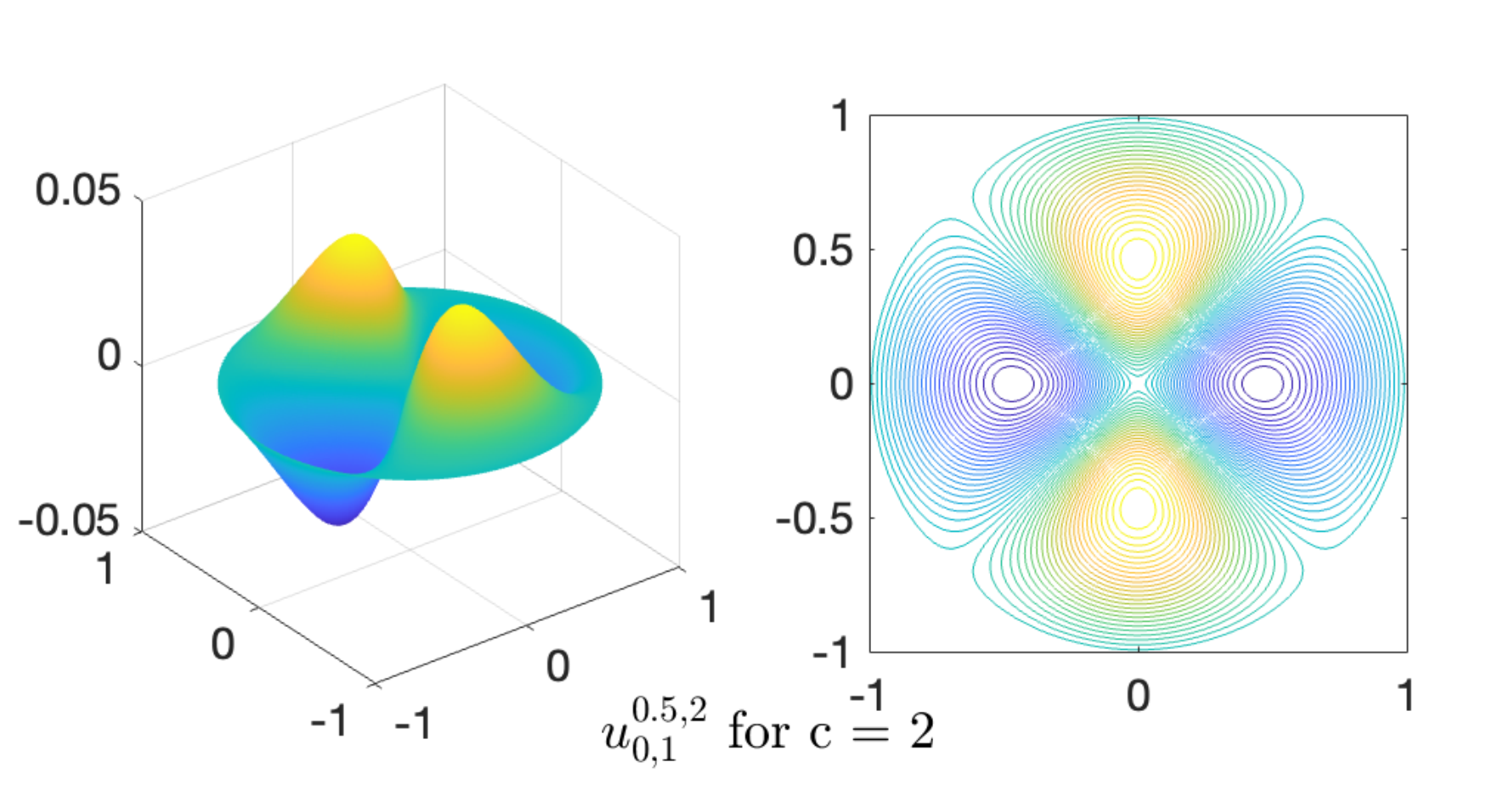}}
	\end{minipage}}%
	\subfigure[$(\mu,n,k,\ell)=(0.5,2,0,2)$]{
		\begin{minipage}[t]{0.45\textwidth}
			\centering \rotatebox[origin=cc]{-0}{\includegraphics[width=2.8in,height=1.45in]{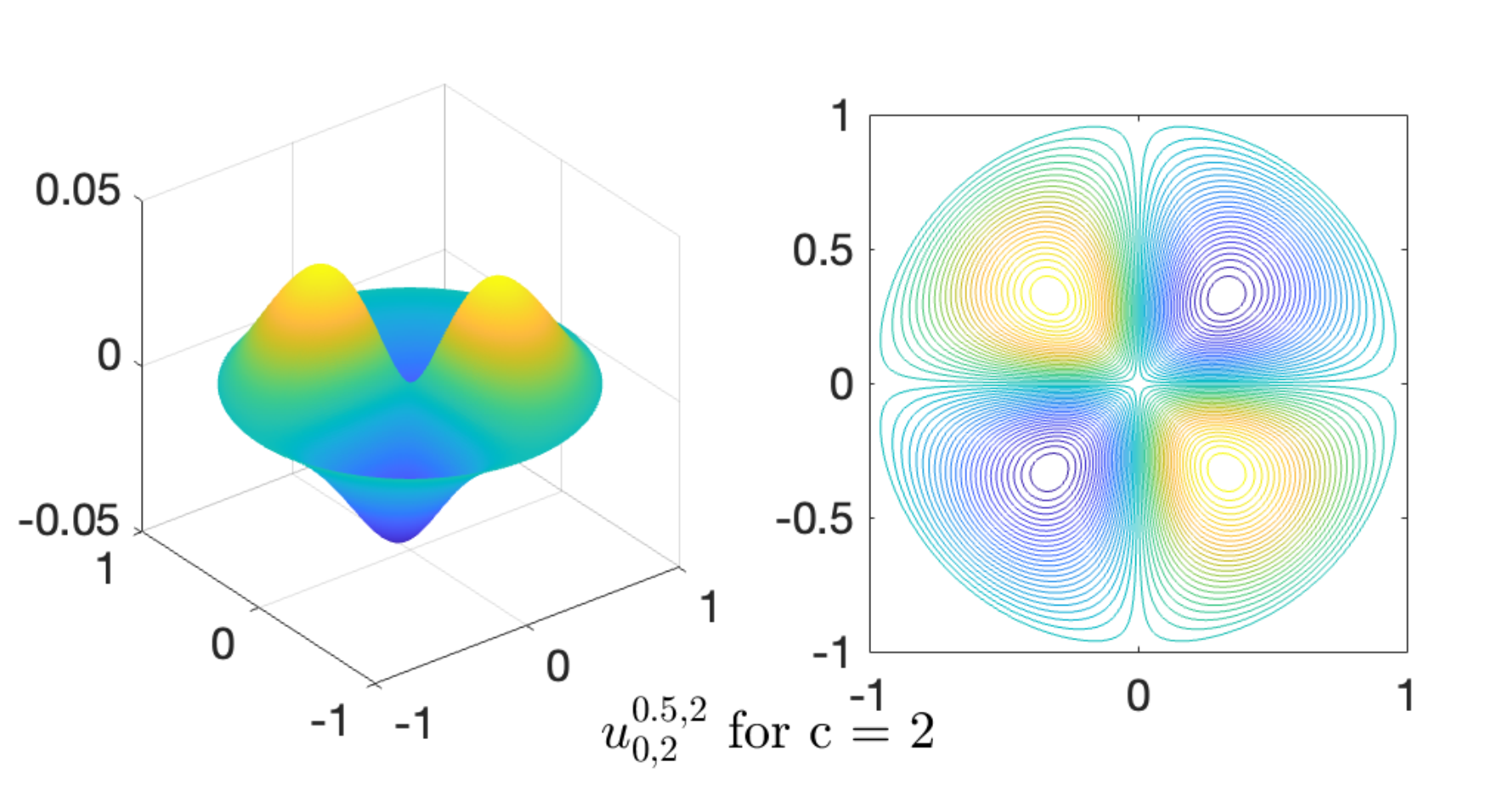}}
	\end{minipage}}%
	\caption {\small Normalized eigenfunctions $\{u_{k,\ell}^{\mu, n}(\bm{x},c)\}$ with  $c=2$, $\mu=1/2$ in $d=2$}\label{GGPs2d-2}
\end{figure}

\begin{figure*}[!htb]
	\centering
	\subfigure[$(\mu,n,k,\ell)=(0.5,0,0,1)$]
	{\includegraphics[width=1.75in,height=1.45in]{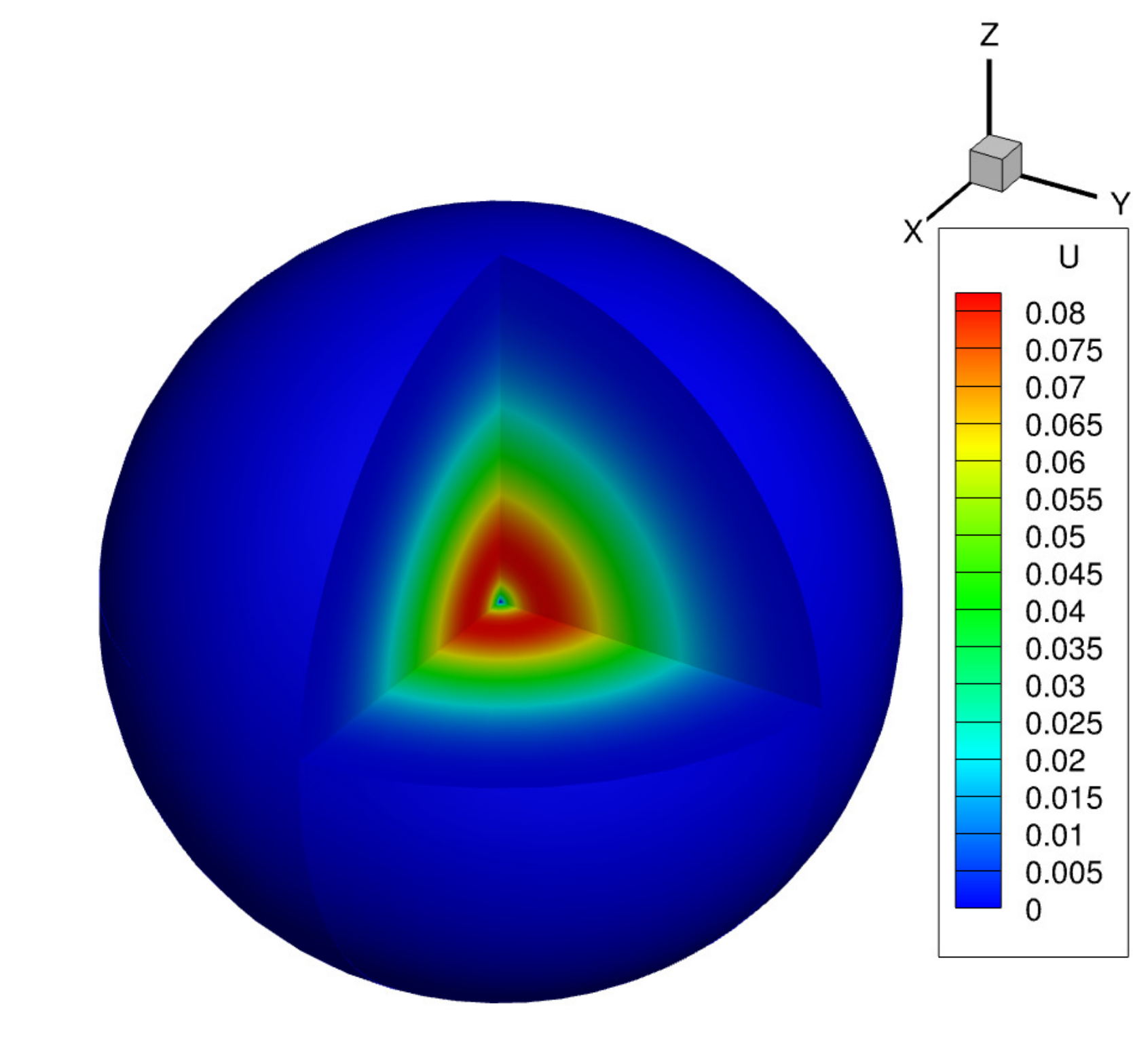}}
	\subfigure[$(\mu,n,k,\ell)=(0.5,0,1,1)$]
	{\includegraphics[width=1.75in,height=1.45in]{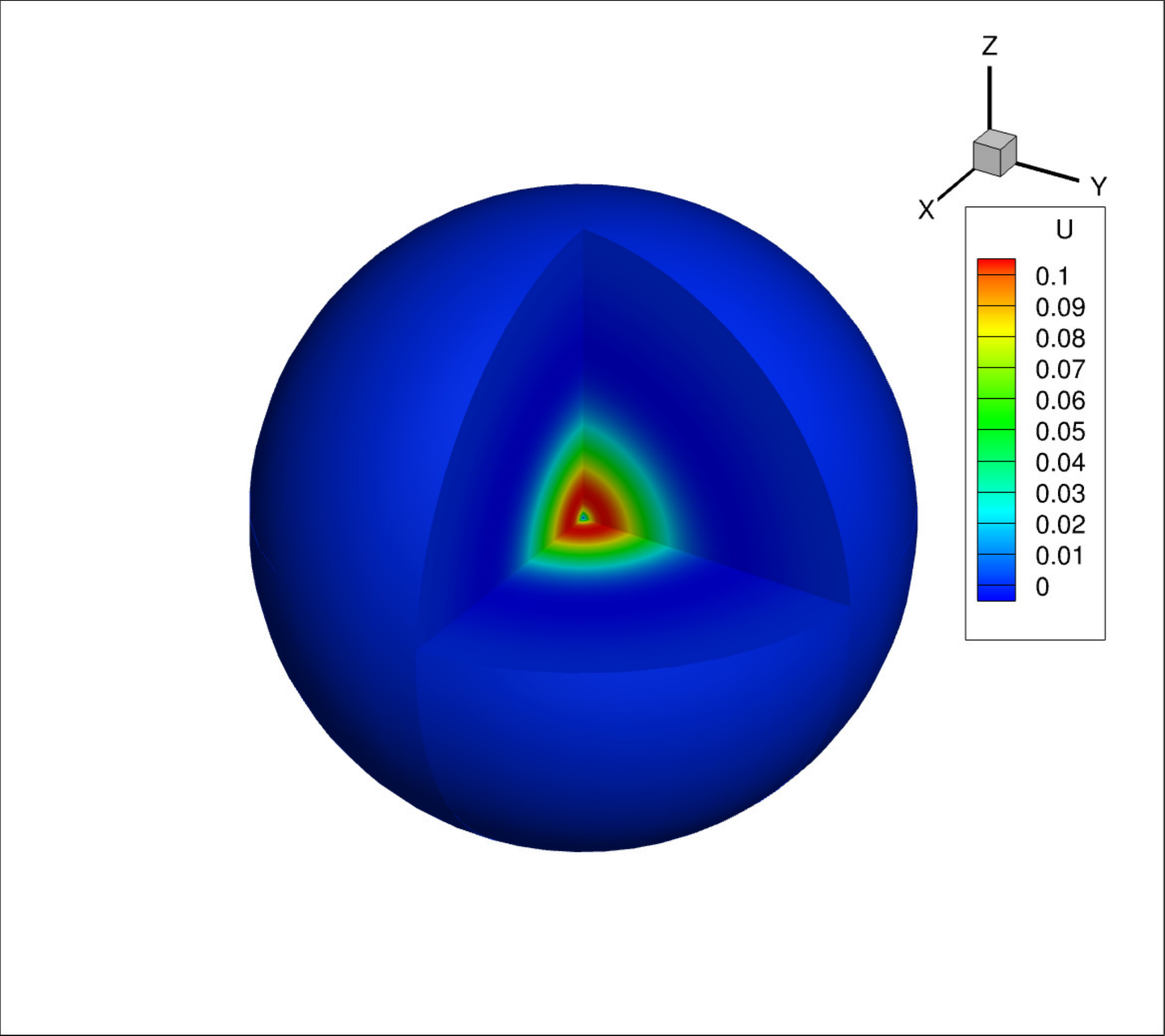}}
	\subfigure[$(\mu,n,k,\ell)=(0.5,0,2,1)$]
	{\includegraphics[width=1.75in,height=1.45in]{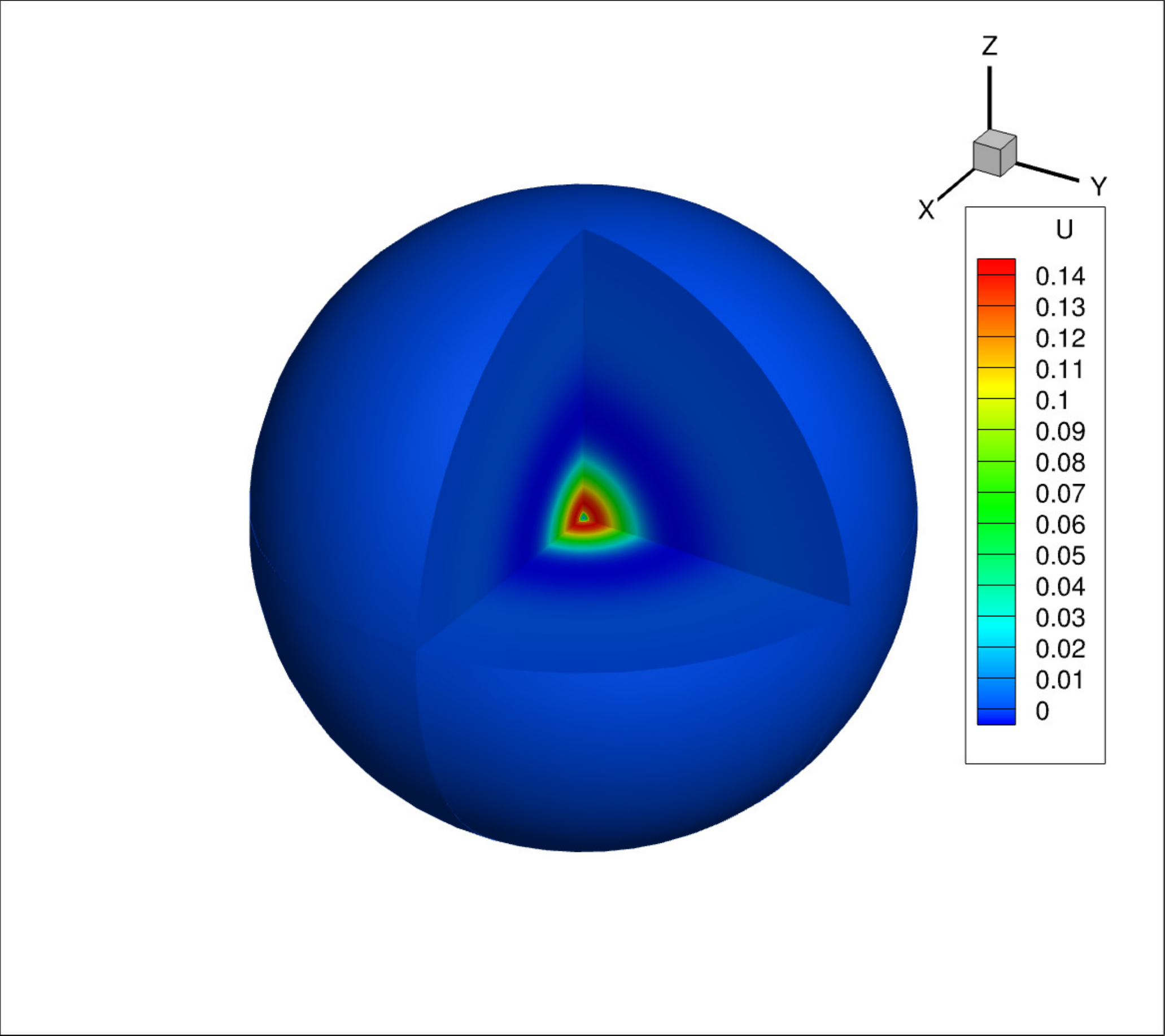}}
	\subfigure[$(\mu,n,k,\ell)=(0.5,2,0,1)$]
	{\includegraphics[width=1.75in,height=1.45in]{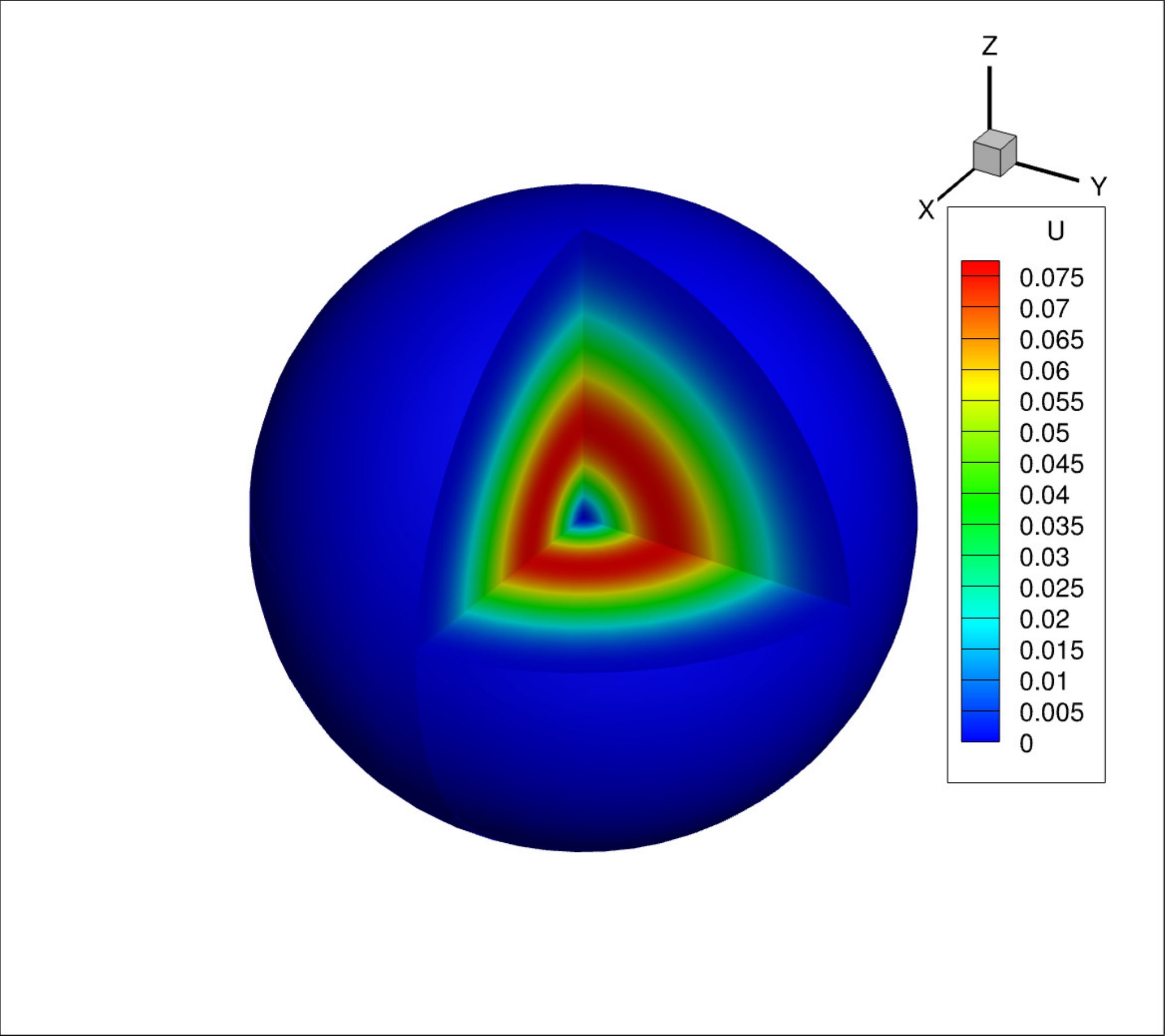}}
	\subfigure[$(\mu,n,k,\ell)=(0.5,2,1,1)$]
	{\includegraphics[width=1.75in,height=1.45in]{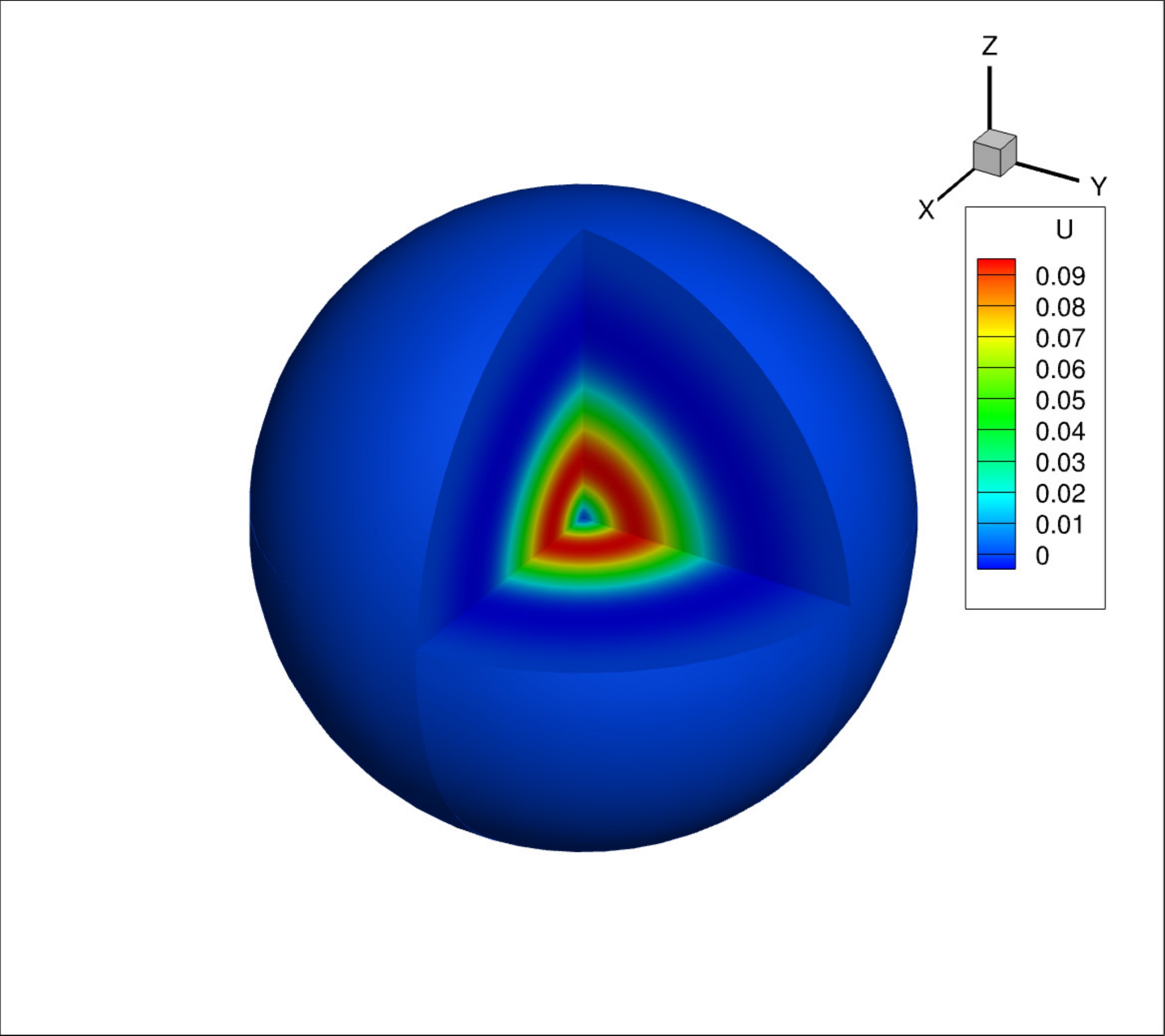}}
	\subfigure[$(\mu,n,k,\ell)=(0.5,2,2,1)$]
	{\includegraphics[width=1.75in,height=1.45in]{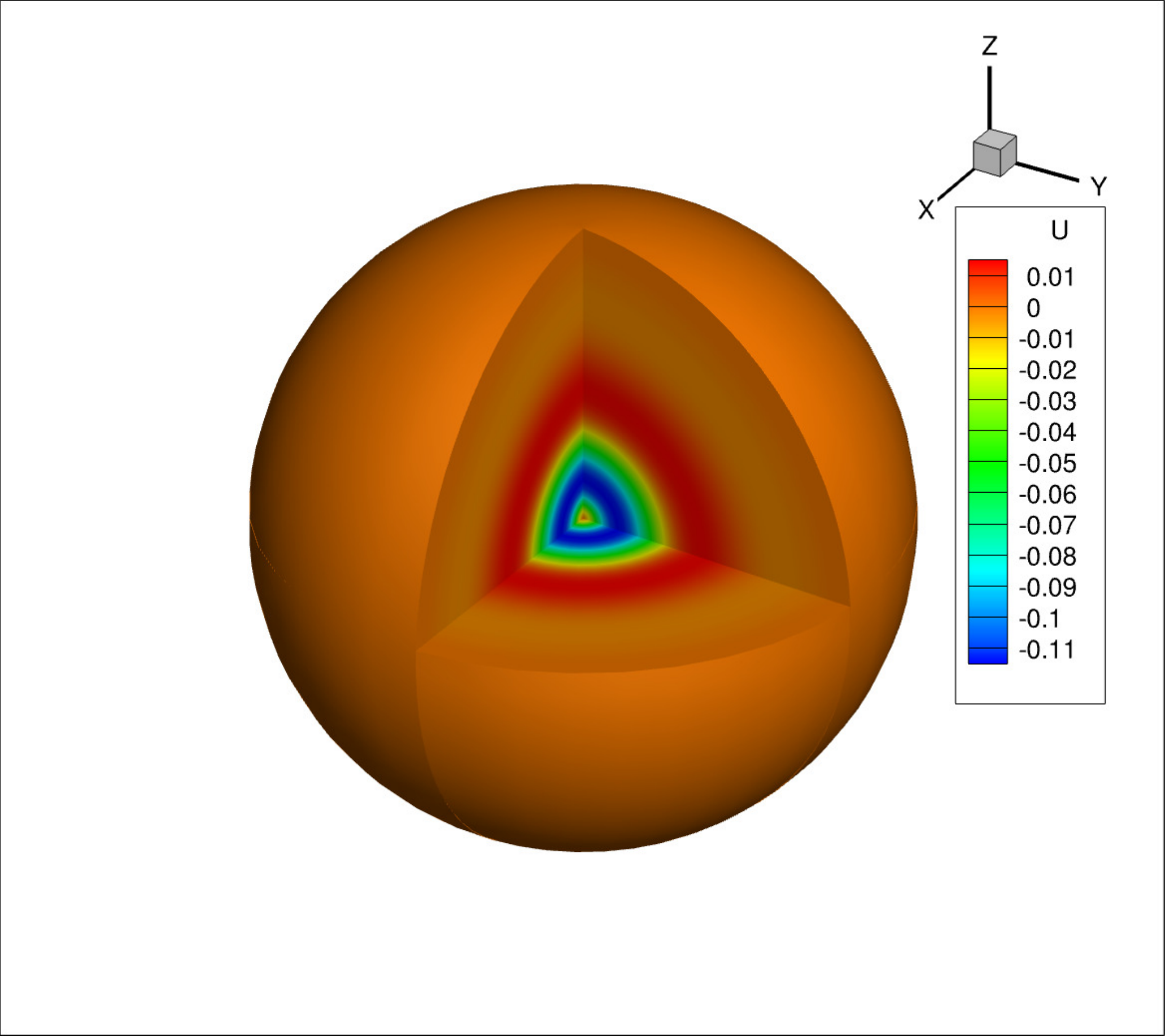}}
	\caption {\small Eigenfunctions $\{u_{k,\ell}^{\mu, n}(\bm{x},c)\}$ with $c=2$, $\mu=1/2$  in $d=3$ }\label{GGPs3d-1}
\end{figure*}

\begin{figure*}[!htb]
	\centering
	\subfigure[$(\mu,n,k,\ell)=(0,0,0,1)$]
	{\includegraphics[width=1.75in,height=1.45in]{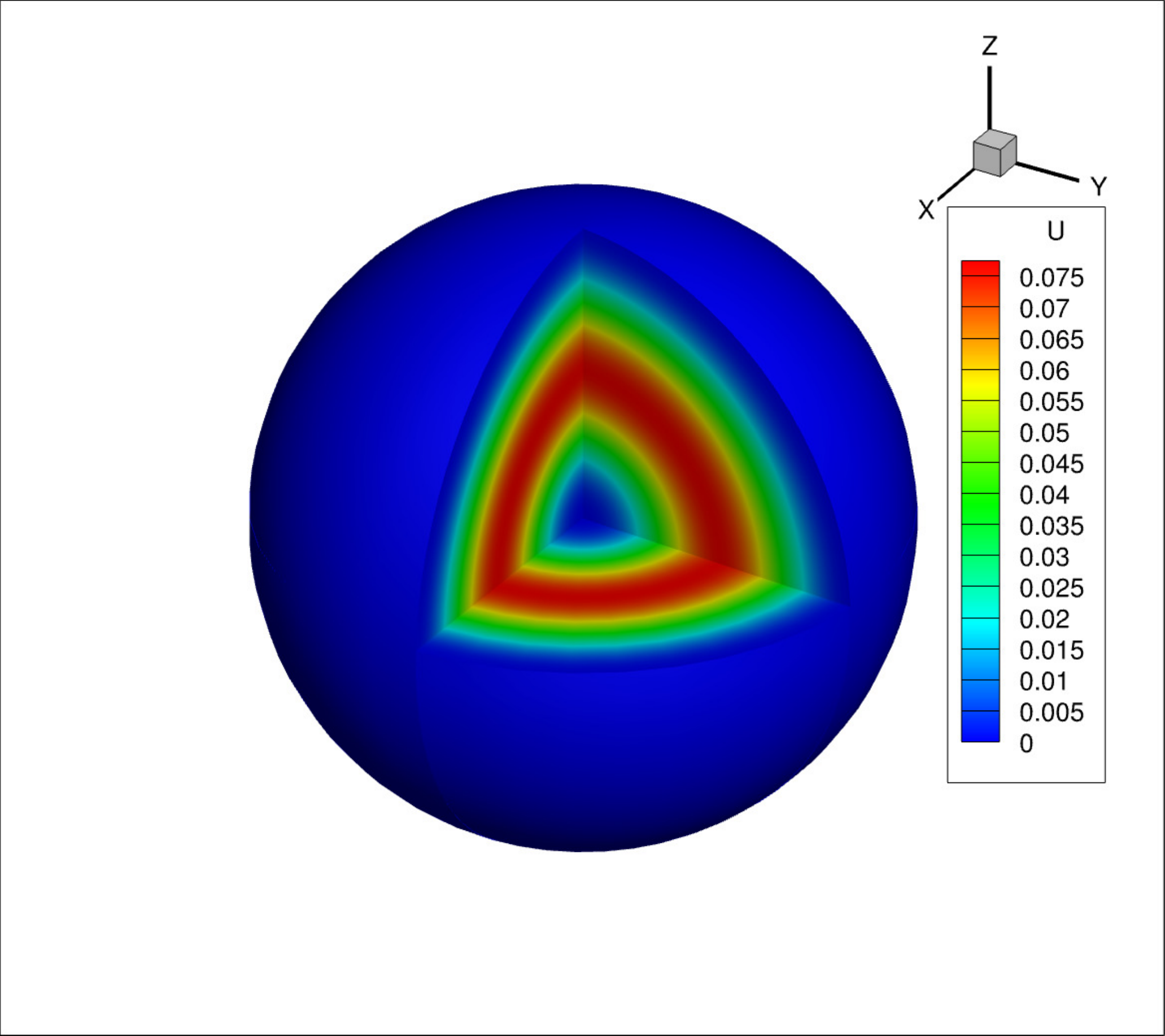}}
	\subfigure[$(\mu,n,k,\ell)=(0,0,1,1)$]
	{\includegraphics[width=1.75in,height=1.45in]{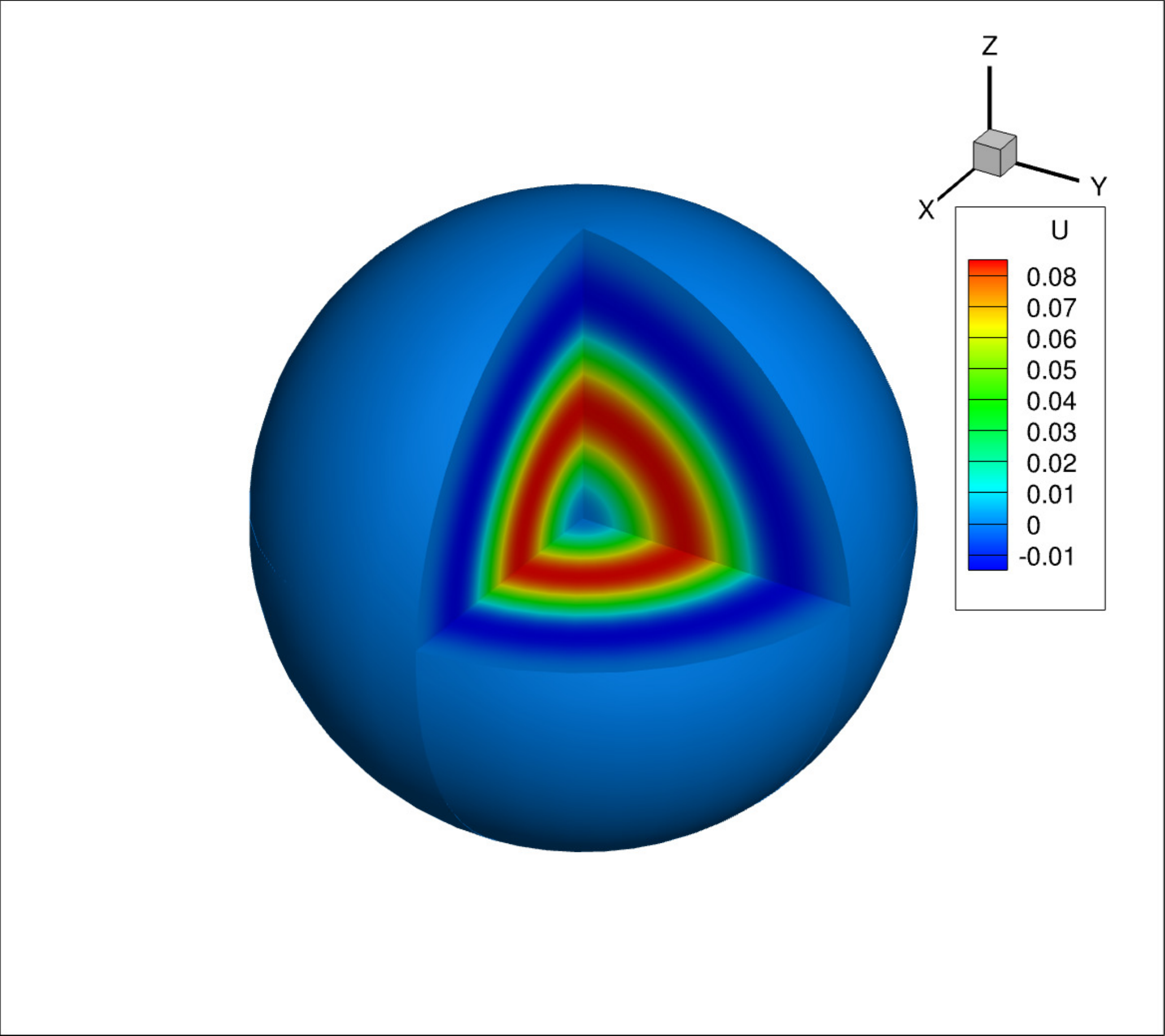}}
	\subfigure[$(\mu,n,k,\ell)=(0,0,2,1)$]
	{\includegraphics[width=1.75in,height=1.45in]{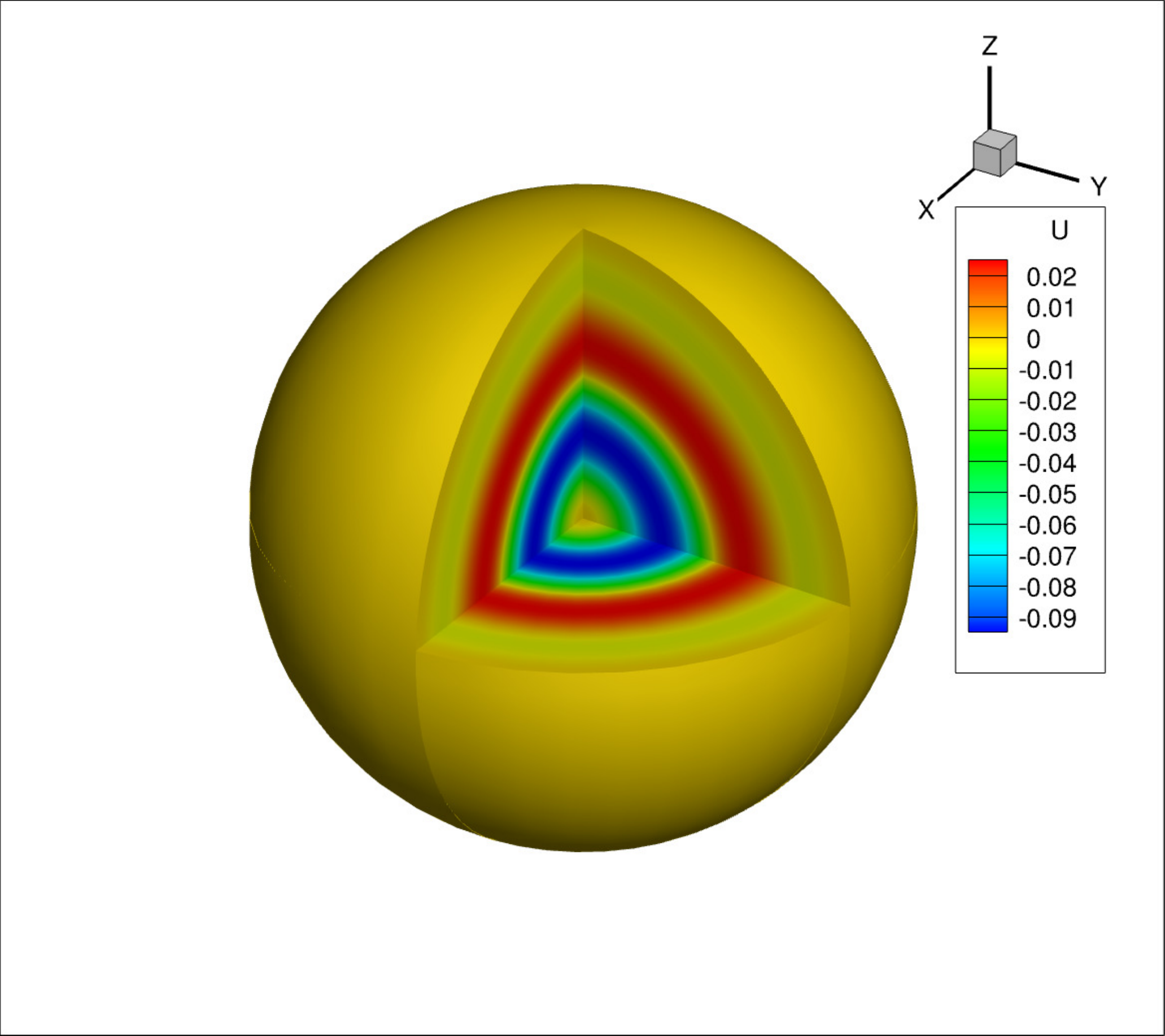}}
	\subfigure[$(\mu,n,k,\ell)=(0,1,0,2)$]
	{\includegraphics[width=1.75in,height=1.45in]{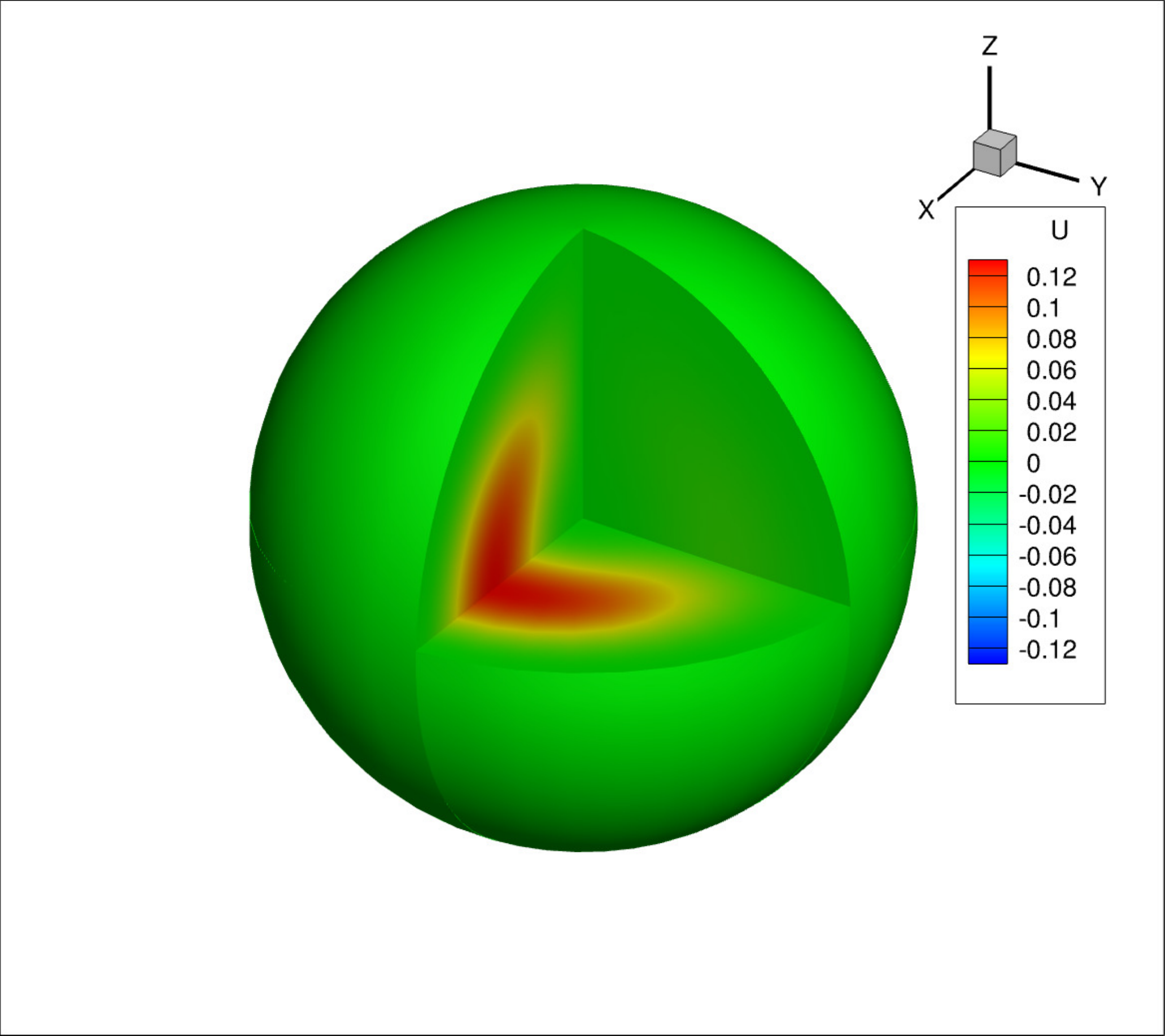}}
	\subfigure[$(\mu,n,k,\ell)=(0,1,1,2)$]
	{\includegraphics[width=1.75in,height=1.45in]{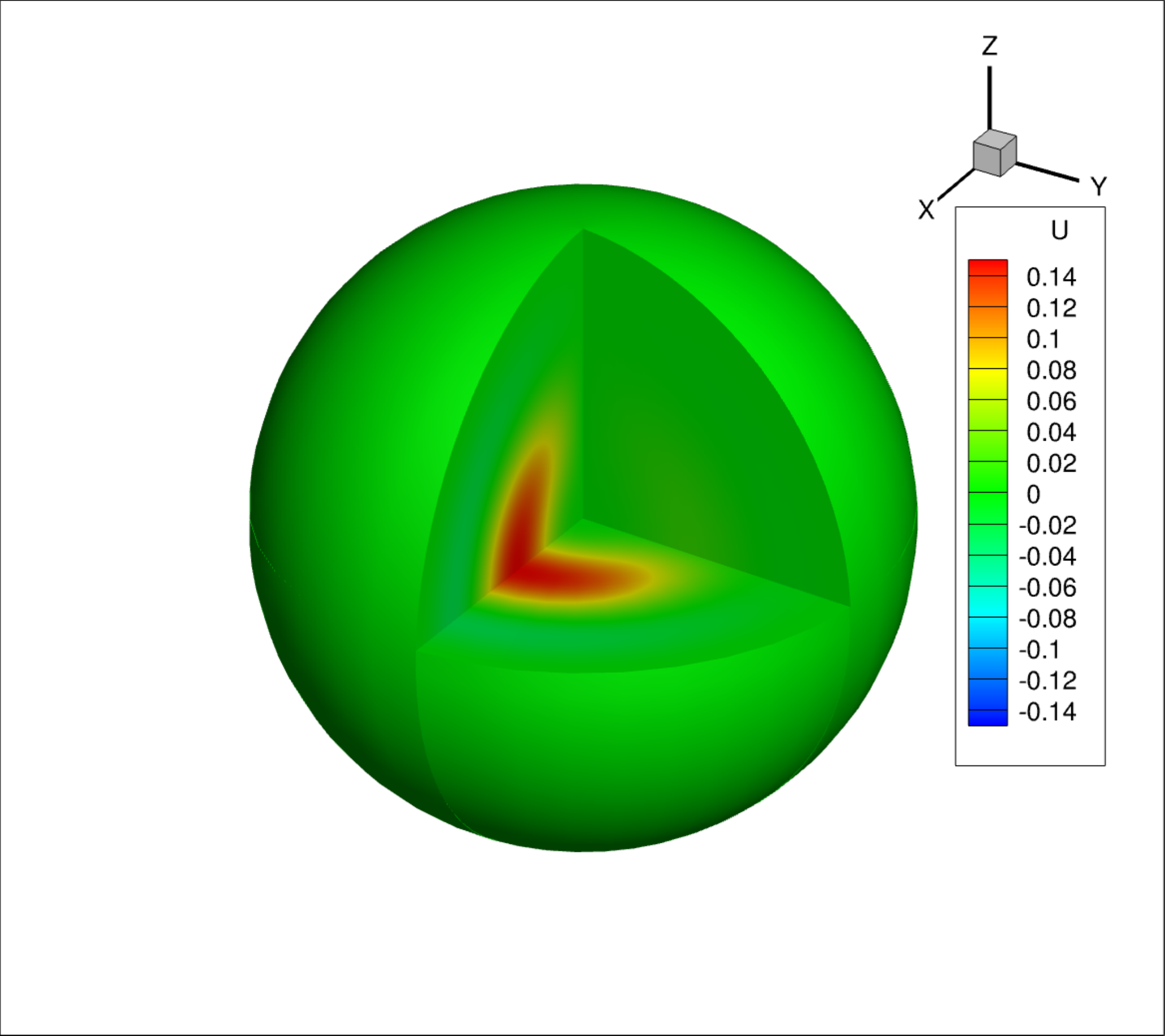}}
	\subfigure[$(\mu,n,k,\ell)=(0,1,2,2)$]
	{\includegraphics[width=1.75in,height=1.45in]{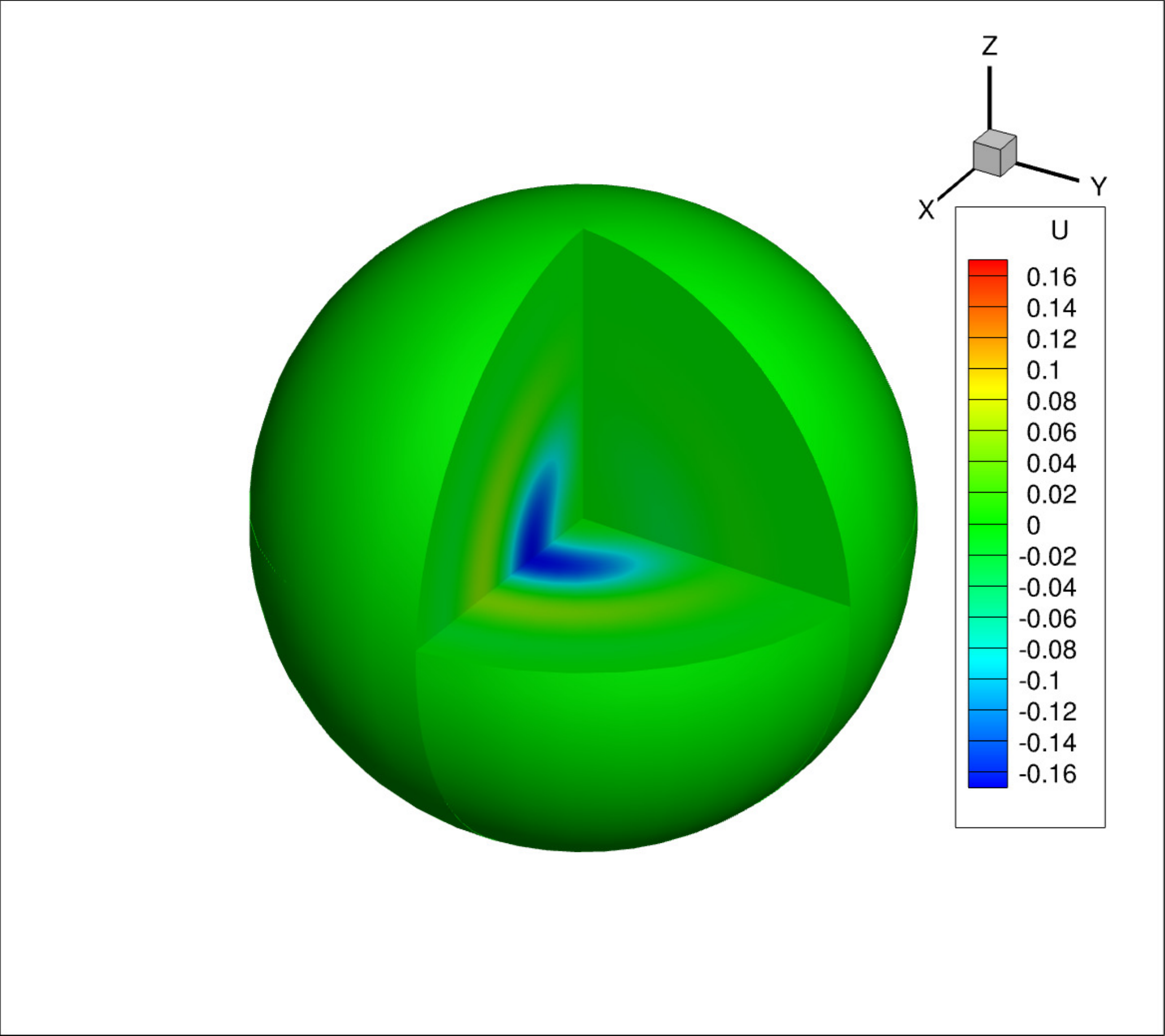}}
	\caption {\small Eigenfunctions $\{u_{k,\ell}^{\mu, n}(\bm{x},c)\}$ with  $c=10$, $\mu=0$ in $d=3$}\label{GGPs3d-2}
\end{figure*}

 In Table \ref{tableeigen}, we list the values of the normalized eigenfunction  in radial direction corresponding to the eigenvalues for several values of $r$ with $\mu=1/2$ in $d=2$.
 One can observe that numerical result obtained by the MBP spectral method match well with the  expression \eqref{anaeigfun}. We denote by $\{u_{k,\ell}^{\mu, n}(\bm{x})\}$ as the normalized eigenfunction of  \eqref{eigensingular}.  In  Figures \ref{GGPs2d-1} and \ref{GGPs2d-2}, we depict  the surfaces and contours of  $\{u_{k,\ell}^{\mu, n}(\bm{x},c)\}$ with different $c$, $k$, $\ell$, $n$  and $\mu =1/2$ in $d =2$.  In Figures \ref{GGPs3d-1} and \ref{GGPs3d-2}, we  intend to visualize  $\{u_{k,\ell}^{\mu, n}(\bm{x})\}$ with different $c$, $k$, $\ell$, $n$  and $\mu =0$ in $d =3$.

\subsection{Schr\"{o}dinger eigenvalue problems with  fractional power potential}
In the sequel, we consider the following Schr\"odinger equation with an inverse and a fractional power potential as follows
	\begin{equation}\label{eigfrac}
		\left\{\begin{array}{ll}
			\Big[-\Delta +\dfrac{c}{\|\bm{x}\|^2}\Big]u(\bm{x}) + z\|\boldsymbol{x}\|^{\frac{2 \nu-2 \eta}{\eta+1}} u(\boldsymbol{x})=\lambda u(\boldsymbol{x}), &\quad  \bm{x} \in \mathbb{B}^{d},\\[4pt]	
			u(\bm{x})=0, &\quad  \bm{x} \in \mathbb{S}^{d-1},	
		\end{array}\right.
	\end{equation}
where $\eta, \nu \in \mathbb{N}_0$. For any given rational number $\frac{q}{p}>-2$ with $p \in \mathbb{N}$ and $q \in \mathbb{Z}$, we can always rewrite it as
$$
\frac{q}{p}=\frac{2 \nu-2 \eta}{\eta+1} \quad \text { with } \quad \eta=2 p-1 \in \mathbb{N}_0, \quad \nu=2 p+q-1 \in \mathbb{N}_0 .
$$
%

Note that the left-hand side of  \eqref{eigensingular} in the spherical-polar coordinates can be reformulated as	
	\[\begin{aligned}
	-\Big(\partial^2_r+ \frac{d-1}{r} \partial_r + \frac{\Delta_0-c}{r^2}\Big)u+z r^{\frac{2 \nu-2 \eta}{\eta+1}}u=-\frac{1}{r^{d-1}}\partial_r[r^{d-1} \partial_r]u - \frac{\Delta_0-c}{r^{2}}u+z r^{\frac{2 \nu-2 \eta}{\eta+1}}u.
\end{aligned}\]
Thus in the radial direction, we have
\begin{equation}
	\begin{aligned} &-\frac{1}{r^{d-1}} \partial_{r}(r^{d-1} \partial_{r} u_{\ell}^{n})+\frac{n(n+d-2)+c}{r^{2}} u_{\ell}^{n}+z r^{\frac{2 \nu-2 \eta}{\eta+1}}u_{\ell}^{n}\\
		=&-\Big(\partial_{r}^{2} u_{\ell}^{n}+\frac{d-1}{r} \partial_{r} u_{\ell}^{n}-\frac{n(n+d-2)+c}{r^{2}} u_{\ell}^{n}\Big) +z r^{\frac{2 \nu-2 \eta}{\eta+1}}u_{\ell}^{n}\\
		=&-\frac{1}{r^{\frac{d}{2}+\theta\beta_n}} \partial_{r}[r^{2\theta\beta_n+1} \partial_{r} (r^{\frac{d}{2}-\theta\beta_n-1}u_{\ell}^{n})]+z r^{\frac{2 \nu-2 \eta}{\eta+1}}u_{\ell}^{n}.
	\end{aligned}
\end{equation}

Using the variable substitution $\rho=r^{2\theta}$, we can write
\begin{equation}\label{anafrac}
	\begin{aligned} &-\frac{1}{r^{d-1}} \partial_{r}(r^{d-1} \partial_{r} u_{\ell}^{n})+\frac{n(n+d-2)+c}{r^{2}} u_{\ell}^{n}+z r^{\frac{2 \nu-2 \eta}{\eta+1}}u_{\ell}^{n}\\
	=&-\frac{4\theta^2}{r^{\frac{d}{2}+\theta\beta_n-2\theta+1}}\partial_{\rho}\big[\rho^{\beta_n+1} \partial_{\rho}\big(r^{\frac{d}{2}-\theta\beta_n-1} u_{\ell}^{n}\big)\big]+z r^{\frac{2 \nu-2 \eta}{\eta+1}}u_{\ell}^{n}.
\end{aligned} 
\end{equation}
If we take  $\theta=\frac{1}{\eta+1}$, then the last term of \eqref{anafrac} becomes $z \rho^qu_{\ell}^{n}$, so we choose the MBP approximation with $\alpha=-1$, $\mu=0$, $\theta=\frac{1}{\eta+1}$ to account for both the accuracy and efficiency.
Accordingly, we introduce the approximation space
\begin{equation}\label{WNKsps}
\mathcal{W}_{N, K}=\operatorname{span}\big\{S_{k, \ell,n,c}^{-1,0,\theta}(\bm{x}), \, (\ell,n)\in \Upsilon_{N}^d,\, 1\leq k\leq K, \, k \in \mathbb{N}_{0}\big\},
\quad \theta=\frac{1}{\eta+1}.
\end{equation}
Then the spectral scheme for \eqref{eigfrac} is to find $\lambda_{N,K} \in \mathbb{R}$ and $u_{N,K}\in \mathcal{W}_{\!N,K} \backslash \{0\}$, such that
\begin{equation}\label{eigfrc1}
   	 \mathcal{B}(u_{N,K},v_{N,K})=\lambda_{N,K}(u_{N,K}, v_{N,K}), \quad \forall\, v_{N,K} \in \mathcal{W}_{\!N,K},
\end{equation}
where $$\mathcal{B}(u_{N,K},v_{N,K})=(\nabla u_{N,K},\nabla v_{N,K})+c(\|\bm x\|^{-2}u_{N,K},v_{N,K})+z(\|\bm x\|^{\frac{2 \nu-2 \eta}{\eta+1}}u_{N,K},v_{N,K}).$$

In implementation, we write
\begin{equation*}
	u_{N,K}(\bm{x})=\sum_{n=0}^N\sum_{\ell=1}^{a_n^d}\sum_{k=1}^{K}\hat{u}_{k,\ell}^nS_{k, \ell,n,c}^{-1,0,\theta}(\bm{x}),
\end{equation*}
and denote
\begin{equation*}
	\begin{array}{l}
		\boldsymbol{u}=\big( \hat{\boldsymbol{u}}_{1}^{0}, \hat{\boldsymbol{u}}_{2}^{0}, \cdots, \hat{\boldsymbol{u}}_{a_{0}^{d}}^{0}, \cdots, \hat{\boldsymbol{u}}_{1}^{N}, \hat{\boldsymbol{u}}_{2}^{N}, \cdots, \hat{\boldsymbol{u}}_{a_{N}^{d}}^{N}\big)^{t}, \quad \hat{\boldsymbol{u}}_{\ell}^{n}=\big(\hat{u}_{1, \ell}^{n}, \hat{u}_{2, \ell}^{n}, \ldots, \hat{u}^n_{K, \ell}\big)^{t}.
	\end{array}
\end{equation*}
Corresponding to this ordering, we denote the stiffness and mass matrices by $\bm{S}$ and $\bm{M}$, respectively.  The algebraic eigen-system  of  \eqref{eigfrc1} reads
\[\bm{S}\boldsymbol{u}=\lambda_N\bm{M}\boldsymbol{u}.\]
Moreover, we can explicitly evaluate their entries. For fixed $k,j \in\mathbb{N}_0$, $(\ell,n),(\iota,m)\in \Upsilon_{\infty}^d$, we derive from \eqref{ballint} that
\begin{equation*}
	\begin{split}
		&\big(S_{k,  \ell,n,c}^{-1,0,\theta}, S_{j, \iota,m,c}^{-1,0,\theta}\big)=\int_{\mathbb{S}^{d-1}} Y_{\ell}^{n}(\bm{\hat{x}}) Y_{\iota}^{m}(\bm{\hat{x}}) \mathrm{d} \sigma(\bm{\hat{x}}) \\
		&\quad \times  \int_{0}^{1} P_{k}^{(-1,\beta_n)}(2 r^{\frac{2}{\eta+1}}-1) P_{j}^{(-1,\beta_n)}(2 r^{\frac{2}{\eta+1}}-1) r^{\frac{2\beta_n}{\eta+1}+1} \mathrm{d} r \\
		&=\delta_{k j}  \delta_{n m}  \int_{0}^{1} P_{k}^{(-1,\beta_n)}(2 r^{\frac{2}{\eta+1}}-1) P_{j}^{(-1,\beta_n)}(2 r^{\frac{2}{\eta+1}}-1) r^{\frac{2\beta_n}{\eta+1}+1} \mathrm{d} r \\
		&=\delta_{k j}  \delta_{n m} \dfrac{\eta+1}{2^{\beta_n+\eta+2}} \int_{-1}^{1} P_{k}^{(-1,\beta_n)}(\rho) P_{j}^{(-1,\beta_n)}(\rho)(\rho+1)^{\beta_n+\eta} \mathrm{d} \rho\qquad (\text{note: } \rho=2r^{2\theta}-1) \\
		&=\delta_{k j}  \delta_{n m} \dfrac{\eta+1}{2^{\beta_n+\eta+4}}\dfrac{k+\beta_n}{k} \dfrac{j+\beta_n}{j} \int_{-1}^{1} P_{k}^{(1,\beta_n)}(\rho) P_{j}^{(1,\beta_{n})}(\rho)(1-\rho)^2(1+\rho)^{\beta_n+\eta} \mathrm{d} \rho,
	\end{split}	
\end{equation*}
and
	\begin{align}
		&\big({\|\boldsymbol{x}\|^{\frac{2 \nu-2 \eta}{\eta+1}} }S_{k, \ell,n,c}^{-1,0,\theta}, S_{j, \iota,m,c}^{-1,0,\theta}\big)=\int_{\mathbb{S}^{d-1}} Y_{\ell}^{n}(\bm{\hat{x}}) Y_{\iota}^{m}(\bm{\hat{x}}) \mathrm{d} \sigma(\bm{\hat{x}}) \nonumber\\
		&\quad \times  \int_{0}^{1} P_{k}^{(-1,\beta_n)}(2 r^{\frac{2}{\eta+1}}-1) P_{j}^{(-1,\beta_{m})}(2 r^{\frac{2}{\eta+1}}-1) r^{\frac{2\beta_n+2 \nu-2 \eta}{\eta+1}+1} \mathrm{d} r \nonumber\\
		&=\delta_{k j}  \delta_{n m}  \int_{0}^{1} P_{k}^{(-1,\beta_n)}(2 r^{\frac{2}{\eta+1}}-1) P_{j}^{(-1,\beta_n)}(2 r^{\frac{2}{\eta+1}}-1) r^{\frac{2\beta_n+2 \nu-2 \eta}{\eta+1}+1} \mathrm{d} r \nonumber\\
		&=\delta_{k j}  \delta_{n m} \dfrac{\eta+1}{2^{\beta_n+\nu+2}} \int_{-1}^{1} P_{k}^{(-1,\beta_n)}(\rho) P_{j}^{(-1,\beta_n)}(\rho)(\rho+1)^{\beta_n+\nu} \mathrm{d} \rho \nonumber\\
		&=\delta_{k j}  \delta_{n m} \dfrac{\eta+1}{2^{\beta_n+\nu+4}}\dfrac{k+\beta_n}{k} \dfrac{j+\beta_n}{j} \int_{-1}^{1} P_{k}^{(1,\beta_n)}(\rho) P_{j}^{(1,\beta_n)}(\rho)(1-\rho)^2(1+\rho)^{\beta_n+\nu} \mathrm{d} \rho.\nonumber
	\end{align}
Using  Lemma \ref{lemmuntzorth} and  \eqref{GUPeq50}, we find that the stiffness matrix $\bm{S}$ is a banded matrix with a bandwidth $2\nu+1$ and the mass matrix $\bm{M}$ is also banded with a bandwidth $2\eta+1$.

\begin{figure}[!htb]
	\centering
	\subfigure[$d=1,\eta=3,\nu=2,c=2,z=-3$]
	{\includegraphics[width=7.2cm,height=5.3cm]{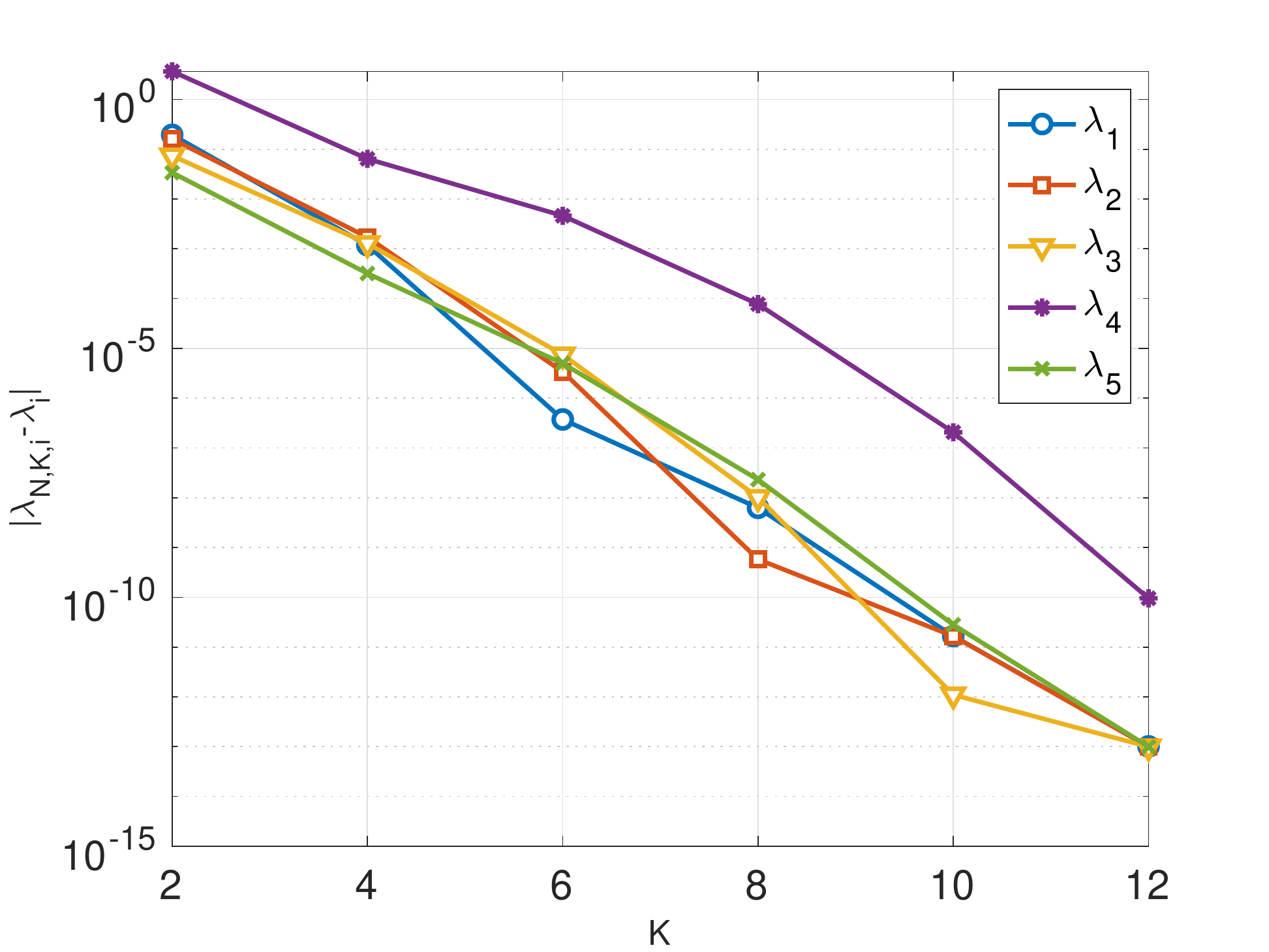}}
	\subfigure[$d=2,\eta=1,\nu=4,c=5,z=3$]
	{\includegraphics[width=7.2cm,height=5.3cm]{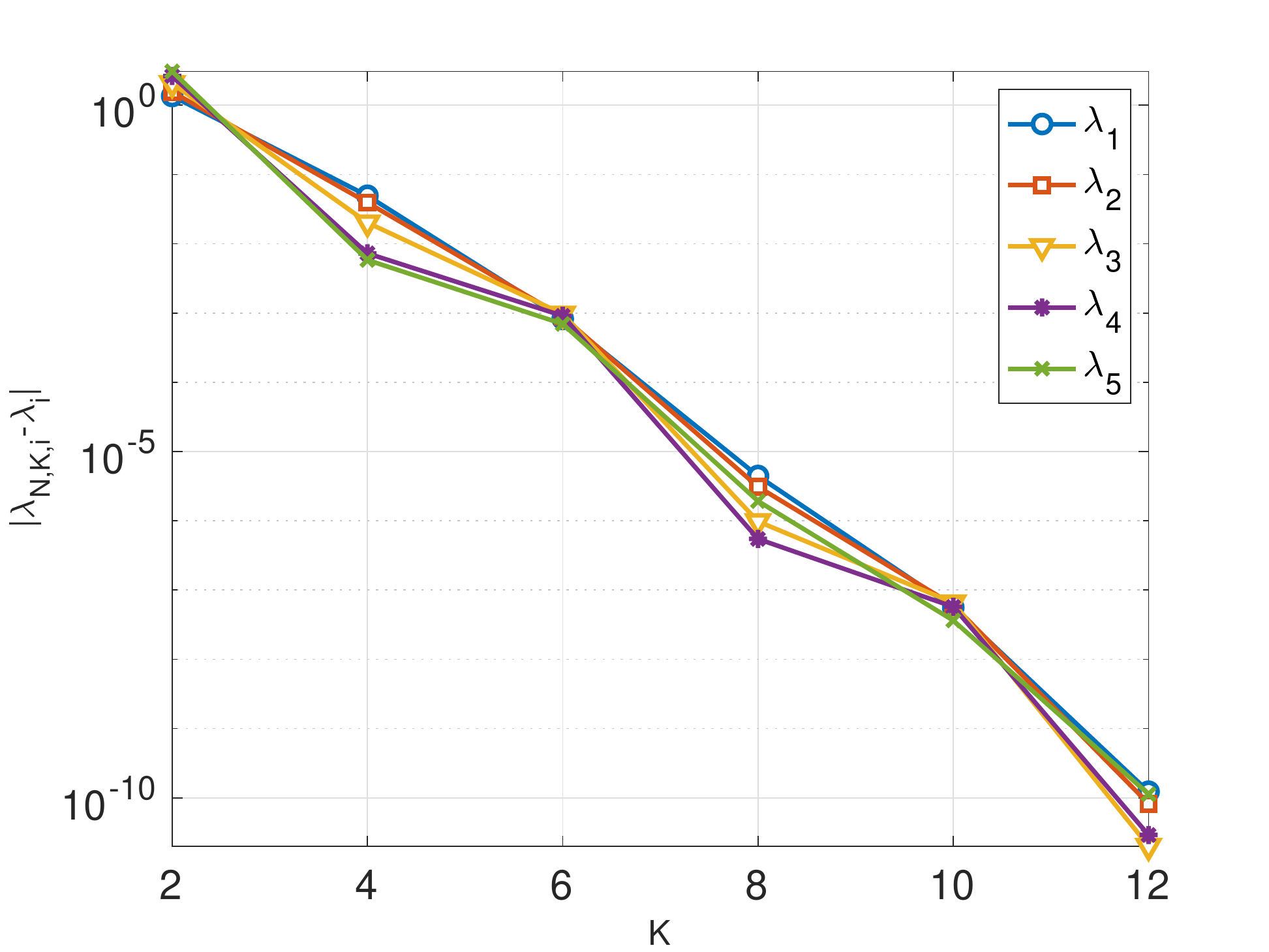}}
	\subfigure[$d=3,\eta=1,\nu=2,c=10,z=1$]
	{\includegraphics[width=7.2cm,height=5.3cm]{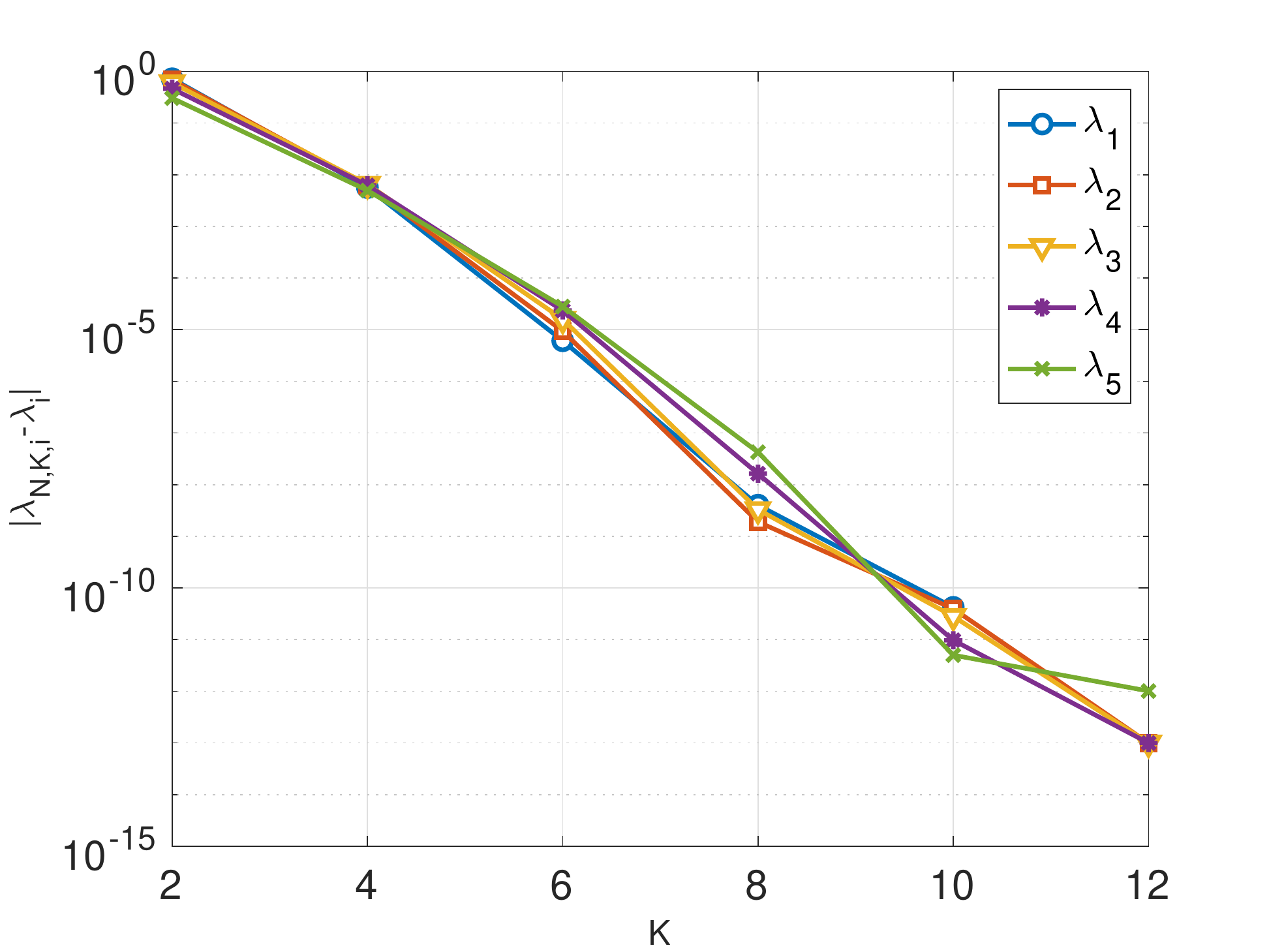}}
	\subfigure[$d=4,\eta=3,\nu=5,c=0.1,z=1$]
	{\includegraphics[width=7.2cm,height=5.3cm]{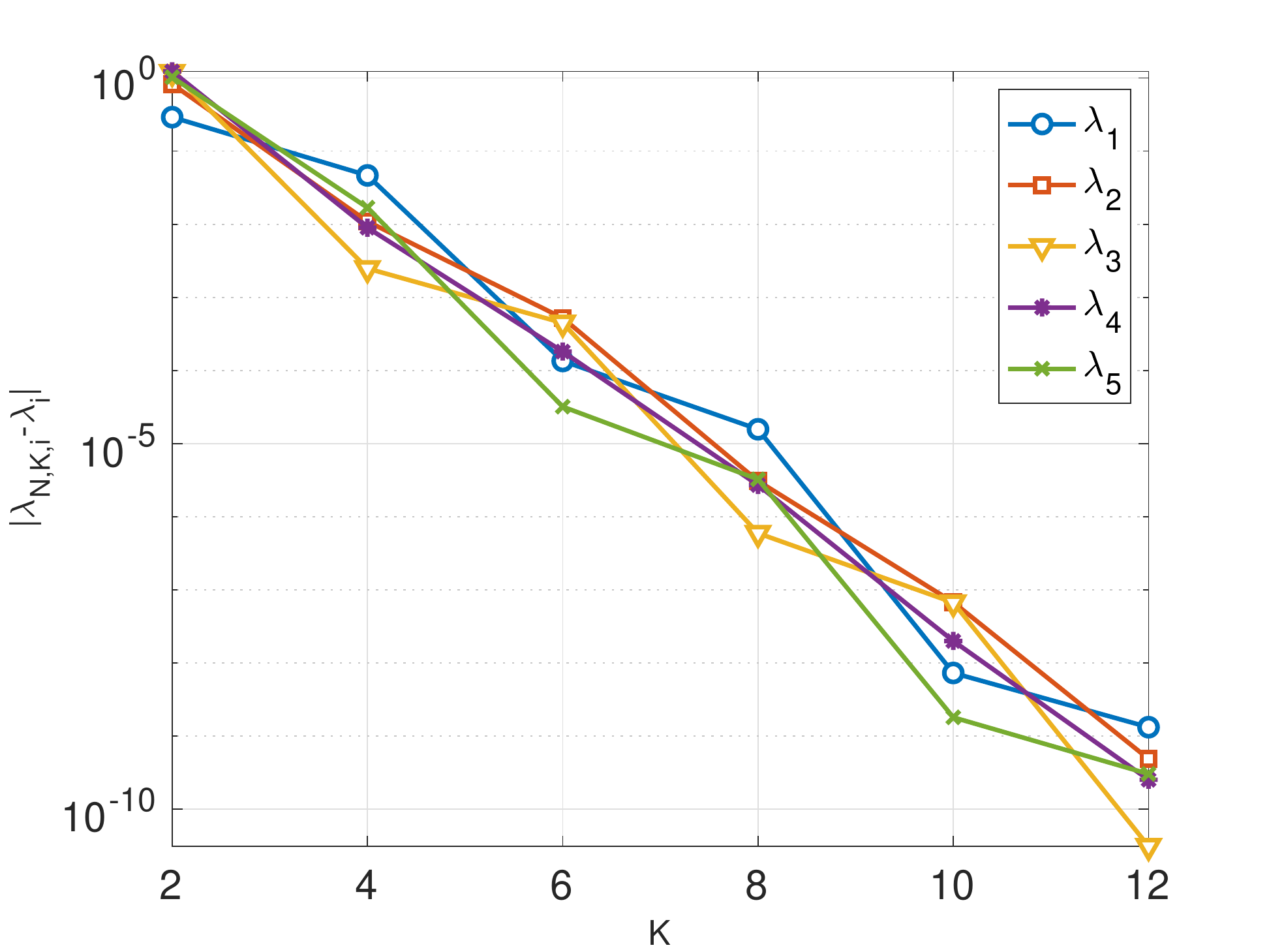}}
	\caption {\small The errors of the smallest $5$ eigenvalues versus $K$ for solving \eqref{eigfrac} with $N=10$}\label{frac_error}
\end{figure}

In the numerical tests, we fix $N=10$, and choose different  values of $c, z, \eta, \nu, d$. in Figure \ref{frac_error}, we depict the numerical errors between the first several  eigenvalues by MBP spectral methods and the reference eigenvalues (obtained by the scheme with large $N$ and $K$). Exponential orders of convergence are clearly observed in all cases, which demonstrate the effectiveness of the MBP spectral method.

\section{Concluding remarks}

In this paper, we introduced a new family of orthogonal  M\"{u}ntz ball polynomials and presented various appealing properties. 
We then developed  efficient and accurate  MBP  spectral-Galerkin methods for a class of  degenerating eigenvalue problems  with singular potentials and Schr\"{o}dinger eigenvalue problems with fractional power potentials. 
 The proposed approximation tools should have a much wider capability for numerical solutions of  singular PDEs.

\medskip 

\noindent{\bf Declarations}

\begin{itemize}
  \item {\bf Availability of data and materials:}  The datasets generated during and/or analysed during the current study are available from the corresponding author on reasonable request. 
  \item {\bf Authors' contributions:} All authors contributed to this study. The computations and the first draft were prepared by the first author. All authors read and approved the final manuscript.
\item {\bf Conflict of interest statement:}   We have no conflicts of interest to disclose.
\end{itemize}

\begin{appendix}
\renewcommand{\theequation}{A.\arabic{equation}}	


\section{Proof of Lemma \ref{lemmuntzorth}}\label{prolemmamuntzorth}
	First, for the Laplace-Beltrami operator $\Delta_{0}$, it holds that (cf.\! \cite[pp. 16, 26]{Dai2013approximation})
	\begin{equation}\label{weaknabla0}
		\quad\left(\nabla_{\!0} u, \nabla_{\!0} v\right)_{\mathbb{S}^{d-1}}=-\left(\Delta_{0} u, v\right)_{\mathbb{S}^{d-1}}, \quad u \in H^{2}\left(\mathbb{S}^{d-1}\right), v \in H^{1}\left(\mathbb{S}^{d-1}\right).
	\end{equation}
	We next prove that 
	\begin{align}
			&\int_{0}^{1} \big(r^{d / 2-\theta\beta_n-1+\mu} u\big)' \big(r^{d / 2-\theta\beta_n-1+\mu} v\big)'  r^{2\theta\beta_n+1} \mathrm{d} r \nonumber\\
			&=\dint_{0}^{1}\Big[r^{2\mu+d-1}  u'  v'+\Big((\theta\beta_n)^2-\big(\frac{d}{2}-(1-\mu)\big)^{2}\Big)r^{2\mu+d-3} u v\Big] \mathrm{d} r\nonumber \\
			&\quad+\big(d / 2-\theta\beta_n-1+\mu)\big)\big(u(1) v(1)-\delta_{2\mu+d, 2} u(0) v(0)\big).\nonumber
		\end{align}
	We derive from direct calculation and integration by parts that
	\begin{eqnarray}\label{intertrans}
		&\dint_{0}^{1} \big(r^{d / 2-\theta\beta_n-1+\mu} u\big)' \big(r^{d / 2-\theta\beta_n-1+\mu} v\big)'  r^{2\theta\beta_n+1} \mathrm{d} r  \nonumber
		\\&=\dint_{0}^{1}\Big[r^{2\mu+d-1}  u' v'+r^{2\mu+d-3}\big(d / 2-\theta\beta_n-1+\mu\big)^{2} u v \nonumber\\&+\big(d / 2-\theta\beta_n-1+\mu\big)r^{2\mu+d-2} (u v)'\Big] \mathrm{d} r \nonumber\\
		&=\dint_{0}^{1}\Big[r^{2\mu+d-1}  u'  v'+r^{2\mu+d-3}\big(d / 2-\theta\beta_n-1+\mu\big)^{2} u v\\
		&-(2\mu+d-2)\big(d / 2-\theta\beta_n-1+\mu)\big) r^{2\mu+d-3} u v\Big] \mathrm{d} r \nonumber \\
		&\quad +\big(d / 2-\theta\beta_n-1+\mu)\big)\left[u(1) v(1)-0^{2\mu+d-2} u(0) v(0)\right] \nonumber\\
		&\hspace{50pt}=\dint_{0}^{1}\Big[r^{2\mu+d-1} u'v'+\Big((\theta\beta_n)^2-\big(\frac{d}{2}-(1-\mu)\big)^{2}\Big)r^{2\mu+d-3} u v\Big] \mathrm{d} r\nonumber \\
		&\quad+\big(d / 2-\theta\beta_n-1+\mu)\big)\big(u(1) v(1)-\delta_{2\mu+d, 2} u(0) v(0)\big).\nonumber
	\end{eqnarray}
	From  \eqref{nabla0}, we know that 
	$$\bm{\hat{x}}\cdot \nabla_{0}=\bm{x}\cdot \nabla-r\partial_{r}=r\bm{\hat{x}}\cdot\nabla-r\partial_{r}=0,$$ so
\begin{equation}\label{weaknabla}
	\quad(\nabla u, \nabla v)_{\Omega}=\left(\nabla_{0} u, \nabla_{0} v\right)_{r^{-2}, \Omega}+\left(\partial_{r} u, \partial_{r} v\right)_{\Omega}, \quad\forall\, u, v \in H^{1}(\Omega).
\end{equation}
Then from \eqref{weaknabla0}, \eqref{weaknabla} and \eqref{GGPeq38}, we obtain that	
	\[\begin{aligned}
		&\big(\nabla \mathcal{S}_{k, \ell,n,c}^{-1,\mu,\theta}, \nabla \mathcal{S}_{j, \iota,m,c}^{-1,\mu,\theta}\big)_{r^{2\mu}}+c\big(\mathcal{S}_{k, \ell,n,c}^{-1,\mu,\theta}, \nabla \mathcal{S}_{j, \iota,m,c}^{-1,\mu,\theta}\big)_{r^{2\mu-2}}  \\
		&=\big(\partial_{r} \mathcal{S}_{k, \ell,n,c}^{-1,\mu,\theta}, \partial_{r} \mathcal{S}_{j, \iota,m,c}^{-1,\mu,\theta}\big)_{r^{2\mu}}+\big(\nabla_{\!0} \mathcal{S}_{k, \ell,n,c}^{-1,\mu,\theta}, \nabla_{\!0} \mathcal{S}_{j, \iota,m,c}^{-1,\mu,\theta}\big)_{r^{2\mu-2}}+c\big(\mathcal{S}_{k, \ell,n,c}^{-1,\mu,\theta}, \mathcal{S}_{j, \iota,m,c}^{-1,\mu,\theta}\big)_{r^{2\mu-2}} \\
		&=\big(\partial_{r} \mathcal{S}_{k, \ell,n,c}^{-1,\mu,\theta}, \partial_{r} \mathcal{S}_{j, \iota,m,c}^{-1,\mu,\theta}\big)_{r^{2\mu}}-\big(\Delta_{0} \mathcal{S}_{k, \ell,n,c}^{-1,\mu,\theta}, \mathcal{S}_{j, \iota,m,c}^{-1,\mu,\theta}\big)_{r^{2\mu-2}}+c\big(\mathcal{S}_{k, \ell,n,c}^{-1,\mu,\theta}, \mathcal{S}_{j, \iota,m,c}^{-1,\mu,\theta}\big)_{r^{2\mu-2}} \\
		&=\big(\partial_{r} \mathcal{S}_{k, \ell,n,c}^{-1,\mu,\theta}, \partial_{r} \mathcal{S}_{j, \iota,m,c}^{-1,\mu,\theta}\big)_{r^{2\mu}}+\left(c+n(n+d-2)\right)\big(\mathcal{S}_{k, \ell,n,c}^{-1,\mu,\theta}, \mathcal{S}_{j, \iota,m,c}^{-1,\mu,\theta}\big)_{r^{2\mu-2}} .
	\end{aligned}\]
	For notational convenience,
 denote $$q_{k, n}(r)= r^{\theta\beta_n+1-d/2-\mu} P_{k}^{(-1, \beta_{n})}(2 r^{2\theta}-1) .$$
We further obtain from \eqref{intertrans}  that	
	\begin{align}\label{stiffmatrix1}
		&\big(\nabla \mathcal{S}_{k, \ell,n,c}^{-1,\mu,\theta}, \nabla \mathcal{S}_{j, \iota,m,c}^{-1,\mu,\theta}\big)_{r^{2\mu}}+c\big(\mathcal{S}_{k, \ell,n,c}^{-1,\mu,\theta}, \nabla \mathcal{S}_{j, \iota,m,c}^{-1,\mu,\theta}\big)_{r^{2\mu-2}} \nonumber\\
		&= \int_{\mathbb{S}^{d-1}} Y_{\ell}^{n}(\bm{\hat{x}}) Y_{\iota}^{m}(\bm{\hat{x}}) d \sigma(\bm{\hat{x}}) \int_{0}^{1}\left[r^{2\mu+d-1}  q_{k, n}'  q_{j, n}'+\left(c+n(n+d-2)\right) r^{2\mu+d-3} q_{k, n}  q_{j, n}\right] \mathrm{d} r \nonumber\\
		&= \delta_{nm} \delta_{\ell\iota} \int_{0}^{1}  \big(r^{d / 2-\theta\beta_n-1+\mu} q_{k, n} \big)' \big(r^{d / 2-\theta\beta_n-1+\mu} q_{j, n} \big)'  r^{2\theta\beta_n+1} \mathrm{d} r \nonumber\\
		&\quad+\delta_{nm} \delta_{\ell \iota}\big(\theta\beta_n+1-d/2-\mu\big)\big[q_{k, n} (1) q_{j, n} (1)-\delta_{2\mu+d, 2}\, q_{k, n} (0) q_{j, n}(0)\big] \nonumber\\
		&= \delta_{nm} \delta_{\ell \iota}\int_{0}^{1}  \big(P_{k}^{(-1,\beta_{n})}(2 r^{2\theta}-1)\big)'\big( P_{j}^{(-1, \beta_{n})}(2 r^{2\theta}-1) \big)' r^{2\theta\beta_n+1} \mathrm{d} r \nonumber\\
		&\quad+ \delta_{nm} \delta_{\ell \iota}\big(\theta\beta_n+1-d/2-\mu\big)\left[\delta_{k0} \delta_{kj}-\delta_{2\mu+d, 2}\delta_{\theta\beta_n+1-d/2-\mu,0} \right] \nonumber\\
		&= \delta_{nm} \delta_{\ell \iota}\frac{2\theta}{2^{\beta_{n}}}\int_{-1}^{1} \frac{\rm d}{{\rm d}{\rho}} P_{k}^{(-1,\beta_n)}(\rho)\frac{\rm d}{{\rm d}{\rho}}  P_{j}^{(-1,\beta_n)}(\rho) (1+\rho)^{\beta_{n}+1}  \mathrm{d} \rho \nonumber\\
		&\quad+ \delta_{nm} \delta_{\ell \iota}\big(\theta\beta_n+1-d/2-\mu\big)\delta_{k0} \delta_{kj} , \nonumber
	\end{align}
	where in the last equatlity, we use the variable transformation $\rho=r^{2\theta}-1$.
	Using the property 
	\begin{equation}
\frac{\rm d}{{\rm d}x}  P_k^{(-1,\beta_n)}(x)=\frac{k+\beta_n}{2}P_{k-1}^{(0,\beta_n+1)}(x),\label{Jacobipro1}
\end{equation}
We derive   from \eqref{Jacobipro1}  and  \eqref{GGPeq9} that
 \begin{align}
 &	\int_{-1}^{1} \frac{\rm d}{{\rm d}{\rho}} P_{k}^{(-1,\beta_n)}(\rho)\frac{\rm d}{{\rm d}{\rho}} P_{j}^{(-1,\beta_n)}(\rho) (1+\rho)^{\beta_{n}+1}  \mathrm{d} \rho \nonumber\\
 =& \frac{(k+\beta_{n})(j+\beta_{n})}{4}\int_{-1}^{1} P_{k-1}^{(0,\beta_{n}+1)}(\rho) P_{j-1}^{(0, \beta_{n}+1)}(\rho) (1+\rho)^{\beta_{n}+1}  \mathrm{d} \rho\nonumber\\
 =&\frac{2^{\beta_{n}}(k+\beta_{n})^2}{2k+\beta_{n}}\delta_{kj}.
 	\end{align}
 This completes the proof.

\renewcommand{\theequation}{B.\arabic{equation}}
\section{Proof of Proposition \ref{anasolution}}\label{proanasolution}
\noindent
	
Using the Leibniz rule for gradient and divergence, we can reformulate the problem in the spherical-polar coordinates as follows:
	\[\begin{aligned}
		-r^{2\mu}\Big(\partial^2_r+ \frac{2\mu+d-1}{r} \partial_r + \frac{\Delta_0-c}{r^2}\Big)u=-\frac{1}{r^{d-1}}\partial_r(r^{2\mu+d-1} \partial_r u) - \frac{\Delta_0-c}{r^{2-2\mu}}u=\lambda u.
	\end{aligned}\]
	We now represent the unknown eigenfunction $u$ as an expansion of spherical harmonic function,
	\[
	u(\bm{x})=\sum_{n=0}^{\infty} \sum_{\ell=1}^{a_{n}^{d}} u_{\ell}^{n}(r) Y_{\ell}^{n}(\bm{\hat{x}}),
	\]
	and then obtain the radial eigenvalue problem 
	\begin{equation}\label{eigequation}
		-\frac{1}{r^{d-1}}\partial_r\big(r^{2\mu+d-1} \partial_r u_{\ell}^{n}(r)\big) +\frac{n(n+d-2)+c}{r^{2-2\mu}}u_{\ell}^{n}(r)=\lambda u_{\ell}^{n}(r),\quad 1 \leq \ell \leq a_{n}^{d},\;\; n \geq 0.
	\end{equation}
	Let $r=\rho^\frac{1}{1-\mu}$ and set $v_{\ell}^{n}(\rho)=u_{\ell}^{n}(r)$. We can reformulate \eqref{eigequation} as
	\begin{equation}\label{eigequation1}
		-\frac{1}{\rho^{\frac{d}{1-\mu}-1}}\partial_\rho\big[\rho^{\frac{d}{1-\mu}-1}\partial_\rho v_{\ell}^{n}(\rho)\big]+ \frac{n(n+d-2)+c}{(1-\mu)^2\rho^2}v_{\ell}^{n}(\rho)=\frac{\lambda}{(1-\mu)^2}
		v_{\ell}^{n}(\rho),
	\end{equation}
	and  \eqref{eigequation} can be written as
	\[
	\rho^{2} \partial_{\rho}^{2}\big[\rho^{\frac{d }{2(1-\mu)}-1} v_{\ell}^{n}(\rho)\big]+\rho \partial_{\rho}\big[\rho^{\frac{d }{2(1-\mu)}-1} v_{\ell}^{n}(\rho) \big] +\big[\frac{\lambda}{(1-\mu)^2} \rho^{2}-\frac{\beta_n^{2}}{4}\big]\rho^{\frac{d }{2(1-\mu)}-1} v_{\ell}^{n}(\rho)=0.
	\]
	Making the variable transformation $\eta=\frac{\sqrt{\lambda}}{1-\mu}\rho$ and setting $\chi_{\ell}^{n}(\eta)=\rho^{\frac{d }{2(1-\mu)}-1}v_{\ell}^{n}(\rho)$, one obtains
	\[\eta^{2} \partial_{\eta}^{2} \chi_{\ell}^{n}(\eta)+\eta \partial_{\eta} \chi_{\ell}^{n}(\eta)+\Big(\eta^{2}-\frac{\beta_n^{2}}{4}\Big) \chi_{\ell}^{n}(\eta)=0,\]
	which is exactly the Sturm-Liouville equation for the first kind Bessel function and here admits a unique solution $\chi_{\ell}^{n}(\eta)=J_{\frac{\beta_n}{2}}(\eta)$, where $J_{\alpha}(x)$ is the first kind Bessel function which can be expressed as
	\[J_\alpha(x)=\sum_{m=0}^{\infty} \frac{(-1)^m}{m ! \Gamma(m+\alpha+1)}\left(\frac{x}{2}\right)^{2 m+\alpha}.\]
	In return,
	\begin{equation}\label{analysissolution}
		\begin{aligned}
			u_{\ell}^{n}(r) &=v_{\ell}^{n}(\rho)=\rho^{1-\frac{d }{2(1-\mu)}} \chi_{n}(\eta)=\rho^{1-\frac{d }{2(1-\mu)}}J_{\frac{\beta_n}{2}}\Big(\frac{\sqrt{\lambda}}{1-\mu}\rho\Big)=r^{1-\mu-\frac{d}{2}}J_{\frac{\beta_n}{2}}\Big(\frac{\sqrt{\lambda}}{1-\mu}r^{1-\mu}\Big)\\
			&=\sum_{m=0}^{\infty} \frac{(-1)^{m}}{m ! \Gamma\big(m+\frac{\beta_n}{2}+1\big)}\Big(\frac{\sqrt{\lambda} r^{1-\mu}}{2(1-\mu)}\Big)^{2 m+\beta_n/2} r^{1-\mu-\frac{d}{2}} .
		\end{aligned}
	\end{equation}
	Since the homogeneous Dirichlet boundary condition implies $u_{\ell}^{n}(1)=0$, one readily finds that the eigenvalue $\lambda$ of \eqref{eigequation} satisfies
	\[J_{\frac{\beta_n}{2}}\Big(\frac{\sqrt{\lambda}}{1-\mu}\Big)=0.\]
	This completes the derivation.
	
	\renewcommand{\theequation}{C.\arabic{equation}}
	\section{Proof of Lemma \ref{lemmamuntz1}}\label{prolemmamuntz1}		
	\noindent	From \eqref{ballint}, it is easy to see that
	\begin{align}\label{GUPeq11}
		&\big(S_{k,\ell,n,c}^{-1,\mu,1-\mu},  S_{j,\iota,m,c}^{-1,\mu,1-\mu}\big)=\int_{\mathbb{S}^{d-1}} Y_{\ell}^{n}(\bm{\hat{x}}) Y_{\iota}^{m}(\bm{\hat{x}}) \mathrm{d} \sigma(\bm{\hat{x}}) \nonumber\\
		&\quad \times  \int_{0}^{1} P_{k}^{(-1,\beta_n)}(2 r^{2-2\mu}-1) P_{j}^{(-1,\beta_m)}(2 r^{2-2\mu}-1) r^{(1-\mu)(\beta_n+\beta_m)-2\mu+1} \mathrm{d} r \nonumber\\
		&=\delta_{k j}  \delta_{n m}  \int_{0}^{1} P_{k}^{(-1, \beta_n)}(2 r^{2-2\mu}-1) P_{j}^{(-1,\beta_n)}(2 r^{2-2\mu}-1) r^{2(1-\mu)(\beta_n+1)-1} \mathrm{d} r \nonumber\\
		&=\delta_{k j}  \delta_{n m}  \frac{1}{2^{2+\beta_n}(1-\mu)}\int_{-1}^{1} P_{k}^{(-1, \beta_n)}(\rho) P_{j}^{(-1,\beta_n)}(\rho) (1+\rho)^{\beta_n} \mathrm{d} \rho .
	\end{align}
	where in the last inequality, we use the change of variable $\rho=2r^{2-2\mu}-1$. Then, from the  properties of Jacobi polynomials in \eqref{LemEQ2} and \eqref{symEQ1}, we obtain
	\begin{equation}\label{GUPeq50}
		P_{k}^{(\alpha, \beta)}(x)=\frac{k+\alpha+\beta+1}{2 k+\alpha+\beta+1} P_{k}^{(\alpha+1, \beta)}(x)-\frac{k+\beta}{2 k+\alpha+\beta+1} P_{k-1}^{(\alpha+1, \beta)}(x),
	\end{equation}
	which implies
	\begin{equation}\label{GUPeq12}
		P_{k}^{(-1, \beta_n)}= \frac{(k+\beta_n)}{(2k+\beta_n)} P_{k}^{(0, \beta_n)}-\frac{(k+\beta_n)}{(2k+\beta_n)} P_{k-1}^{(0, \beta_n)}.
	\end{equation}
	Then a combination of \eqref{GUPeq11}, \eqref{GUPeq12}, and \eqref{GGPeq9} immediately yields \eqref{muntzmassmatrix1}.

\end{appendix}

\bibliographystyle{siam}
\bibliography{refGUPf}

\begin{thebibliography}{10}

\bibitem{Andrews1999special}
{\sc G.~E. Andrews, R.~Askey, and R.~Roy}, {\em Special functions}, vol.~71,
  Cambridge University Press, Cambridge, 1999.

\bibitem{Atkinson2019Book}
{\sc K.~Atkinson, D.~Chien, and O.~Hansen}, {\em Spectral Methods Using
  Multivariate Polynomials On The Unit Ball}, CRC Press, 2019.

\bibitem{Cagliero2015Explicit}
{\sc L.~Cagliero and T.~H. Koornwinder}, {\em Explicit matrix inverses for
  lower triangular matrices with entries involving {J}acobi polynomials}, J.
  Approx. Theory, 193 (2015), pp.~20--38.

\bibitem{Cao2006Solutions}
{\sc D.~Cao and P.~Han}, {\em Solutions to critical elliptic equations with
  multi-singular inverse square potentials}, J. Differential Equations, 224
  (2006), pp.~332--372.

\bibitem{Cheney1998Introduction}
{\sc E.~W. Cheney}, {\em Introduction to approximation theory}, AMS Chelsea
  Publishing, Providence, RI, 1998.
\newblock Reprint of the second (1982) edition.

\bibitem{Chihara1978an}
{\sc T.~S. Chihara}, {\em An introduction to orthogonal polynomials}, vol.~13,
  Gordon and Breach Science Publishers, New York-London-Paris, 1978.

\bibitem{Dai2013approximation}
{\sc F.~Dai and Y.~Xu}, {\em Approximation theory and harmonic analysis on
  spheres and balls}, Springer, New York, 2013.

\bibitem{Dunkl2014Orthogonal}
{\sc C.~F. Dunkl and Y.~Xu}, {\em Orthogonal polynomials of several variables},
  vol.~155 of Encyclopedia of Mathematics and its Applications, Cambridge
  University Press, Cambridge, second~ed., 2014.

\bibitem{Dyda2017Eigenvalues}
{\sc B.~o. Dyda, A.~Kuznetsov, and M.~Kwa\'{s}nicki}, {\em Eigenvalues of the
  fractional {L}aplace operator in the unit ball}, J. Lond. Math. Soc. (2), 95
  (2017), pp.~500--518.

\bibitem{Dyda2017Fractional}
\leavevmode\vrule height 2pt depth -1.6pt width 23pt, {\em Fractional {L}aplace
  operator and {M}eijer {G}-function}, Constr. Approx., 45 (2017),
  pp.~427--448.

\bibitem{Felli2007Schorodinger}
{\sc V.~Felli, E.~M. Marchini, and S.~Terracini}, {\em On {S}chr\"{o}dinger
  operators with multipolar inverse-square potentials}, J. Funct. Anal., 250
  (2007), pp.~265--316.

\bibitem{Felli2006Elliptic}
{\sc V.~Felli and S.~Terracini}, {\em Elliptic equations with multi-singular
  inverse-square potentials and critical nonlinearity}, Comm. Partial
  Differential Equations, 31 (2006), pp.~469--495.

\bibitem{Kufner1985Weighted}
{\sc A.~Kufner}, {\em Weighted {S}obolev spaces}, A Wiley-Interscience
  Publication, John Wiley \& Sons, Inc., New York, 1985.
\newblock Translated from the Czech.

\bibitem{Li2010Optimal}
{\sc H.~Li and J.~Shen}, {\em Optimal error estimates in {J}acobi-weighted
  {S}obolev spaces for polynomial approximations on the triangle}, Math. Comp.,
  79 (2010), pp.~1621--1646.

\bibitem{Li2014spectral}
{\sc H.~Li and Y.~Xu}, {\em Spectral approximation on the unit ball}, SIAM J.
  Numer. Anal., 52 (2014), pp.~2647--2675.

\bibitem{Li2017efficient}
{\sc H.~Li and Z.~Zhang}, {\em Efficient spectral and spectral element methods
  for eigenvalue problems of {S}chr\"{o}dinger equations with an inverse square
  potential}, SIAM J. Sci. Comput., 39 (2017), pp.~A114--A140.

\bibitem{Ma2018Efficient}
{\sc S.~Ma, H.~Li, and Z.~Zhang}, {\em Efficient spectral methods for some
  singular eigenvalue problems}, J. Sci. Comput., 77 (2018), pp.~657--688.

\bibitem{Olver2020Orthogonal}
{\sc S.~Olver and Y.~Xu}, {\em Orthogonal polynomials in and on a quadratic
  surface of revolution}, Math. Comp., 89 (2020), pp.~2847--2865.

\bibitem{Sheng2021Nontensorial}
{\sc C.~Sheng, S.~Ma, H.~Li, L.-L. Wang, and L.~Jia}, {\em Nontensorial
  generalised {H}ermite spectral methods for {PDE}s with fractional {L}aplacian
  and {S}chr\"{o}dinger operators}, ESAIM Math. Model. Numer. Anal., 55 (2021),
  pp.~2141--2168.

\bibitem{Szego1975orthogonal}
{\sc G.~Szeg\H{o}}, {\em Orthogonal polynomials}, vol.~XXIII, American
  Mathematical Society, Providence, R.I., fourth~ed., 1975.

\bibitem{Zhang2020ACHA}
{\sc J.~Zhang, H.~Li, L.-L. Wang, and Z.~Zhang}, {\em Ball prolate spheroidal
  wave functions in arbitrary dimensions}, Appl. Comput. Harmon. Anal., 48
  (2020), pp.~539--569.

\end{thebibliography}

\end{document}